\documentclass[]{amsart}

\usepackage{amsmath}
\usepackage{amsthm}
\usepackage{amssymb}
\usepackage{latexsym}               
\usepackage{amsgen}
\usepackage[mathscr]{eucal}

\usepackage{array,arydshln}

\numberwithin{equation}{section}

\newcommand{\R}{\mathbb{R}}

\usepackage{color}

\newcommand{\dis}{\displaystyle}

\newcommand{\lag}{\langle}
\newcommand{\rag}{\rangle}

\newcommand{\al}{\alpha}
\newcommand{\be}{\beta}

\newcommand{\la}{\lambda}

\newcommand{\ga}{\gamma}

\newcommand{\pa}{\partial}

\newtheorem{lem}{Lemma}[section]
\newtheorem{thm}[lem]{Theorem}
\newtheorem{pro}[lem]{Proposition}

\newtheorem{defn}[lem]{Definition}

\makeatletter
\newif\if@borderstar
\def\bordermatrix{\@ifnextchar*{%
 \@borderstartrue\@bordermatrix@i}{\@borderstarfalse\@bordermatrix@i*}%
}
\def\@bordermatrix@i*{\@ifnextchar[{\@bordermatrix@ii}{\@bordermatrix@ii[()]}}
\def\@bordermatrix@ii[#1]#2{%
\begingroup
 \m@th\@tempdima8.75\p@\setbox\z@\vbox{%
 \def\cr{\crcr\noalign{\kern 2\p@\global\let\cr\endline }}%
 \ialign {$##$\hfil\kern 2\p@\kern\@tempdima & \thinspace %
  \hfil $##$\hfil && \quad\hfil $##$\hfil\crcr\omit\strut %
  \hfil\crcr\noalign{\kern -\baselineskip}#2\crcr\omit %
  \strut\cr}}%
 \setbox\tw@\vbox{\unvcopy\z@\global\setbox\@ne\lastbox}%
 \setbox\tw@\hbox{\unhbox\@ne\unskip\global\setbox\@ne\lastbox}%
 \setbox\tw@\hbox{%
  $\kern\wd\@ne\kern -\@tempdima\left\@firstoftwo#1%
  \if@borderstar\kern 2pt\else\kern -\wd\@ne\fi%
 \global\setbox\@ne\vbox{\box\@ne\if@borderstar\else\kern 2\p@\fi}%
 \vcenter{\if@borderstar\else\kern -\ht\@ne\fi%
  \unvbox\z@\kern -\if@borderstar2\fi\baselineskip}%
 \if@borderstar\kern-4.5\@tempdima\kern2\p@\else\,\fi\right\@secondoftwo#1 $%
 }\null \;\vbox{\kern\ht\@ne\box\tw@}%
\endgroup
}
\makeatother

\begin{document}

\title[Decay structure of two hyperbolic relaxation models]
{Decay structure of two hyperbolic relaxation models with regularity-loss}

\author[Y. Ueda]{Yoshihiro Ueda}
\address{(YU) 
Faculty of Maritime Sciences, Kobe University, Kobe 658-0022, Japan}
\email{ueda@maritime.kobe-u.ac.jp}


\author[R.-J. Duan]{Renjun Duan}
\address{(RJD)
Department of Mathematics, The Chinese University of Hong Kong, Shatin, Hong Kong}
\email{rjduan@math.cuhk.edu.hk}

\author[S. Kawashima]{Shuichi Kawashima}
\address{(SK)
Faculty of Mathematics, Kyushu University, Fukuoka 819-0395, Japan}
\email{kawashim@math.kyushu-u.ac.jp}



\def\cal#1{{\fam2#1}}


%
\begin{abstract}
The paper aims at investigating two types of decay structure for 
linear symmetric hyperbolic systems with non-symmetric
 relaxation. Precisely, the system is of the type $(p,q)$ if the real
 part of all eigenvalues admits an upper bound
 $-c|\xi|^{2p}/(1+|\xi|^2)^{q}$, where $c$ is a generic positive
 constant and $\xi$ is the frequency variable, and the system enjoys the
 regularity-loss property if $p<q$. It is well known that the standard
 type $(1,1)$ can be assured by  the classical Kawashima-Shizuta
 condition.  A new structural condition was introduced in  \cite{UDK} to
 analyze the regularity-loss type $(1,2)$ system with non-symmetric relaxation. In the paper, we construct two more complex models of the  regularity-loss type corresponding to $p=m-3$, $q=m-2$ and $p=(3m-10)/2$, $q=2(m-3)$, respectively, where $m$ denotes phase dimensions. The proof is based on the delicate Fourier energy method as well as the suitable linear combination of series of energy inequalities. Due to arbitrary higher dimensions, it is not obvious to capture the energy dissipation rate with respect to the degenerate components. Thus, for each model, the analysis always starts from the case of low phase dimensions in order to understand the basic dissipative structure in the general case, and in the mean time, we also give the explicit construction of the compensating symmetric matrix $K$ and skew-symmetric matrix $S$.

\end{abstract}

\maketitle


Keywords: Decay structure, Regularity-loss, Symmetric hyperbolic system, 
Energy method


MSC 2010: 35B35, 
35B40, 
35L40 

\tableofcontents
\thispagestyle{empty}

\section{Introduction}

In the paper, we consider the Cauchy problem on
the following linear symmetric hyperbolic system with relaxation (cf.~\cite{Da}):
\begin{gather}
u_t + A_m u_{x} + L_m u = 0 \label{sys1}
\end{gather}
%
%
%
with
\begin{equation}\label{ID}
    u|_{t=0}=u_0.
\end{equation}
Here
$u=u(t,x) = (u_1, \cdots, u_m)^T(t,x) \in \R^m$ over  $t >0$, $x \in \mathbb{R}$ is an unknown function,
$u_0=u_0(x)\in \R^m$ over $x\in \R$ is a given function,
and $A_m$ and $L_m$ are $m\times m$ real constant matrices. In general we assume $A_m$ is symmetric and $L_m$ is degenerately dissipative in the sense of $1\leq \dim\, (\ker L_m)\leq m-1$.
As pointed out in \cite{UDK}, for a general linear degenerately dissipative system it is interesting to study
its decay structure under additional conditions on the coefficient matrices
and further investigate the corresponding time-decay property of solutions to the Cauchy problem at the linear level. 
The purpose of the paper is to present two concrete models of $A_m$ and $L_m$,
which do not satisfy the dissipative condition in \cite{UDK}, 
to derive the decay structure of the corresponding linear systems. 
We remark that the similar issue has been extensively investigated in Villani \cite{Vi} for an infinite-dimensional dynamical system, for instance, in the content of kinetic theory.  

In what follows let us explain the motivation of dealing with the problem considered here. More generally one may consider the system in multidimensional space $\R^n$:
\begin{equation}
A^0_m u_t + \sum_{j=1}^n A^j_m u_{x_j} + L_mu = 0, \label{In.sys1}
\end{equation}
where $u=u(t,x) \in \mathbb{R}^m$ over  $t \geq 0$, $x \in \mathbb{R}^n$.
When the degenerate relaxation matrix $L_m$ is symmetric,
Umeda-Kawashima-Shizuta \cite{UKS84} proved the large-time asymptotic stability
of solutions for a class of equations of hyperbolic-parabolic type with applications
to both electro-magneto-fluid dynamics and magnetohydrodynamics.
The key idea in \cite{UKS84} and the later generalized work \cite{SK85}
that first introduced the so-called Kawashima-Shizuta (KS) condition is
to construct the compensating matrix to capture the dissipation
of systems over the degenerate kernel space of $L_m$. The typical
feature of the time-decay property of solutions established in those work
is that the high frequency part decays exponentially while the low frequency part
decays polynomially with the same rate as  the heat kernel. 
To precisely state these results, 
we apply Fourier transform to \eqref{In.sys1} (or \eqref{sys1}).
Then we can obtain
\begin{equation}
A^0_m \hat{u}_t + i |\xi| A_m(\omega) \hat{u} + L_m \hat{u} = 0, \label{Fsys1}
\end{equation}
where
$\xi\in \R^n$ denote the Fourier variable of $x\in \R^n$,
$\omega=\xi/|\xi|\in S^{n-1}$, and $A_m(\omega) := \sum_{j=1}^n A^j_m\omega_j$. 
Moreover we prepare some notations. 
Given a real matrix $X$, we use
$X^{\rm sy}$ and $X^{\rm asy}$ to denote the symmetric and skew-symmetric parts of $X$,
respectively, namely, $X^{\rm sy}=(X+X^T)/2$ and $X^{\rm asy}=(X-X^T)/2$. 
%
Then the decay result in \cite{UKS84,SK85} is stated as follows.

\begin{pro}
[Decay property of the standard type (\cite{UKS84,SK85})]\label{pro1}
%
Consider \eqref{In.sys1} with the following condition:
\begin{description}
  \item[Condition (A)$_0$] $A^0_m$ is real symmetric and positive definite, $A^j_m$ for each $1\leq j\leq n$ is
real symmetric, and $L_m$ is real symmetric and nonnegative definite with
the nontrivial kernel. 
\end{description}
For this problem, assume that the following condition hold:
\begin{description}
\item[Condition (K)]
There is a real compensating matrix
$K(\omega) \in C^{\infty}(S^{n-1})$ with the properties:
$K(-\omega) = - K(\omega)$, $(K(\omega)A^0_m)^T = - K(\omega)A^0_m$ and
\begin{equation*}
[K(\omega)A_m(\omega)]^{\rm sy} >0 \quad \text{on} \quad \ker L_m
\end{equation*}
for each $\omega\in S^{n-1}$. 
\end{description}
Then the Fourier image $\hat u$ of the solution $u$ to the equation \eqref{In.sys1} with initial data $u(0,x)=u_0(x)$ satisfies the pointwise estimate:
\begin{equation}\label{std-point}
|\hat u(t,\xi)| \le Ce^{-c \la(\xi)t}|\hat u_0(\xi)|,
\end{equation}
where $\la(\xi) := |\xi|^2/(1+|\xi|^2)$.
Furthermore, let $s \ge 0$ be an integer and suppose that the initial data
$u_0$ belong to $H^s \cap L^1$.
Then the solution $u$ satisfies the decay estimate: 
\begin{equation}\label{std-decay}
\|\partial_x^{k} u(t)\|_{L^2} \le C(1+t)^{-n/4-k/2}\|u_0\|_{L^1}
	+ Ce^{-ct}\|\partial_x^{k} u_0\|_{L^2}
\end{equation}
for $k\le s$. Here $C$ and $c$ are positive constants.
\end{pro}

Under the conditions {\rm (A$)_0$} and {\rm (K)}, we can construct the following
energy inequality:
 \begin{equation*}
 \frac{d}{dt} E + c D \le 0,
 \end{equation*}
where
 \begin{equation}\label{1K}
E = \langle A_m^0 \hat{u},\hat{u} \rangle 
- \frac{\alpha|\xi|}{1+|\xi|^2}\delta \langle iK(\omega)A_m^0 \hat{u},\hat{u}
\rangle, \quad
D = \frac{|\xi|^2}{1+|\xi|^2}|\hat{u}|^2 + |(I-P)\hat{u}|^2,
 \end{equation}
$\alpha$ and $\delta$ are suitably small constants, 
and $P$ denotes the orthogonal projection onto $\ker L_m$.

\smallskip

For the nonlinear system, the global existence of small-amplitude classical solutions  was proved by Hanouzet-Natalini \cite{HN} in one space dimension  and  by Yong \cite{Yo} in several space dimensions, provided that the system is strictly entropy dissipative and satisfies the KS condition.
And later on, the large time behavior of solutions was obtained by 
Bianchini-Hanouzet-Natalini \cite{BHN} and Kawashima-Yong \cite{KY2}
basing on the analysis of the Green function of the 
linearized problem. Those results show that solutions to such nonlinear system 
will not develop singularities (e.g., shock waves) in finite time for small 
smooth initial perturbations, cf.~\cite{Da,Liu}. 
Notice that the $L^2$-stability of a constant equilibrium state in a one-dimensional 
system of dissipative hyperbolic balance laws endowed with a convex entropy was 
also studied by Ruggeri-Serre \cite{RS}. 
Moreover, it would be an interesting and important topic to study the relaxation 
limit of general hyperbolic conservation laws with relaxations, see \cite{CLL,KY}  
and reference therein.

Recently it has been found that there exist physical systems which violate the 
KS condition but still have some kind of time-decay properties. For instance, for 
the dissipative Timoshenko system \cite{IHK08,IK08} and the Euler-Maxwell system
\cite{D,USK,UK}, the linearized relaxation matrix $L_m$  has
a nonzero skew-symmetric part while
it was still proved that solutions decay in time in some different way. 
Besides those, there are two related works dealing with general
partially dissipative hyperbolic systems with zero-order source when the
KS condition is not satisfied. Beauchard-Zuazua \cite{BZ11} first
observed the equivalence of the KS condition with the Kalman rank
condition in the context of the control theory. 
They extended the previous
analysis to some other situations beyond the KS condition, and
established the explicit estimate on the solution semigroup in terms of
the frequency variable and also the global existence of near-equilibrium
classical solutions for some nonlinear balance laws without the KS
condition. 
In the mean time, Mascia-Natalini \cite{MN} also made a general study of 
the same topic for a class of systems without the KS condition. 
The typical situation considered in \cite{MN} is that the non-dissipative 
components are linearly degenerate which indeed does not hold under the KS condition (see also \cite{KO}). 
Notice that both in \cite{BZ11} and \cite{MN}, the rate of convergence of solutions 
to the equilibrium states for the nonlinear Cauchy problem is still left unknown.

In \cite{UDK},  the same authors of this paper introduced a new structural condition
which is a generalization of  the KS condition, 
and also analyzed the corresponding weak dissipative
structure called the regularity-loss type for general systems with
non-symmetric relaxation which includes the Timoshenko system and the Euler-Maxwell system as two concrete examples. Precisely, one has the following result.

\begin{pro}
[Decay property of the regularity-loss type (\cite{UDK})]\label{pro2}
Consider \eqref{In.sys1} with the condition:
\begin{description}
  \item[Condition (A)] $A^0_m$ is real symmetric and positive definite, $A^j_m$  for each $1\leq j\leq n$ is
real symmetric, while $L_m$ is not necessarily real symmetric but is nonnegative
definite with the nontrivial kernel.
\end{description}
For this problem, assume the previous condition {\rm (K)} and the following condition hold:
\begin{description}
\item[Condition (S)]
There is a real  matrix $S$ such that
$(SA^0_m)^T = SA^0_m$, and
\begin{gather*}
[SL_m]^{\rm sy}+[L_m]^{\rm sy}\geq 0 \ \ {\rm on} \ \ {\mathbb C}^m, \quad
\ker\big([SL_m]^{\rm sy}+[L_m]^{\rm sy}\big) = 
\ker\,L_m,
\end{gather*}
and moreover, for each $\omega\in S^{n-1}$, 
\begin{gather}
i[SA_m(\omega)]^{\rm asy}\geq 0 \quad \text{on} \quad \ker\,[L_m]^{\rm sy}. \label{assump-3}
\end{gather}
\end{description}
Then the Fourier image $\hat u$ of the solution $u$ to the equation \eqref{In.sys1} with initial data $u(0,x)=u_0(x)$ satisfies the pointwise estimate:
\begin{equation}\label{loss-point}
|\hat u(t,\xi)| \le Ce^{-c \la(\xi)t}|\hat u_0(\xi)|,
\end{equation}
where $\la(\xi) := |\xi|^2/(1+|\xi|^2)^2$.
Moreover, let $s \ge 0$ be an integer and suppose that the initial data
$u_0$ belong to $H^s \cap L^1$.
Then the solution $u$ satisfies the decay estimate: 
\begin{equation}\label{loss-decay}
\|\partial_x^{k} u(t)\|_{L^2} \le C(1+t)^{-n/4-k/2}\|u_0\|_{L^1}
	+ C(1+t)^{-\ell/2}\|\partial_x^{k+\ell} u_0\|_{L^2}
\end{equation}
for $k+\ell \le s$. Here $C$ and $c$ are positive constants.
\end{pro}

Observe that $\la(\xi)$ in \eqref{loss-point} behaves as $|\xi|^2$ as $|\xi|\to 0$ but behaves as $1/|\xi|^2$ as $|\xi|\to \infty$. Thus  those estimates \eqref{loss-point} and \eqref{loss-decay} are weaker than
\eqref{std-point} and \eqref{std-decay}, respectively.
In particular, the decay estimate \eqref{loss-point}  is said to be of the regularity-loss type. 
Similar decay properties of the regularity-loss type have been recently
observed for several interesting systems. We refer the reader to
\cite{IHK08,IK08,LK4} (cf. \cite{ABMR,MRR}) for the dissipative Timoshenko system,
\cite{D,USK,UK} for the Euler-Maxwell system,
\cite{HK,KK} for a hyperbolic-elliptic system in radiation gas dynamics,
\cite{LK1,LK2,LK3,LC,SK} for a dissipative plate equation, and
\cite{Du-1VMB,DS-VMB} for various kinetic-fluid models.

In fact, one can show that Proposition \ref{pro1} can be regarded as a corollary of 
Proposition  \ref{pro2} after replacing \eqref{assump-3} in condition {\rm (S)} 
by a stronger condition:
\begin{equation*}
i[SA_m(\omega)]^{\rm asy} \ge 0 \quad \text{on} \quad \mathbb{C}^m.
\end{equation*}
 for each $\omega\in S^{n-1}$. 
 The key point for the proof of \eqref{loss-point} 
 is to derive the matrices $S$ and $K(\omega)$ such that the 
 coercive estimate:
 \begin{equation}\label{core}
 \delta [K(\omega)A_m(\omega)]^{\rm sy} + [SL_m]^{\rm sy} +  [L_m]^{\rm sy} > 0 \quad \text{on} \quad \mathbb{C}^m
 \end{equation}
 holds true for suitably small $\delta>0$.
 Indeed, under the conditions (A), (S) and (K), the estimate
 \eqref{core} is satisfied.
Then, using \eqref{core}, we get the following energy equality
 \begin{equation}\label{known-energy}
 \frac{d}{dt} E + c D \le 0,
 \end{equation}
where
 \begin{equation}\label{1K1S}
\begin{split}
E &= \langle A_m^0 \hat{u},\hat{u} \rangle 
+ \frac{\alpha_1}{1+|\xi|^2}\Big(\langle SA_m^0 \hat{u},\hat{u}
 \rangle  -
 \frac{\alpha_2 |\xi|}{1+|\xi|^2} \delta \langle iK(\omega)A_m^0 \hat{u},\hat{u} \rangle \Big),\\
D &= \frac{|\xi|^2}{(1+|\xi|^2)^2}|\hat{u}|^2 +
 \frac{1}{1+|\xi|^2}|(I-P)\hat{u}|^2 + |(I-P_1)\hat{u}|^2,
\end{split}
 \end{equation}
$\alpha_1$ and $\alpha_2$ are suitably small constants,
and $P$ and $P_1$ denote the orthogonal projections onto $\ker L_m$
and  $\ker \, [L_m]^{\rm sy}$.
 Interested readers may refer to \cite{UDK} for more  details of this issue and also for the construction of $S$ and $K(\omega)$ for the Timoshenko system and the Euler-Maxwell system. Therefore, those conditions in Proposition  \ref{pro2} are generalizations of  the classical KS  conditions. 
We finally remark that it should be interesting to further investigate the nonlinear stability of constant equilibrium states of the system of the regularity-loss type under the structural condition postulated in Proposition \ref{pro2}.
 
Inspired by the previous work \cite{UDK}, the goal of the paper is to construct much more complex models \eqref{sys1} with given $A_m$ and $L_m$ such that they enjoy  some new dissipative structure of the regularity-loss type. Here we recall  a notion of the {\it uniform dissipativity} of
the system \eqref{sys1} introduced in \cite{UDK}. Consider the eigenvalue problem for the system
\eqref{sys1}:
\begin{equation*}
(\eta A^0_m + i\xi A_m + L_m) \phi = 0,
\end{equation*}
where $\eta \in \mathbb{C}$ and $\phi \in \mathbb{C}^m$.
The corresponding characteristic equation is given by
\begin{equation}\label{ce}
{\rm det}(\eta A^0_m + i\xi A_m + L_m) = 0.
\end{equation}
The solution $\eta = \eta(i \xi)$ of \eqref{ce} is called the
eigenvalue of the system \eqref{sys1}.

\begin{defn}
The system \eqref{sys1} is called uniformly dissipative
of the type $(p,q)$ if the eigenvalue $\eta=\eta(i\xi)$ satisfies
$$
\Re\,\eta(i\xi) \leq -c|\xi|^{2p}/(1+|\xi|^2)^q
$$
for all $\xi\in \mathbb{R}^n$, where $c$ is a positive constant and $(p,q)$ is a pair
of positive integers.
\end{defn}

Note that as proved in \cite[Theorem 4.2]{UDK}, one has 
$\Re\,\eta(i\xi)\leq -c\la(\xi)$ whenever the pointwise estimates in the form of 
\eqref{std-point} or \eqref{loss-point} hold true. 
Therefore, we can determine the type $(p,q)$ for a uniformly dissipative system 
\eqref{sys1} in terms of the function $\la(\xi)$ obtained from the pointwise 
estimate on $\hat u(t,\xi)$:
\begin{equation}
\label{D.key}
|\hat u(t,\xi)| \le Ce^{-c \la(\xi)t}|\hat u_0(\xi)|.
\end{equation}
For example, under the assumption in Propositions \ref{pro1} or \ref{pro2},
the system \eqref{sys1} is uniformly dissipative of the type $(1,1)$ or $(1,2)$, respectively. 
Notice that the regularity-loss type corresponds to 
the situation when $p$ is strictly less than $q$, i.e., $p<q$.

Historically, Shizuta-Kawashima \cite{SK} showed that, under the condition 
{\rm(A$)_0$}, the strict dissipativity $\Re \, \eta(i\xi) < 0$
for $\xi \neq 0$ is equivalent to the uniform dissipativity of the type $(1,1)$. 
Moreover, they showed the pointwise estimate \eqref{std-point} by using 
only one compensating skew-symmetric matrix $K(\omega)$ (see \eqref{1K}). 
On the other hand, the authors formulated in \cite{UDK} a class of systems whose 
dissipativity is of the type $(1,2)$ and got Proposition \ref{pro2}. 
Notice that, in this cace, we need to use one compensating symmetric matrix $S$ 
and one compensating skew-symmetric matrix $K(\omega)$ 
to get the desired pointwise estimate \eqref{loss-point} 
(see \eqref{1K1S}). 
We note that the dissipative Timoshenko system and the Euler-Maxwell system 
studied in \cite{IHK08} and \cite{UK}, respectively, are included in the 
class of systems with the type $(1,2)$ which was formulated in \cite{UDK}. 
However, to get the optimal dissipative estimate for these two examples, 
we need to use one $S$ and two different $K(\omega)$ (see \cite{MK1,UK}). 
%

On the other hand, more complicated concrete models are found in these years. 
Indeed, Mori-Kawashima \cite{MK2} considered the Timoshenko-Cattaneo 
system with heat conduction and showed that its dissipativity is of the 
type $(2,3)$. 
Moreover, they proved the optimal dissipative estimate by using four 
different $S$ and four different $K(\omega)$. 
This means that Proposition \ref{pro2} and the class formulated in 
\cite{UDK} is not enough to analyze the dissipativity of general 
systems \eqref{In.sys1}, and we have to study other concrete models.

In this paper, we will present a study of two concrete models of the system 
\eqref{sys1} related to the above general issue. For the Model I, one has
\begin{equation*}
p=m-3,\quad q=m-2,
\end{equation*}
see \eqref{point1} in Theorem \ref{thm1}.  While for the Model II, we let 
$m$ be even and one has
\begin{equation*}
p=\frac{1}{2}(3m-10),\quad q=2(m-3), 
\end{equation*}
see  \eqref{point2} in Theorem \ref{thm2}.  
In both cases we see $p<q$ and hence two models that we consider are of the 
regularity-loss type.
Compared with the energy inequality \eqref{known-energy}, the energy
inequalities of the Model I and II are much more complicated.
More precisely, to control the dissipation term, we must employ a lot of
compensating symmetric matrices and skew-symmetric matrices 
whose numbers 
depend on the dimension $m$ of the coefficient matrices.
Therefore we can not apply Proposition \ref{pro2} to the Model I and II,
and need direct calculations (see in Section 2 and 3).

The proof of  the estimate in the form of \eqref{D.key} is based on the Fourier energy method, and  in the mean time we also give the explicit construction of matrices $S$ and $K$ as used in Proposition \ref{pro2}.
 As seen later on, a series of energy estimates are derived and their appropriate linear combination leads to a Lyapunov-type inequality of the time-frequency functional equivalent with $|\hat u(t,\xi)|^2$, which hence implies \eqref{D.key}. The most difficult point is that it is priorly unclear to justify whether one choice of $(p,q)$ is optimal; see more discussions in Section 4.1. For that purpose, we also present an alternative approach to find out the value of $(p,q)$ for  both Model I and Model II, and the detailed strategy of the approach is to be given later on.


The rest of the paper is organized  as follows. In Section 2 and Section 3, we study Model I and Model II, respectively. In each section, for the given model, we first state the main results on the dissipative structure and the decay property of the system \eqref{sys1}, give the proof by the energy method in the case $m=6$ which indeed corresponds to some existing physical models, show the proof in the general case $m\geq 6$ still using the energy method, and finally give the explicit construction of matrices $S$ and $K$. 
The matrices $S$ and $K$ constructed in Subsection 2.3 and 3.3 have a very important role in obtaining the coercive estimate similar to \eqref{core}. 
Consequently, by employing these matrices, we can derive the desired pointwise estimates through \eqref{1estSK} and \eqref{2estSK} to be verified later on.
In the last Section 4, we provide another approach to justify the dissipative structure of the system \eqref{sys1}.

\medskip
\noindent{\it Notations.}\ \ 
%
For a nonnegative integer $k$, we denote by $\partial^k_x$
the totality of all the $k$-th order derivatives with respect
to $x = (x_1, \cdots ,x_n)$.
Let $1\leq p\leq\infty$. Then $L^p=L^p({\mathbb R}^n)$ denotes
the usual Lebesgue space over ${\mathbb R}^n$ with the norm
$\|\cdot\|_{L^p}$. For a nonnegative integer $s$,
$H^{s}=H^{s}({\mathbb R}^n)$ denotes
the $s$-th order Sobolev space over ${\mathbb R}^n$ in the $L^2$
sense, equipped with the norm $\|\cdot\|_{H^{s}}$.
We note that 
$L^2=H^0$.
%
Finally, in this paper, we use $C$ or $c$ to denote various positive
constants without confusion.


\section{Model I}

\subsection{Main result I}

In this section, we consider the Cauchy problem \eqref{sys1}, \eqref{ID} with coefficient matrices given by
\begin{equation}\label{Tmat}
\begin{split}
A_m &=
\left(
\begin{array}{cccc:cccc}
 {0} & {1} & {0} & {0} &    &     &     &      \\
 {1} & {0} & {0} & {0} &    &     &     &           \\
 {0} & {0} & {0} & {a_4} & {0} &      &       \mbox{\smash{\huge\textit{O}}}   &       \\
 {0} & {0} & {a_4} & {0} & {a_5} & {0} &           &      \\
 \hdashline
       &      & {0} & {a_5} & {0} & {a_6} &              &      \\
       &       &       & {0} & {a_6} & \ddots &           &      \\
       &  \mbox{\smash{\huge\textit{O}}}  &   &    &   &  &  & {a_{m}} \\
       &    &    &    &    &      &   {a_{m}} & {0} \\
\end{array}
\right),  \ \
L_m =
\left(
\begin{array}{cccc:cccc}
 {0} & {0} & {0} & {1}  &    &    &     &      \\
 {0} & {0} & {0} & {0}  &    &    &     &       \\
 {0} & {0} & {0} & {0}  &    &    &    \mbox{\smash{\huge\textit{O}}}  &   \\
 {-1} & {0} & {0} & {0} &    &    &    &    \\
 \hdashline
       &      &      &      & {0} &  &    &     \\
       &      &      &      &     & \ddots &    &   \\
       &     \mbox{\smash{\huge\textit{O}}}    &     &     &    &     &{0} &     \\
       &      &      &       &     &    &     & {\gamma} \\
 \end{array}
\right),
\end{split}
\end{equation}
where integer $m\geq 6$ is even, $\gamma > 0$, and all elements $a_j$ $(4\leq j\leq m)$ are nonzero.
We note that the system \eqref{sys1}, \eqref{Tmat} with $m=6$ is the Timoshenko system with the heat conduction via Cattaneo law (cf.~\cite{FSR,SAR}).
For this problem, we can derive the following decay structure. 

\begin{thm}
\label{thm1}
%
The Fourier image $\hat u$ of the solution $u$ to the Cauchy problem
\eqref{sys1}-\eqref{ID} with \eqref{Tmat} satisfies the pointwise estimate:
\begin{equation}\label{point1}
|\hat u(t,\xi)| \le Ce^{-c \lambda(\xi)t}|\hat u_0(\xi)|,
\end{equation}
where $\lambda(\xi) := \xi^{2(m-3)}/(1+\xi^2)^{m-2}$.
Furthermore, let $s \ge 0$ be an integer and suppose that the initial data
$u_0$ belong to $H^s \cap L^1$.
Then the solution $u$ satisfies the decay estimate: 
\begin{equation}\label{decay1}
\|\partial_x^{k} u(t)\|_{L^2} \le C(1+t)^{-\frac{1}{2(m-3)}(\frac{1}{2}+k)}\|u_0\|_{L^1} + C(1+t)^{-\frac{\ell}{2}}\|\partial_x^{k+\ell} u_0\|_{L^2}
\end{equation}
for $k+\ell \le s$. Here $C$ and $c$ are positive constants.
\end{thm}
%

We remark that the estimates \eqref{point1} and \eqref{decay1} with $m=6$ is not optimal.
Indeed, Mori-Kawashima \cite{MK2} showed more sharp estimates.

The decay estimate \eqref{decay1} is derived by the pointwise estimate \eqref{point1} in Fourier space immediately.
Thus readers may refer to the same authors' paper \cite{UDK} (see also \cite{D-AA}) and we omit the proof of \eqref{decay1} for brevity. In order to make the proof more precise, we first consider the special case $m=6$ in Section 2.2, and then generalize it to the case $m\geq 6$ in Section 2.3. The proof of \eqref{point1} is given in the following two subsections.


\subsection{Energy method in the case $m=6$}
In this subsection we first consider the case $m=6$.
In such case, the system \eqref{sys1} with \eqref{Tmat} is described as 
\begin{equation}\label{equations6}
\begin{split}
&\pa_t \hat u_1 + i \xi \hat u_2 + \hat u_4 = 0, \\ 
&\pa_t \hat u_2 + i \xi \hat u_1 = 0, \\
&\pa_t \hat u_3 + i \xi a_4 \hat u_4 = 0, \\
&\pa_t \hat u_4 + i \xi (a_4 \hat u_3 + a_5 \hat u_5) - \hat u_1 = 0, \\
&\pa_t \hat u_5 + i \xi (a_{5} \hat u_{4} + a_{6} \hat u_{6}) = 0, \\
&\pa_t \hat u_6 + i \xi a_{6} \hat u_{5} + \gamma \hat u_6 = 0. 
\end{split}
\end{equation}
For this system we are going to apply the energy method to derive Theorem \ref{thm1} in the case $m=6$. The proof is organized by the following three steps.

\medskip

\noindent
{\bf Step 1.}\ \
We first derive the basic energy equality for the system \eqref{equations6} in the Fourier space.
We multiply the all equations of \eqref{equations6} by 
$\bar{\hat{u}} = (\bar{\hat{u}}_1, \bar{\hat{u}}_2,\bar{\hat{u}}_3,\bar{\hat{u}}_4, \bar{\hat{u}}_5,\bar{\hat{u}}_6)^T$, respectively,
and combine the resultant equations.
Then we obtain
\begin{equation*}
\sum_{j=1}^6 \bar{\hat{u}}_j \pa_t    \hat u_j  
+ 2  i \xi  \Re  (\hat u_1 \bar{\hat{u}}_{2}) + 2 i \xi  \sum_{j=3}^5 a_{j+1} \Re  (\hat u_j \bar{\hat{u}}_{j+1}) 
+ 2 i {\rm Im} (\hat{u}_4 \bar{\hat{u}}_1)+ \gamma |\hat u_6|^2 = 0.
\end{equation*}
Thus, taking the real part for the above equality, we arrive at the basic energy equality
\begin{equation}\label{eq6}
\frac{1}{2} \pa_t    |\hat u|^2 + \gamma |\hat u_6|^2 = 0.
\end{equation}
Here we use the simple relation $\partial_t(\hat{u}_j^2) = 2 \Re (\bar{\hat{u}}_j \partial_t \hat{u}_j)$
for any $j$.
Next we create the dissipation terms. 

\medskip

\noindent
{\bf Step 2.}\ \
We first construct the dissipation for $\hat{u}_1$.
We multiply the first and fourth equations in \eqref{equations6} 
by $- \bar{\hat{u}}_{4}$ and $ - \bar{\hat{u}}_{1}$, respectively.
Then, combining the resultant equations and taking the real part, we have
\begin{equation}\label{dissipation6-a}
- \pa_t  \Re (\hat u_{1} \bar{\hat{u}}_4) +   |\hat u_{1}|^2 - |\hat u_{4}|^2  \\
- \xi \, \Re (i \hat{u}_{2} \bar{\hat{u}}_4) + a_{4} \xi \, \Re (i \hat{u}_{1} \bar{\hat{u}}_{3}) 
+ a_{5} \xi \, \Re (i \hat{u}_{1} \bar{\hat{u}}_{5}) =  0.
\end{equation}
On the other hand,
we multiply the second and third equations in \eqref{equations6} 
by $- a_4 \bar{\hat{u}}_{3}$ and $ - a_4 \bar{\hat{u}}_{2}$, respectively.
Then, combining the resultant equations and taking the real part, we have
\begin{equation*}
- a_4 \pa_t  \Re (\hat u_{2} \bar{\hat{u}}_3)
- a_{4} \xi \, \Re (i \hat{u}_{1} \bar{\hat{u}}_{3}) 
+ a_{4}^2 \xi \, \Re (i \hat{u}_{2} \bar{\hat{u}}_{4}) =  0.
\end{equation*}
Therefore, combining the above two equalities, we obtain
\begin{multline}\label{dissipation6-4}
- \pa_t  \Re (\hat u_{1} \bar{\hat{u}}_4+ a_4 \hat u_{2} \bar{\hat{u}}_3) +   |\hat u_{1}|^2 - |\hat u_{4}|^2  \\
+ (a_4^2 - 1) \xi \, \Re (i \hat{u}_{2} \bar{\hat{u}}_4) + a_{5} \xi \, \Re (i \hat{u}_{1} \bar{\hat{u}}_{5}) =  0.
\end{multline}
Furthermore,
we multiply the second equation and fifth equation in \eqref{equations6} 
by  $- \bar{\hat{u}}_{5}$ and $ - \bar{\hat{u}}_{2}$, respectively.
Then, combining the resultant equations and taking the real part, we have
\begin{equation}\label{dissipation6-5}
- \pa_t  \Re (\hat u_{2} \bar{\hat{u}}_5) 
- \xi \, \Re (i \hat{u}_{1} \bar{\hat{u}}_5) + a_{5} \xi \, \Re (i \hat{u}_{2} \bar{\hat{u}}_{4}) 
+ a_{6} \xi \, \Re (i \hat{u}_{2} \bar{\hat{u}}_{6}) =  0.
\end{equation}
Finally, multiplying \eqref{dissipation6-4} and \eqref{dissipation6-5} by $a^2_{5}$ and $ -a_5 (a_{4}^2-1)$, respectively,
and combining the resultant equations, we have
\begin{multline}\label{dissipation6-6}
 \pa_t   E_1 
+ a_5^2( |\hat u_{1}|^2 - |\hat u_{4}|^2)  \\
+ a_5 (a_4^2 + a_5^2 - 1) \xi \, \Re (i \hat{u}_{1} \bar{\hat{u}}_{5}) 
- a_{5}  a_6 (a_{4}^2-1) \xi \, \Re (i \hat{u}_{2} \bar{\hat{u}}_{6}) =  0,
\end{multline}
where we have defined that
$E_1 :=
- \Re \big\{ a_5^2(\hat u_{1} \bar{\hat{u}}_4+ a_4 \hat u_{2} \bar{\hat{u}}_3)
 - a_5 (a_{4}^2-1)\hat u_{2} \bar{\hat{u}}_5\big\}. 
$


%
Next, we multiply the first and second equations in \eqref{equations6} 
by $- i \xi \bar{\hat{u}}_{2}$ and $i \xi \bar{\hat{u}}_{1}$, respectively.
Then, combining the resultant equations and taking the real part, we have
\begin{equation}\label{dissipation6-7'}
\xi \pa_t  E_2 + \xi^2( |\hat u_{2}|^2 - |\hat u_{1}|^2)  
+ \xi \, \Re (i \hat{u}_{2} \bar{\hat{u}}_4) = 0,
\end{equation}
where
$E_2 := - \Re (i \hat u_{1} \bar{\hat{u}}_2)$.
Therefore, by Young inequality, the above equation becomes
\begin{equation}\label{dissipation6-7}
\xi \pa_t  E_2 + \frac{1}{2}\xi^2|\hat u_{2}|^2  \le \xi^2 |\hat u_{1}|^2 +\frac{1}{2} |\hat{u}_4|^2.
\end{equation}


%
We multiply the third and fourth equations in \eqref{equations6} 
by $i \xi a_4 \bar{\hat{u}}_{4}$ and $- i \xi a_4 \bar{\hat{u}}_{3}$, respectively.
Then, combining the resultant equations and taking the real part, we have
\begin{equation*}
 a_4 \xi \pa_t  \Re (i \hat u_{3} \bar{\hat{u}}_4) + a_4^2 \xi^2( |\hat u_{3}|^2 - |\hat u_{4}|^2)  
+ a_4 a_5  \xi^2 \, \Re (\hat{u}_{3} \bar{\hat{u}}_5) + a_4 \xi \, \Re (i \hat{u}_{1} \bar{\hat{u}}_3)= 0.
\end{equation*}
On the other hand, we multiply the second and third equations in \eqref{equations} 
by $- a_4 \bar{\hat{u}}_{3}$ and $-a_4 \bar{\hat{u}}_{2}$, respectively.
Then, combining the resultant equations and taking the real part, we have
\begin{equation*}
- a_4 \pa_t  \Re (\hat u_{2} \bar{\hat{u}}_3) 
- a_4  \xi \, \Re (i \hat{u}_{1} \bar{\hat{u}}_3) + a_4^2 \xi \, \Re (i \hat{u}_{2} \bar{\hat{u}}_4)= 0.
\end{equation*}
Finally, combining the above two equations, we get
\begin{equation}\label{dissipation6-8'}
\pa_t   \big\{ \xi E_3 + F_1 \big\} + a_4^2 \xi^2( |\hat u_{3}|^2 - |\hat u_{4}|^2)  
+ a_4 a_5  \xi^2 \, \Re (\hat{u}_{3} \bar{\hat{u}}_5) + a_4^2 \xi \, \Re (i \hat{u}_{2} \bar{\hat{u}}_4)= 0.
\end{equation}
where 
$E_3 :=  a_4 \Re (i \hat u_{3} \bar{\hat{u}}_4)$ and
$F_1 := - a_4  \Re (\hat u_{2} \bar{\hat{u}}_3)$.
By using Young inequality, we can obtain the following inequality:
\begin{equation}\label{dissipation6-8}
\pa_t   \big\{ \xi E_3 + F_1 \big\} + \frac{1}{2}a_4^2 \xi^2 |\hat u_{3}|^2  \le   a_4^2 \xi^2 |\hat u_{4}|^2
+ \frac{1}{2}a_5^2  \xi^2 |\hat{u}_5|^2 + a_4^2 |\xi| |\hat{u}_{2}| |\hat{u}_4|.
\end{equation}


%
Multiplying the fourth equation and fifth equation in \eqref{equations} 
by $i \xi a_{5} \bar{\hat{u}}_{5}$ and $ - i \xi a_{5} \bar{\hat{u}}_{4}$, respectively,
combining the resultant equations, and taking the real part,  then we have
\begin{multline}\label{dissipation6-9'}
\xi \pa_t  E_4 
+  a_{5}^2 \xi^2( |\hat u_{4}|^2 - |\hat u_{5}|^2)  \\
- a_{4} a_{5} \xi^2 \, \Re (\hat{u}_{3} \bar{\hat{u}}_5)
+ a_{5} a_{6} \xi^2 \, \Re (\hat{u}_{4} \bar{\hat{u}}_{6})
- a_{5}  \xi \, \Re (i \hat{u}_{1} \bar{\hat{u}}_{5}) = 0,
\end{multline}
where
$E_4 := a_{5} \Re (i \hat u_{4} \bar{\hat{u}}_5)$.
Here, by using Young inequality, we obtain 
\begin{multline}\label{dissipation6-9}
\xi \pa_t  E_4 
+  \frac{1}{2}a_{5}^2 \xi^2 |\hat u_{4}|^2  \\ 
 \le  a_{5}^2 \xi^2 |\hat u_{5}|^2 + \frac{1}{2} a_{6}^2 \xi^2 |\hat{u}_{6}|^2
+ a_{4} a_{5} \xi^2 \, \Re (\hat{u}_{3} \bar{\hat{u}}_5)
+ a_{5}  \xi \, \Re (i \hat{u}_{1} \bar{\hat{u}}_{5}).
\end{multline}


%
On the other hand, we multiply the fifth equation and the last equation in \eqref{equations6} 
by $i \xi a_{6} \bar{\hat{u}}_{6}$ and $ - i \xi a_{6} \bar{\hat{u}}_{5}$, respectively.
Then, combining the resultant equations and taking the real part, we obtain
\begin{equation*}
a_{6} \xi \pa_t  \Re (i \hat u_{5} \bar{\hat{u}}_6) +  a_{6}^2 \xi^2( |\hat u_{5}|^2 - |\hat u_{6}|^2)  \\
- a_{5} a_{6} \xi^2 \, \Re (\hat{u}_{4} \bar{\hat{u}}_6) + \gamma a_{6} \xi \, \Re (i \hat{u}_{5} \bar{\hat{u}}_6) = 0.
\end{equation*}
Using Young inequality, this yields
\begin{equation}\label{dissipation6-10}
a_{6} \xi \pa_t  \Re (i \hat u_{5} \bar{\hat{u}}_6) +  \frac{1}{2}a_{6}^2 \xi^2 |\hat u_{5}|^2 \\
\le a_{6}^2 \xi^2 |\hat u_{6}|^2  + \frac{1}{2} \gamma^2 |\hat{u}_6|^2
+ a_{5} a_{6} \xi^2 \, \Re (\hat{u}_{4} \bar{\hat{u}}_6). 
\end{equation}

\medskip

\noindent
{\bf Step 3.}\ \
In this step, we sum up the energy inequalities derived in the previous step, 
and then get the desired energy estimate.
Throughout this step, $\beta_j$ with $j \in \mathbb{N}$ denote the real numbers determined later. 
We first multiply \eqref{dissipation6-6} and \eqref{dissipation6-7} by $\xi^2$ and $\beta_1$, respectively.
Then we combine the resultant equation, obtaining
\begin{multline*}
 \pa_t  \big\{ \xi^2 E_1 + \beta_1 \xi E_2 \big\}
+ (a_5^2 - \beta_1) \xi^2 |\hat u_{1}|^2 + \frac{\beta_1}{2} \xi^2  |\hat u_{2}|^2 \\
\le  \Big(\frac{\beta_1}{2} + a_5^2 \xi^2 \Big)|\hat u_{4}|^2 
- a_5 (a_4^2 + a_5^2 - 1) \xi^3 \, \Re (i \hat{u}_{1} \bar{\hat{u}}_{5}) 
 + a_{5}  a_6 (a_{4}^2-1) \xi^3 \, \Re (i \hat{u}_{2} \bar{\hat{u}}_{6}).
\end{multline*}
Moreover, combining \eqref{dissipation6-6}, \eqref{dissipation6-8} and the above inequality, we have 
\begin{equation*}
\begin{split}
& \pa_t  \big\{ (1+ \xi^2) E_1 + \beta_1 \xi E_2 + \xi E_3 + F_1\big\} \\
&\qquad + \big\{a_5^2 + (a_5^2 - \beta_1) \xi^2\big\} |\hat u_{1}|^2 + \frac{\beta_1}{2} \xi^2  |\hat u_{2}|^2+ \frac{1}{2}a_4^2 \xi^2 |\hat u_{3}|^2 \\
&\le  \Big\{a_5^2+ \frac{\beta_1}{2} + (a_4^2 + a_5^2) \xi^2 \Big\}|\hat u_{4}|^2 + \frac{1}{2}a_5^2  \xi^2 |\hat{u}_5|^2  + a_4^2 |\xi| |\hat{u}_{2}| |\hat{u}_4| \\
&\qquad - a_5 (a_4^2 + a_5^2 - 1) \xi (1+\xi^2) \, \Re (i \hat{u}_{1} \bar{\hat{u}}_{5}) 
+ a_{5}  a_6 (a_{4}^2-1) \xi (1+\xi^2)\, \Re (i \hat{u}_{2} \bar{\hat{u}}_{6}).
\end{split}
\end{equation*}
For this inequality, letting $\beta_1$ suitably small and employing Young inequality, we can get
\begin{multline}\label{eneq6-1}
 \pa_t  \big\{(1+ \xi^2) E_1 + c \xi E_2 + \xi E_3 + F_1\big\} 
+  c(1 +  \xi^2) |\hat u_{1}|^2 + \beta_1 \xi^2  (|\hat u_{2}|^2+ |\hat u_{3}|^2) \\
\le  C(1+\xi^2)|\hat u_{4}|^2 + C  \xi^2 |\hat{u}_5|^2  \\
\qquad +  |a_4^2 + a_5^2 - 1| C |\xi|^3 |\hat{u}_{1}| |\hat{u}_{5}| 
+ |a_{4}^2-1| C |\xi| (1+\xi^2) |\hat{u}_{2}| |\hat{u}_{6}|.
\end{multline}
Similarly, multiplying \eqref{dissipation6-9} and \eqref{eneq6-1} by $1+\xi^2$ and $\beta_2 \xi^2$, respectively.
Then we combine the resultant equation, obtainig
\begin{equation*}
\begin{split}
& \pa_t  \big\{ \beta_2 \xi^2 ((1+ \xi^2) E_1 + \beta_1 \xi E_2 + \xi E_3 + F_1) + \xi (1+\xi^2)E_4 \big\}  \\
&\qquad  +  \beta_2 c \xi^2 (1 +  \xi^2) |\hat u_{1}|^2 + \beta_2 c \xi^4  (|\hat u_{2}|^2+ |\hat u_{3}|^2) + \Big(\frac{1}{2}a_{5}^2 - \beta_2 C \Big) \xi^2(1+\xi^2) |\hat u_{4}|^2  \\
&\le  \beta_2 C  \xi^4 |\hat{u}_5|^2  +  a_{5}^2 \xi^2(1+\xi^2) |\hat u_{5}|^2 + \frac{1}{2} a_{6}^2 \xi^2(1+\xi^2) |\hat{u}_{6}|^2 \\
&\qquad + a_{4} a_{5} \xi^2(1+\xi^2) \, \Re (\hat{u}_{3} \bar{\hat{u}}_5)
+ a_{5}  \xi(1+\xi^2) \, \Re (i \hat{u}_{1} \bar{\hat{u}}_{5}) \\
&\qquad +  \beta_2 |a_4^2 + a_5^2 - 1| C |\xi|^5 |\hat{u}_{1}| |\hat{u}_{5}| 
+ \beta_2 |a_{4}^2-1| C |\xi|^3 (1+\xi^2) |\hat{u}_{2}| |\hat{u}_{6}|.
\end{split}
\end{equation*}
Letting $\beta_2$ suitably small and using Young inequality derive that 
\begin{multline}\label{eneq6-2}
 \pa_t  \big\{ \beta_2 \xi^2 ((1+ \xi^2) E_1 + \beta_1 \xi E_2 + \xi E_3 + F_1) + \xi (1+\xi^2)E_4 \big\}  \\
+  c \xi^2 (1 +  \xi^2) (|\hat u_{1}|^2 + |\hat u_{4}|^2) + c \xi^4  (|\hat u_{2}|^2+ |\hat u_{3}|^2)   \\
\le   C (1+\xi^2)^2 |\hat u_{5}|^2 + C \xi^2(1+\xi^2) |\hat{u}_{6}|^2 \\
+  |a_4^2 + a_5^2 - 1| C \xi^6 |\hat{u}_{5}|^2 
+ |a_{4}^2-1| C |\xi|^2 (1+\xi^2)^2 |\hat{u}_{2}| |\hat{u}_{6}|.
\end{multline}
If we assume that $a_{4}^2-1 = 0$, 
the estimate \eqref{eneq6-2} can be  rewritten as
\begin{multline}\label{eneq6-2'}
 \pa_t  \big\{ \beta_2 \xi^2 ((1+ \xi^2) E_1 + \beta_1 \xi E_2 + \xi E_3 + F_1) + \xi (1+\xi^2)E_4 \big\}  \\
 +  c \xi^2 (1 +  \xi^2) (|\hat u_{1}|^2 + |\hat u_{4}|^2) + c \xi^4  (|\hat u_{2}|^2+ |\hat u_{3}|^2) \\
\le   C (1+\xi^2)^3 |\hat u_{5}|^2 + C \xi^2(1+\xi^2) |\hat{u}_{6}|^2. 
\end{multline}
Then, multiplying \eqref{dissipation6-10} and the above inequality by $(1+\xi^2)^3$ and $\beta_3 \xi^2$, respectively,
and combining the resultant equation, we have
\begin{multline*}
 \pa_t  \big\{ \beta_3 \xi^2 (\beta_2 \xi^2 ((1+ \xi^2) E_1 + \beta_1 \xi E_2 + \xi E_3 + F_1) + \xi (1+\xi^2)E_4) + \xi (1+\xi^2)^3 E_5 \big\}  \\
 +   \beta_3 c \xi^4 (1 +  \xi^2) (|\hat u_{1}|^2 + |\hat u_{4}|^2)+ \beta_3 c \xi^6  (|\hat u_{2}|^2+ |\hat u_{3}|^2)\\
 +  \Big( \frac{1}{2}a_{6}^2 - \beta_3 C\Big)  \xi^2 (1+\xi^2)^3 |\hat u_{5}|^2 \le  \beta_3  C \xi^4(1+\xi^2) |\hat{u}_{6}|^2 \\
 +  \Big( a_{6}^2 \xi^2  + \frac{1}{2} \gamma^2 \Big) (1+\xi^2)^3|\hat{u}_6|^2
+ a_{5} a_{6} \xi^2 (1+\xi^2)^3 \, \Re (\hat{u}_{4} \bar{\hat{u}}_6). 
\end{multline*}
Hence we arrive at
\begin{equation*}
\begin{split}
& \pa_t  \big\{ \beta_3 \xi^2 (\beta_2 \xi^2 ((1+ \xi^2) E_1 + \beta_1 \xi E_2 + \xi E_3 + F_1) \\
&\qquad\qquad
 + \xi (1+\xi^2)E_4) + \xi (1+\xi^2)^3 E_5 \big\}  \\
& +    c \xi^4 (1 +  \xi^2) (|\hat u_{1}|^2 + |\hat u_{4}|^2)+ c \xi^6  (|\hat u_{2}|^2+ |\hat u_{3}|^2)  +  c \xi^2 (1+\xi^2)^3 |\hat u_{5}|^2 \\
&\le   C (1+\xi^2)^4 |\hat{u}_{6}|^2 + C \xi^2 (1+\xi^2)^3 |\hat{u}_{4}| | |\hat{u}_6| .
\end{split}
\end{equation*}
Moreover, we multiply \eqref{dissipation6-8} and \eqref{dissipation6-9} by $\beta_4 \xi^6$ and $\beta_5 \xi^6$, respectively,
and combining the resultant equations and the above inequality.
Then, letting $\beta_4$ and $\beta_5$ suitably small,  this yields
\begin{multline}\label{eneq6-Eq}
 \partial_t E
 +    c \xi^4 (1 +  \xi^2) |\hat u_{1}|^2 + c \xi^6  |\hat u_{2}|^2+ c \xi^6 (1 +  \xi^2) |\hat u_{3}|^2    \\
+c \xi^4 (1 +  \xi^2)^2 |\hat u_{4}|^2  +  c \xi^2 (1+\xi^2)^3 |\hat u_{5}|^2 
\le   C (1+\xi^2)^4 |\hat{u}_{6}|^2,
\end{multline}
where we have defined 
\begin{multline}\label{eneq6-E}
E    = \beta_2 \beta_3 \xi^4 (1+ \xi^2) E_1 + \beta_1\beta_2 \beta_3 \xi^5 E_2 +  \xi^4 (\beta_2 \beta_3  + \beta_4 \xi^2)(\xi E_3 +  F_1) \\
 + \xi^3 ( \beta_3 (1+\xi^2) + \beta_5 \xi^4) E_4 + \xi (1+\xi^2)^3 E_5.
\end{multline}
Finally, combining the basic energy \eqref{eq6} with the above estimate, this yields
\begin{multline}\label{eneq6-3}
\partial _t \Big\{ \frac{1}{2}(1+\xi^2)^{4}|\hat{u}|^2 + \beta_{7} E \Big\} 
  +    c \xi^{4} (1 +  \xi^2) |\hat u_{1}|^2  \\
+ c \xi^{6}  |\hat u_{2}|^2 + c \sum_{j=3}^{6}\xi^{2(6-j)} (1+\xi^2)^{j-2} |\hat{u}_{j}|^2 \le 0.
\end{multline}
Thus, integrating the above estimate with respect to $t$,
we obtain the following energy estimate
\begin{multline}\label{eneq6-4}
 |\hat{u}(t,\xi)|^2 
+ \int^t_0   \Big\{
\frac{\xi^{4}}{(1+\xi^2)^{3}} |\hat u_{1}|^2  
+  \frac{\xi^{6}}{(1+\xi^2)^{4}} |\hat u_{2}|^2 \\
 + \sum_{j=3}^{6}  \frac{\xi^{2(6-j)}}{(1+\xi^2)^{6-j}} |\hat u_{j}|^2 
  \Big\} d\tau \le  C|\hat{u}(0,\xi)|^2.
\end{multline}
Here we used the following inequality
\begin{equation}\label{eneq6-5}
c |\hat{u}|^2 \le 
\frac{1}{2} |\hat{u}|^2  +   \frac{\beta_{7}}{(1+\xi^2)^{4}} E
\le C |\hat{u}|^2
\end{equation}
for suitably small $\beta_{7}$.
Furthermore the estimate \eqref{eneq6-3} with \eqref{eneq6-5} gives us the following pointwise estimate
\begin{equation}\label{pt6-1}
|\hat{u}(t,\xi)|
\le C e^{- c \lambda(\xi)t} |\hat{u}(0,\xi)|,
\qquad
\lambda(\xi) = 
\frac{\xi^{6}}{(1+\xi^2)^{4}}.
\end{equation}

On the other hand, 
if we assume that $a_{4}^2 + a_5^2 -1 = 0$, 
the estimate \eqref{eneq6-2} is rewritten as
\begin{equation}\label{eneq6-2''}
\begin{split}
& \pa_t  \big\{ \beta_2 \xi^2 ((1+ \xi^2) E_1 + \beta_1 \xi E_2 + \xi E_3 + F_1) + \xi (1+\xi^2)E_4 \big\}  \\
&\qquad  +  c \xi^2 (1 +  \xi^2) (|\hat u_{1}|^2 + |\hat u_{4}|^2) + c \xi^4  (|\hat u_{2}|^2+ |\hat u_{3}|^2) \\
&\le   C (1+\xi^2)^2 |\hat u_{5}|^2 + C (1+\xi^2)^4 |\hat{u}_{6}|^2. \\
\end{split}
\end{equation}
Then, multiplying \eqref{dissipation6-10} and the above inequality by $(1+\xi^2)^2$ and $\beta_3 \xi^2$, respectively,
and combining the resultant equation, we have
\begin{multline*}
 \pa_t  \big\{ \beta_3 \xi^2 (\beta_2 \xi^2 ((1+ \xi^2) E_1 + \beta_1 \xi E_2 + \xi E_3 + F_1) + \xi (1+\xi^2)E_4) + \xi (1+\xi^2)^2 E_5 \big\}  \\
+   \beta_3 c \xi^4 (1 +  \xi^2) (|\hat u_{1}|^2 + |\hat u_{4}|^2)+ \beta_3 c \xi^6  (|\hat u_{2}|^2+ |\hat u_{3}|^2)
 +  \Big( \frac{1}{2}a_{6}^2 - \beta_3 C\Big)  \xi^2 (1+\xi^2)^2 |\hat u_{5}|^2 \\
\le  \beta_3  C (1+\xi^2)^4 |\hat{u}_{6}|^2 
 +  \Big( a_{6}^2 \xi^2  + \frac{1}{2} \gamma^2 \Big) (1+\xi^2)^2|\hat{u}_6|^2
+ a_{5} a_{6} \xi^2 (1+\xi^2)^2 \, \Re (\hat{u}_{4} \bar{\hat{u}}_6). 
\end{multline*}
Hence we arrive at
\begin{equation*}
\begin{split}
& \pa_t  \big\{ \beta_3 \xi^2 (\beta_2 \xi^2 ((1+ \xi^2) E_1 + \beta_1 \xi E_2 + \xi E_3 + F_1) \\
&\qquad\qquad
 + \xi (1+\xi^2)E_4) + \xi (1+\xi^2)^2 E_5 \big\}  \\
& +    c \xi^4 (1 +  \xi^2) (|\hat u_{1}|^2 + |\hat u_{4}|^2)+ c \xi^6  (|\hat u_{2}|^2+ |\hat u_{3}|^2)  +  c \xi^2 (1+\xi^2)^2 |\hat u_{5}|^2 \\
&\le   C (1+\xi^2)^4 |\hat{u}_{6}|^2.
\end{split}
\end{equation*}
Moreover, we multiply \eqref{dissipation6-8}, \eqref{dissipation6-9} and \eqref{dissipation6-10} by $\beta_4 \xi^6$, $\beta_5 \xi^6$ and $\beta_6 \xi^6$, respectively,
and combine the resultant equations and the above inequality.
Then, letting $\beta_4$ and $\beta_5$ suitably small,  this yields
\begin{equation*}
\begin{split}
& \pa_t  \big\{ \beta_2 \beta_3 \xi^4 (1+ \xi^2) E_1 + \beta_1\beta_2 \beta_3 \xi^5 E_2 +  \xi^4 (\beta_2 \beta_3  + \beta_4 \xi^2)(\xi E_3 +  F_1) \\
&\qquad\qquad
 + \xi^3 ( \beta_3 (1+\xi^2) + \beta_5 \xi^4) E_4 + \xi ((1+\xi^2)^2 + \beta_6 \xi^6) E_5 \big\}  \\
& +    c \xi^4 (1 +  \xi^2) |\hat u_{1}|^2 + c \xi^6  |\hat u_{2}|^2+ c \xi^6 (1 +  \xi^2) |\hat u_{3}|^2    \\
&+c \xi^4 (1 +  \xi^2)^2 |\hat u_{4}|^2  +  c \xi^2 (1+\xi^2)^3 |\hat u_{5}|^2 
\le   C (1+\xi^2)^4 |\hat{u}_{6}|^2.
\end{split}
\end{equation*}
We note that this estimate is essentially the same as \eqref{eneq6-Eq}.
Hence we can obtain the energy estimate \eqref{eneq6-4} and the pointwise estimate \eqref{pt6-1}.
Eventually, we arrive at the estimate for both cases $a_4^2 -  1 = 0$ and $a_4^2 + a_5^2 - 1 =0$. 
Moreover, by using the similar argument, we can derive the same estimates in the case $a_4^2 - 1 \neq  0$, $a_4^2 + a_5^2 - 1 \neq 0$. 
Thus we complete the proof of Theorem 2.1 with $m=6$.


\subsection{Energy method for model I}

Inspired by the concrete computation in Subsection 2.2,
we consider the more general case $m\geq 6$.
Now, our system \eqref{sys1} with \eqref{Tmat} is described as 
\begin{equation}\label{equations}
\begin{split}
&\pa_t \hat u_1 + i \xi \hat u_2 + \hat u_4 = 0, \\
&\pa_t \hat u_2 + i \xi \hat u_1 = 0, \\
&\pa_t \hat u_3 + i \xi a_4 \hat u_4 = 0, \\
&\pa_t \hat u_4 + i \xi (a_4 \hat u_3 + a_5 \hat u_5) - \hat u_1 = 0, \\
&\pa_t \hat u_j + i \xi (a_{j} \hat u_{j-1} + a_{j+1} \hat u_{j+1}) = 0, \qquad j = 5, \cdots, m-1, \\
&\pa_t \hat u_m + i \xi a_{m} \hat u_{m-1} + \gamma \hat u_m = 0. 
\end{split}
\end{equation}
We are going to apply the energy method to this system and derive Theorem \ref{thm1}.
The  proof is organized by the following three steps.

\medskip

\noindent
{\bf Step 1.}\ \
We first derive the basic energy equality for the system \eqref{sys1} in the Fourier space.
Taking the inner product of \eqref{sys1} with $\hat{u}$, we have
\begin{equation*}
	\langle \hat u_t, \hat{u} \rangle
	+ i \xi \langle A_m \hat u, \hat u \rangle
	+ \langle L_m \hat u, \hat u \rangle = 0.
\end{equation*}
Taking the real part, we get the basic energy equality
\begin{equation*}
\frac{1}{2}\frac{\partial }{\partial t} |\hat u|^2
+ \langle L_m \hat u, \hat u \rangle = 0,
\end{equation*}
and hence
\begin{equation}\label{eq}
\frac{1}{2}\pa_t  |\hat u|^2
+ \gamma \hat u_m^2 = 0.
\end{equation}
Next we create the dissipation terms by the following two steps.

\medskip

\noindent
{\bf Step 2.}\ \
For $\ell = 6, \cdots , m-1$,
we multiply the fifth equations with $j=\ell-1$ and $j=\ell$ in \eqref{equations} 
by $i \xi a_{\ell} \bar{\hat{u}}_{\ell}$ and $ - i \xi a_{\ell} \bar{\hat{u}}_{\ell-1}$, respectively.
Then, combining the resultant equations and taking the real part, we have
\begin{multline}\label{dissipation-1'}
a_{\ell} \xi \pa_t  \Re (i \hat u_{\ell-1} \bar{\hat{u}}_\ell) +  a_{\ell}^2 \xi^2( |\hat u_{\ell-1}|^2 - |\hat u_{\ell}|^2)  \\
- a_{\ell} a_{\ell-1} \xi^2 \, \Re (\hat{u}_{\ell-2} \bar{\hat{u}}_\ell) + a_{\ell} a_{\ell+1} \xi^2 \, \Re (\hat{u}_{\ell-1} \bar{\hat{u}}_{\ell+1}) = 0.
\end{multline}
Here, by using Young inequality, we obtain 
\begin{equation}\label{dissipation-1}
\xi \partial_t E_{\ell-1} +  \frac{1}{2}a_{\ell}^2 \xi^2 |\hat u_{\ell-1}|^2  
\le  a_{\ell}^2 \xi^2 |\hat u_{\ell}|^2 + \frac{1}{2}a_{\ell+1}^2 \xi^2 |\hat{u}_{\ell+1}|^2 
+ a_{\ell} a_{\ell-1} \xi^2 \, \Re (\hat{u}_{\ell-2} \bar{\hat{u}}_\ell)
\end{equation}
for $\ell = 6, \cdots , m-1$, where we have defined $E_{\ell -1} = a_{\ell} \xi \Re (i \hat u_{\ell-1} \bar{\hat{u}}_\ell)$.
On the other hand, we multiply the fifth equation with $j=m-1$ and the last equation in \eqref{equations} 
by $i \xi a_{m} \bar{\hat{u}}_{m}$ and $ - i \xi a_{m} \bar{\hat{u}}_{m-1}$, respectively.
Then, combining the resultant equations and taking the real part, we obtain
\begin{multline}\label{dissipation-2'}
a_{m} \xi \pa_t  \Re (i \hat u_{m-1} \bar{\hat{u}}_m) +  a_{m}^2 \xi^2( |\hat u_{m-1}|^2 - |\hat u_{m}|^2)  \\
- a_{m} a_{m-1} \xi^2 \, \Re (\hat{u}_{m-2} \bar{\hat{u}}_m) + \gamma a_{m} \xi \, \Re (i \hat{u}_{m-1} \bar{\hat{u}}_m) = 0.
\end{multline}
Using Young inequality, this yields
\begin{multline}\label{dissipation-2}
 \xi \partial_t E_{m-1} +  \frac{1}{2}a_{m}^2 \xi^2 |\hat u_{m-1}|^2 \\
 \le a_{m}^2 \xi^2 |\hat u_{m}|^2  + \frac{1}{2} \gamma^2 |\hat{u}_m|^2
+ a_{m} a_{m-1} \xi^2 \, \Re (\hat{u}_{m-2} \bar{\hat{u}}_m),
\end{multline}
where we have defined $E_{m -1} = a_{m} \xi \Re (i \hat u_{m-1} \bar{\hat{u}}_m)$.

\medskip

\noindent
{\bf Step 3.}\ \
We note that equations \eqref{equations} with $1 \le j \le 5$ are the same as the five equations in \eqref{equations6}. 
Thus we can adopt the useful estimates derived in Subsection 2.2.
More precisely, we employ \eqref{dissipation6-8}, \eqref{dissipation6-9} and \eqref{eneq6-2} again.

For the estimate \eqref{eneq6-2}, if we assume that $a_{4}^2-1 = 0$,  we can obtain \eqref{eneq6-2'}.
Then, multiplying \eqref{dissipation-1} with $\ell = 6$ and \eqref{eneq6-2'} by $(1+\xi^2)^3$ and $\beta_3 \xi^2$, respectively,
and combining the resultant equation, we have
\begin{multline*}
 \partial _t \big\{ \beta_3 \xi^2 (\beta_2 \xi^2 ((1+ \xi^2) E_1 + \beta_1 \xi E_2 + \xi E_3 + F_1) + \xi (1+\xi^2)E_4) + \xi (1+\xi^2)^3 E_5 \big\}  \\
 +   \beta_3 c \xi^4 (1 +  \xi^2) (|\hat u_{1}|^2 + |\hat u_{4}|^2)+ \beta_3 c \xi^6  (|\hat u_{2}|^2+ |\hat u_{3}|^2)\\
 +  \Big( \frac{1}{2}a_{6}^2 - \beta_3 C\Big)  \xi^2 (1+\xi^2)^3 |\hat u_{5}|^2 
\le  \beta_3  C \xi^4(1+\xi^2) |\hat{u}_{6}|^2 
 + a_{6}^2 \xi^2  (1+\xi^2)^3|\hat{u}_6|^2\\
+ \frac{1}{2} a_{7}^2 \xi^2  (1+\xi^2)^3|\hat{u}_7|^2
+ a_{5} a_{6} \xi^2 (1+\xi^2)^3 \, \Re (\hat{u}_{4} \bar{\hat{u}}_6). 
\end{multline*}
Hence we arrive at
\begin{equation*}
\begin{split}
&\partial_t \big\{ \beta_3 \xi^2 (\beta_2 \xi^2 ((1+ \xi^2) E_1 + \beta_1 \xi E_2 + \xi E_3 + F_1) \\
&\qquad\qquad
 + \xi (1+\xi^2)E_4) + \xi (1+\xi^2)^3 E_5 \big\}  \\
& +    c \xi^4 (1 +  \xi^2) (|\hat u_{1}|^2 + |\hat u_{4}|^2)+ c \xi^6  (|\hat u_{2}|^2+ |\hat u_{3}|^2)  +  c \xi^2 (1+\xi^2)^3 |\hat u_{5}|^2 \\
&\le   C \xi^2 (1+\xi^2)^3 (|\hat{u}_{6}|^2 +  |\hat{u}_{7}|^2) + C \xi^2 (1+\xi^2)^3 |\hat{u}_{4}| | |\hat{u}_6| .
\end{split}
\end{equation*}
Moreover, we multiply \eqref{dissipation6-8} and \eqref{dissipation6-9} by $\beta_4 \xi^6$ and $\beta_5 \xi^6$, respectively,
and combining the resultant equations and the above inequality.
Then, letting $\beta_4$ and $\beta_5$ suitably small,  this yields
\begin{multline}\label{eneq-1}
\partial _t E  
 +    c \xi^4 (1 +  \xi^2) |\hat u_{1}|^2 + c \xi^6  |\hat u_{2}|^2+ c \xi^6 (1 +  \xi^2) |\hat u_{3}|^2   +c \xi^4 (1 +  \xi^2)^2 |\hat u_{4}|^2 \\
  +  c \xi^2 (1+\xi^2)^3 |\hat u_{5}|^2 
\le   C (1+\xi^2)^4 |\hat{u}_{6}|^2  + C \xi^2 (1+\xi^2)^3 |\hat{u}_{7}|^2.
\end{multline}
where $E$ is defined in \eqref{eneq6-E}.

On the other hand, 
if we assume that $a_{4}^2 + a_5^2 -1 = 0$, 
we employ \eqref{eneq6-2''}.
Then, multiplying \eqref{dissipation-1} with $\ell = 6$ and \eqref{eneq6-2''} by $(1+\xi^2)^2$ and $\beta_3 \xi^2$, respectively,
and combining the resultant equation, we have
\begin{multline*}
 \partial_t \big\{ \beta_3 \xi^2 (\beta_2 \xi^2 ((1+ \xi^2) E_1 + \beta_1 \xi E_2 + \xi E_3 + F_1) + \xi (1+\xi^2)E_4) + \xi (1+\xi^2)^2 E_5 \big\}  \\
 +   \beta_3 c \xi^4 (1 +  \xi^2) (|\hat u_{1}|^2 + |\hat u_{4}|^2)+ \beta_3 c \xi^6  (|\hat u_{2}|^2+ |\hat u_{3}|^2)\\
 +  \Big( \frac{1}{2}a_{6}^2 - \beta_3 C\Big)  \xi^2 (1+\xi^2)^2 |\hat u_{5}|^2 
\le  \beta_3  C (1+\xi^2)^4 |\hat{u}_{6}|^2 
 +  a_{6}^2 \xi^2  (1+\xi^2)^2|\hat{u}_6|^2\\
  +  \frac{1}{2} a_{7}^2 \xi^2  (1+\xi^2)^2|\hat{u}_7|^2
+ a_{5} a_{6} \xi^2 (1+\xi^2)^2 \, \Re (\hat{u}_{4} \bar{\hat{u}}_6). 
\end{multline*}
Hence we arrive at
\begin{equation*}
\begin{split}
& \partial_t \big\{ \beta_3 \xi^2 (\beta_2 \xi^2 ((1+ \xi^2) E_1 + \beta_1 \xi E_2 + \xi E_3 + F_1) \\
&\qquad\qquad
 + \xi (1+\xi^2)E_4) + \xi (1+\xi^2)^2 E_5 \big\}  \\
& +    c \xi^4 (1 +  \xi^2) (|\hat u_{1}|^2 + |\hat u_{4}|^2)+ c \xi^6  (|\hat u_{2}|^2+ |\hat u_{3}|^2)  +  c \xi^2 (1+\xi^2)^2 |\hat u_{5}|^2 \\
&\le   C (1+\xi^2)^4 |\hat{u}_{6}|^2  +  C\xi^2  (1+\xi^2)^2|\hat{u}_7|^2.
\end{split}
\end{equation*}
Moreover, we multiply \eqref{dissipation6-8}, \eqref{dissipation6-9} and \eqref{dissipation-1} with $\ell = 6$ by $\beta_4 \xi^6$, $\beta_5 \xi^6$ and $\beta_6 \xi^6$, respectively,
and combine the resultant equations and the above inequality.
Then, letting $\beta_4$, $\beta_5$ and $\beta_6$ suitably small,  this yields
\begin{equation*}
\begin{split}
& \pa_t  \big\{ \beta_2 \beta_3 \xi^4 (1+ \xi^2) E_1 + \beta_1\beta_2 \beta_3 \xi^5 E_2 +  \xi^4 (\beta_2 \beta_3  + \beta_4 \xi^2)(\xi E_3 +  F_1) \\
&\qquad\qquad
 + \xi^3 ( \beta_3 (1+\xi^2) + \beta_5 \xi^4) E_4 + \xi ((1+\xi^2)^2 + \beta_6 \xi^6) E_5 \big\}  \\
& +    c \xi^4 (1 +  \xi^2) |\hat u_{1}|^2 + c \xi^6  |\hat u_{2}|^2+ c \xi^6 (1 +  \xi^2) |\hat u_{3}|^2    \\
&+c \xi^4 (1 +  \xi^2)^2 |\hat u_{4}|^2  +  c \xi^2 (1+\xi^2)^3 |\hat u_{5}|^2 
\le   C (1+\xi^2)^4 |\hat{u}_{6}|^2  +  C\xi^2  (1+\xi^2)^3|\hat{u}_7|^2.
\end{split}
\end{equation*}
Consequently, this estimate is essentially the same as \eqref{eneq-1}.
Moreover, by using the similar argument, we can derive the same estimate in the case $a_4^2 - 1 \neq 1$ and $a_4^2 + a_5^2 - 1 \neq 0$. 

\medskip

By using the estimate \eqref{eneq-1}, we construct the desired estimate.
We multiply \eqref{dissipation-1} with $\ell = 7$ and \eqref{eneq-1} by $(1+\xi^2)^4$ and $\beta_7 \xi^2$, respectively,
and combine the resultant equation.
Moreover, letting $\beta_7$ suitably small and using Young inequality, we obtain
\begin{equation*}
\begin{split}
&\partial _t \big\{ \beta_7 \xi^2 E + \xi (1+\xi^2)^4 E_6\big\}
  +    c \xi^{6} (1 +  \xi^2) |\hat u_{1}|^2 + c \xi^{8}  |\hat u_{2}|^2 + c \xi^{8} (1 +  \xi^2) |\hat u_{3}|^2  \\
  &+c \xi^{6} (1 +  \xi^2)^2 |\hat u_{4}|^2+  c \xi^{4} (1+\xi^2)^3 |\hat u_{5}|^2  + c \xi^{2} (1+\xi^2)^4 |\hat{u}_{6}|^2  \\
&\le  C(1+\xi^2)^5 |\hat{u}_{7}|^2 + C \xi^2 (1+\xi^2)^4 |\hat{u}_{8}|^2.
\end{split}
\end{equation*}
Eventually, by the induction argument with respect to $j$ in \eqref{dissipation-1},
we can derive
\begin{multline}\label{eneq-2}
\partial _t \mathcal{E}_{m-2} 
  +    c \xi^{2(m - 5)} (1 +  \xi^2) |\hat u_{1}|^2 + c \xi^{2(m-4)}  |\hat u_{2}|^2 \\
+ c \sum_{j=3}^{m-2}\xi^{2(m-j-1)} (1+\xi^2)^{j-2} |\hat{u}_{j}|^2  \\
\le  C(1+\xi^2)^{m-3} |\hat{u}_{m-1}|^2 + C \xi^2 (1+\xi^2)^{m-4} |\hat{u}_{m}|^2.
\end{multline}
for $m \ge 7$.
Here we define $\mathcal{E}_{m-2}$ as $\mathcal{E}_5 = E$ and 
$$
\mathcal{E}_{m-2} = \beta_{m-1} \xi^2 \mathcal{E}_{m -3} + \xi (1+\xi^2)^{m-4} E_{m-2},
\qquad
m \ge 8.
$$
Now, multiplying \eqref{dissipation-2} and \eqref{eneq-2} by $(1+\xi^2)^{m-3}$ and $\beta_{m} \xi^2$, respectively,
and making the appropriate  combination, we get
\begin{multline}\label{eneq-3}
\partial _t \mathcal{E}_{m-1} 
  +    c \xi^{2(m - 4)} (1 +  \xi^2) |\hat u_{1}|^2 + c \xi^{2(m-3)}  |\hat u_{2}|^2 \\
+ c \sum_{j=3}^{m-1}\xi^{2(m-j)} (1+\xi^2)^{j-2} |\hat{u}_{j}|^2 \le  C(1+\xi^2)^{m-2} |\hat{u}_{m}|^2.
\end{multline}
Finally, combining \eqref{eq} with \eqref{eneq-3}, this yields
\begin{multline}\label{eneq-4}
\partial _t \Big\{ \frac{1}{2}(1+\xi^2)^{m-2}|\hat{u}|^2 + \beta_{m+1} \mathcal{E}_{m-1}\Big\} 
  +    c \xi^{2(m - 4)} (1 +  \xi^2) |\hat u_{1}|^2  \\
+ c \xi^{2(m-3)}  |\hat u_{2}|^2 + c \sum_{j=3}^{m}\xi^{2(m-j)} (1+\xi^2)^{j-2} |\hat{u}_{j}|^2 \le 0.
\end{multline}
Thus, integrating the above estimate with respect to $t$,
we obtain the following energy estimate
\begin{multline}\label{energy-eq1}
 |\hat{u}(t,\xi)|^2 
+ \int^t_0   \Big\{
\frac{\xi^{2(m-4)}}{(1+\xi^2)^{m-3}} |\hat u_{1}|^2  
+  \frac{\xi^{2(m-3)}}{(1+\xi^2)^{m-2}} |\hat u_{2}|^2 \\
 + \sum_{j=3}^{m}  \frac{\xi^{2(m-j)}}{(1+\xi^2)^{m-j}} |\hat u_{j}|^2 
  \Big\} d\tau \le  C|\hat{u}(0,\xi)|^2
\end{multline}
for $m \ge 7$.
Here we have used the following inequality
\begin{equation}
\notag
c |\hat{u}|^2 \le 
\frac{1}{2} |\hat{u}|^2  +   \frac{\beta_{m+1}}{(1+\xi^2)^{m-2}} \mathcal{E}_{m-1}
\le C |\hat{u}|^2
\end{equation}
for suitably small $\beta_{m+1}$.
Furthermore the estimate \eqref{eneq-3} with \eqref{eneq-4} gives us the following pointwise estimate
\begin{equation}
\notag
|\hat{u}(t,\xi)|
\le C e^{- c \lambda(\xi)t} |\hat{u}(0,\xi)|,
\qquad
\lambda(\xi) = 
\frac{\xi^{2(m-3)}}{(1+\xi^2)^{m-2}}
\end{equation}
for $m \ge 7$.
Therefore, together with the proof in Subsection 2.2, \eqref{point1} is proved, and we then complete the proof of Theorem 2.1.


\subsection{Construction of the matrices $K$ and  $S$}

In this section, inspired by the energy method employed in Sections 2.2 and 2.3, we shall derive the matrices $K$ and $S$.

Based on the energy method of Step 2 in Subsection 2.2, we introduce the following $m \times m$ matrices:
\begin{equation*}
S_1 =
\left(
\begin{array}{cccc:c}
 {0} & {0} & {0} & {1}  &    \\
 {0} & {0} & {0} & {0}  &    \\
 {0} & {0} & {0} & {0}  &    \mbox{\smash{\huge\textit{O}}}  \\
 {1} & {0} & {0} & {0}  &    \\
 \hdashline
     &      &      &      &       \\
       &     \mbox{\smash{\huge\textit{O}}}    &     &     &  \mbox{\smash{\huge\textit{O}}}     \\
 \end{array}
\right), \quad
S_2 =
\left(
\begin{array}{cccc:c}
 {0} & {0} & {0} & {0}  &    \\
 {0} & {0} & {1} & {0}  &    \\
 {0} & {1} & {0} & {0}  &    \mbox{\smash{\huge\textit{O}}}  \\
 {0} & {0} & {0} & {0}  &    \\
 \hdashline
     &      &      &      &       \\
       &     \mbox{\smash{\huge\textit{O}}}    &     &     &  \mbox{\smash{\huge\textit{O}}}     \\
 \end{array}
\right), \quad
S_3 =
\left(
\begin{array}{cccc:cc}
 {0} & {0} & {0} & {0}  &    {0} & \\
 {0} & {0} & {0} & {0}  &    {1} &  \\
 {0} & {0} & {0} & {0}  &    {0} & \mbox{\smash{\huge\textit{O}}}  \\
 {0} & {0} & {0} & {0}  &    {0} &  \\
 \hdashline
 {0} & {1} & {0}  & {0}  &   {0} &      \\
     &      &      &      &     &  \\
       &     \mbox{\smash{\huge\textit{O}}}    &     &     &   & \mbox{\smash{\huge\textit{O}}}      \\
 \end{array}
\right),
\end{equation*}
and hence
\begin{equation}
\notag
\begin{split}
\tilde{S} &= - a_5 \big\{ a_5 (S_1 + a_4 S_2) - a_5 (a_4^2-1) S_3 \big\}  \\[2mm]
&= - a_5 
\left(
\begin{array}{cccc:cc}
 {0} & {0} & {0} & {a_5}  &   {0} & \\
 {0} & {0} & {a_4 a_5} & {0}  &    {1-a_4^2} &  \\
 {0} & {a_4 a_5} & {0} & {0}  &    {0} & \mbox{\smash{\huge\textit{O}}}  \\
 {a_5} & {0} & {0} & {0}  &    {0} &  \\
 \hdashline
 {0} & {1-a_4^2} & {0}  & {0}  &   {0} &      \\
     &      &      &      &     &  \\
       &     \mbox{\smash{\huge\textit{O}}}    &     &     &   & \mbox{\smash{\huge\textit{O}}}      \\
 \end{array}
\right). \\
\end{split}
\end{equation}
Then, we multiply \eqref{Fsys1} by $\tilde{S}$ and take the inner product with $\hat u$.
Furthermore, taking the real part of the resultant equation, we obtain 
\begin{equation}\label{FsysS-1}
\frac{1}{2} \partial_t
\langle \tilde{S} \hat u, \hat u \rangle  
+  \xi  \langle i [\tilde{S}A_m]^{\rm asy} \hat u, \hat u \rangle + \langle [\tilde{S}L_m]^{\rm sy} \hat u, \hat u \rangle = 0, 
\end{equation}
where 
\begin{equation*}
\begin{split}
\tilde{S}A_m &= - a_5 
\left(
\begin{array}{cccc:ccc}
 {0} & {0} & {a_4 a_5} & {0}  &   {a_5^2} & {0} & \\
 {0} & {0} & {0} & {a_5}  &  {0} & {a_6(1-a_4^2)} &  \\
 {a_4a_5} & {0} & {0} & {0}  &  {0} &  {0} & \mbox{\smash{\huge\textit{O}}}  \\
 {0} & {a_5} & {0} & {0}  &  {0}  & {0} &  \\
 \hdashline
 {1-a_4^2} & {0} & {0}  & {0}  &   {0} & {0} &   \\
  {0} & {0} & {0}  & {0}  &   {0} & {0} &   \\
     &      &      &      &     &   & \\
       &   \mbox{\smash{\huge\textit{O}}}  &    &     &   &   & \mbox{\smash{\huge\textit{O}}}      \\
 \end{array}
\right), \\[2mm]
\tilde{S} L_m &= a_5^2 
\left(
\begin{array}{cccc:c}
 {1} & {0} & {0} & {0}  &    \\
 {0} & {0} & {0} & {0}  &     \\
 {0} & {0} & {0} & {0}  &     \mbox{\smash{\huge\textit{O}}}  \\
 {0} & {0} & {0} & {-1}  &     \\
 \hdashline
     &      &      &      &   \\
       &     \mbox{\smash{\huge\textit{O}}}    &     &     &  \mbox{\smash{\huge\textit{O}}}      \\
 \end{array}
\right). \\
\end{split}
\end{equation*}
The equality \eqref{FsysS-1} is equivalent to \eqref{dissipation6-6}.
We note that the symmetric matrix $S_1 + a_4 S_2$ is the key matrix for $4 \times 4$ Timoshenko system (see \cite{IHK08,IK08}).
The symmetric matrix $\tilde{S}$ is the one of the key matrix for the system \eqref{Fsys1}.


On the other hand, 
we introduce the following $m \times m$ matrix:
\begin{equation*}
K_1 =
\left(
\begin{array}{cccc:c}
 {0} & {-1} & {0} & {0}  &    \\
 {1} & {0} & {0} & {0}  &    \\
 {0} & {0} & {0} & {0}  &    \mbox{\smash{\huge\textit{O}}}  \\
 {0} & {0} & {0} & {0}  &    \\
 \hdashline
     &      &      &      &       \\
       &     \mbox{\smash{\huge\textit{O}}}    &     &     &  \mbox{\smash{\huge\textit{O}}}     \\
 \end{array}
\right).
\end{equation*}
Then, we multiply \eqref{Fsys1} by $-i\xi K_1$ and take the inner product with $\hat u$.
Moreover, taking the real part of the resultant equation, we have 
\begin{equation}\label{FsysK1-1}
-\frac{1}{2} \xi \partial_t
\langle i K_1 \hat u, \hat u \rangle 
 + \xi^2 \langle [K_1A_m]^{\rm sy} \hat u, \hat u \rangle
 -  \xi \langle i [K_1L_m]^{\rm asy} \hat u, \hat u \rangle = 0, 
\end{equation}
where 
\begin{equation*}
K_1A_m =  
\left(
\begin{array}{cccc:cc}
 {-1} & {0} & {0} & {0}   & \\
 {0} & {1} & {0} & {0}   &  \\
 {0} & {0} & {0} & {0}   & \mbox{\smash{\huge\textit{O}}}  \\
 {0} & {0} & {0} & {0}  &  \\
 \hdashline
     &      &      &      &     &  \\
       &   \mbox{\smash{\huge\textit{O}}}  &   &   &   \mbox{\smash{\huge\textit{O}}}      \\
 \end{array}
\right), \qquad
K_1 L_m =  
\left(
\begin{array}{cccc:c}
 {0} & {0} & {0} & {0}  &    \\
 {0} & {0} & {0} & {1}  &     \\
 {0} & {0} & {0} & {0}  &     \mbox{\smash{\huge\textit{O}}}  \\
 {0} & {0} & {0} & {0}  &     \\
 \hdashline
     &      &      &      &   \\
       &     \mbox{\smash{\huge\textit{O}}}    &     &     &  \mbox{\smash{\huge\textit{O}}}      \\
 \end{array}
\right). \\
\end{equation*}
The equality \eqref{FsysK1-1} is equivalent to \eqref{dissipation6-7'}. 


We next introduce the following $m \times m$ matrices:
\begin{equation*}
K_4 = a_4
\left(
\begin{array}{cccc:c}
 {0} & {0} & {0} & {0}  &    \\
 {0} & {0} & {0} & {0}  &    \\
 {0} & {0} & {0} & {1}  &    \mbox{\smash{\huge\textit{O}}}  \\
 {0} & {0} & {-1} & {0}  &    \\
 \hdashline
     &      &      &      &       \\
       &     \mbox{\smash{\huge\textit{O}}}    &     &     &  \mbox{\smash{\huge\textit{O}}}     \\
 \end{array}
\right),  \qquad
S_4 = - a_4
\left(
\begin{array}{cccc:c}
 {0} & {0} & {0} & {0}  &    \\
 {0} & {0} & {1} & {0}  &    \\
 {0} & {1} & {0} & {0}  &    \mbox{\smash{\huge\textit{O}}}  \\
 {0} & {0} & {0} & {0}  &    \\
 \hdashline
     &      &      &      &       \\
       &     \mbox{\smash{\huge\textit{O}}}    &     &     &  \mbox{\smash{\huge\textit{O}}}     \\
 \end{array}
\right). 
\end{equation*}
Then, we multiply \eqref{Fsys1} by $-i\xi K_2$ and $S_4$, 
and take the inner product with $\hat u$, respectively.
Moreover, taking the real part of the resultant equations and combining these, then we have 
\begin{multline}\label{FsysKS-1}
\frac{1}{2} \partial_t
\langle (S_4 - i \xi K_4) \hat u, \hat u \rangle 
+ \xi^2 \langle [K_4A_m]^{\rm sy} \hat u, \hat u \rangle+ \langle [S_4L_m]^{\rm sy} \hat u, \hat u \rangle  \\
 +  \xi  \langle i [S_4A_m - K_4L_m]^{\rm asy} \hat u, \hat u \rangle = 0, 
\end{multline}
where $S_4 L_m =O$ and
\begin{equation*}
\begin{split}
K_4A_m &=  
\left(
\begin{array}{cccc:cc}
 {0} & {0} & {0} & {0}  &   {0} & \\
 {0} & {0} & {0} & {0}  &  {0}  &  \\
 {0} & {0} & {a_4^2} & {0}  &  {a_4a_5}  & \mbox{\smash{\huge\textit{O}}}  \\
 {0} & {0} & {0} & {-a^2_4}  &  {0}  &  \\
 \hdashline
 {0} & {0} & {0}  & {0}  &   {0} &    \\
     &      &      &      &     &  \\
       &   \mbox{\smash{\huge\textit{O}}}  &   &   &   & \mbox{\smash{\huge\textit{O}}}      \\
 \end{array}
\right), \quad
S_4A_m - K_4 L_m = 
\left(
\begin{array}{cccc:c}
 {0} & {0} & {0} & {0}  &    \\
 {0} & {0} & {0} & {- a_4^2 }  &     \\
 {0} & {0} & {0} & {0}  &     \mbox{\smash{\huge\textit{O}}}  \\
 {0} & {0} & {0} & {0}  &     \\
 \hdashline
     &      &      &      &   \\
       &     \mbox{\smash{\huge\textit{O}}}    &     &     &  \mbox{\smash{\huge\textit{O}}}      \\
 \end{array}
\right). \\
\end{split}
\end{equation*}
The equality \eqref{FsysKS-1} is equivalent to \eqref{dissipation6-8'}. 
%

Similarly 
we introduce the following $m \times m$ matrix:
\begin{equation*}
K_5 = a_5
\left(
\begin{array}{cccc:cc}
 {0} & {0} & {0} & {0}  &    {0} & \\
 {0} & {0} & {0} & {0}  &    {0} &  \\
 {0} & {0} & {0} & {0}  &    {0} & \mbox{\smash{\huge\textit{O}}}  \\
 {0} & {0} & {0} & {0}  &    {1} &  \\
 \hdashline
 {0} & {0} & {0}  & {-1}  &   {0} &      \\
     &      &      &      &     &  \\
       &     \mbox{\smash{\huge\textit{O}}}    &     &     &   & \mbox{\smash{\huge\textit{O}}}      \\
 \end{array}
\right).
\end{equation*}
Then, we multiply \eqref{Fsys1} by $-i \xi K_5$ and take the inner product with $\hat u$.
Furthermore, taking the real part of the resultant equation, we obtain 
\begin{equation}\label{FsysK3-1}
-\frac{1}{2} \xi \partial_t
\langle i K_5 \hat u, \hat u \rangle 
 + \xi^2 \langle [K_5A_m]^{\rm sy} \hat u, \hat u \rangle
 -  \xi \langle i [K_5L_m]^{\rm asy} \hat u, \hat u \rangle = 0, 
\end{equation}
where 
\begin{equation*}
K_5A_m = 
\left(
\begin{array}{cccc:ccc}
 {0} & {0} & {0} & {0}  &  {0} & {0} & \\
 {0} & {0} & {0} & {0}  &  {0} & {0} &  \\
 {0} & {0} & {0} & {0}  &  {0} &  {0} & \mbox{\smash{\huge\textit{O}}}  \\
 {0} & {0} & {0} & {a_5^2}  &  {0}  & {a_5a_6} &  \\
 \hdashline
 {0} & {0} & {-a_4a_5}  & {0}  &   {-a_5^2} & {0} &   \\
  {0} & {0} & {0}  & {0}  &   {0} & {0} &   \\
     &      &      &      &     &   & \\
       &   \mbox{\smash{\huge\textit{O}}}  &    &     &   &   & \mbox{\smash{\huge\textit{O}}}      \\
 \end{array}
\right), \quad
K_5L_m = 
\left(
\begin{array}{cccc:cc}
 {0} & {0} & {0} & {0}  &  {0} &  \\
 {0} & {0} & {0} & {0}  &  {0} &   \\
 {0} & {0} & {0} & {0}  &  {0} &   \mbox{\smash{\huge\textit{O}}}  \\
 {0} & {0} & {0} & {0}  &  {0}  &   \\
 \hdashline
 {a_5} & {0} & {0}  & {0}  &  {0} &   \\
     &      &      &      &     &  \\
       &   \mbox{\smash{\huge\textit{O}}}  &    &    &   & \mbox{\smash{\huge\textit{O}}}      \\
 \end{array}
\right).
\end{equation*}
The equality \eqref{FsysK3-1} is equivalent to \eqref{dissipation6-9'}.


Based on the energy method of Step 2 in Subsection 2.3, we introduce the following $m \times m$ matrices:
$$
K_{\ell} = a_\ell \hskip-3mm
\bordermatrix*[{(  )}]{
       &           & &0         & 0          &  &  & &  \cr
       & \mbox{\smash{\huge\textit{O}}} 
       & &\vdots &\vdots  &  & \mbox{\smash{\huge\textit{O}}}  & & \cr
       &           & &0         & 0          &  &  & &  \cr
0       & \cdots  &0 &0         & 1          & 0  &\cdots   &0\ \ &  \ell-1 \cr
0       & \cdots  &0 &-1        &0           &0  & \cdots   &0\ \  & \ell \cr
       &           & &0          & 0         &     &  &   & \cr
       & \mbox{\smash{\huge\textit{O}}}& &\vdots  & \vdots &    
       & \mbox{\smash{\huge\textit{O}}}    &   &  \cr
       &           & &0         & 0 &   & & &  \cr
 &     &   &\ell -1 &  \ell      & & &   &
}
$$

\vskip4mm

\noindent
for $\ell = 6, \cdots , m-1$.
Then, we multiply \eqref{Fsys1} by $-i K_{\ell}$ and take the inner product with $\hat u$.
Furthermore, taking the real part of the resultant equation, we obtain 
\begin{equation}\label{FsysKL-1}
-\frac{1}{2} \xi \partial_t
\langle i K_{\ell} \hat u, \hat u \rangle 
 + \xi^2 \langle [K_{\ell} A_m]^{\rm sy} \hat u, \hat u \rangle
=0
\end{equation}
for $\ell = 6, \cdots , m-1$, 
where 
$$
\hspace{-12mm} K_{\ell} A_m =  \hskip-3mm
\bordermatrix*[{(  )}]{
       &           & &0 &0         & 0      &0    &  &  & &  \cr
       & \mbox{\smash{\huge\textit{O}}} 
       & &\vdots &\vdots &\vdots  &\vdots &  & \mbox{\smash{\huge\textit{O}}}  & & \cr
       &           & &0 &0         & 0   &0       &  &  & &  \cr
0       & \cdots  &0 &0 &a_\ell^2         & 0          & a_\ell a_{\ell+1} &0&\cdots   &0\ \   &  \ell-1 \cr
0       & \cdots  &0 & -a_{\ell-1} a_{\ell} &0         &-a_\ell^2           &0  &0& \cdots   &0 \ \  & \ell \cr
       &           & &0 &0          & 0    &0     &     &  &   & \cr
       & \mbox{\smash{\huge\textit{O}}}& &\vdots &\vdots  & \vdots & \vdots &    
       & \mbox{\smash{\huge\textit{O}}}    &   &  \cr
       &           & &0 &0         & 0 &0  &   & & &  \cr
 &     &  &\ell-2  &\ell -1 &  \ell   & \ell+1 & &&   &
}
$$

\vskip4mm

\noindent
Moreover we have
\begin{equation}\label{FsysKM-1}
-\frac{1}{2} \xi \partial_t
\langle i K_{m} \hat u, \hat u \rangle 
 + \xi^2 \langle [K_{m} A_m]^{\rm sy} \hat u, \hat u \rangle
 -  \xi \langle i [K_{m} L_m]^{\rm asy} \hat u, \hat u \rangle = 0,
\end{equation}
where
\begin{equation*}
K_{m}A_m = 
\left(
\begin{array}{cccccc}
       &               &       & {0} & {0}  &  {0}  \\
       &   \mbox{\smash{\huge\textit{O}}}  &       & {\vdots} & {\vdots}  &  {\vdots}   \\
       &               &       & {0}                       & {0}           &  {0}   \\
 {0} & {\cdots} & {0} & {0}                       & {a_m^2}  &  {0}    \\
 {0} & {\cdots} & {0} & {-a_{m-1}a_m}  & {0}  &   {-a_m^2}   \\
 \end{array}
\right), \quad
K_{m}L_m = 
\left(
\begin{array}{cccccc}
       &                     &  &  {0}  \\
       &   \mbox{\smash{\huge\textit{O}}}         &  &  {\vdots}   \\
       &                                           &        &  {0}   \\
 {0} & {\cdots}                        & {0}  &  {a_m \gamma}    \\
 {0} & {\cdots}   & {0}  &   {0}   \\
 \end{array}
\right).
\end{equation*}
The equalities \eqref{FsysKL-1} and  \eqref{FsysKM-1} 
are equivalent to \eqref{dissipation-1'} and \eqref{dissipation-2'}, respectively.


For the rest of this subsection, we construct the desired matrices.
According to the strategy of Step 3 in Subsection 2.2, we first combine \eqref{FsysS-1} and \eqref{FsysK1-1}.
More precisely, multiplying \eqref{FsysS-1}, \eqref{FsysKS-1} and \eqref{FsysK1-1} by $(1+\xi^2)$, $(1+\xi^2)$ and $\delta_1$, respectively, 
and combining the resultant equations, we obtain
\begin{equation*}
\begin{split}
&\frac{1}{2} \partial_t
\big\langle \big\{(1+\xi^2)\mathcal{S}  - i \xi( \delta_1 K_1 + (1+\xi^2) K_4)\big\} \hat u, \hat u \big\rangle  \\
& +(1+\xi^2) \langle [\mathcal{S}L_m]^{\rm sy} \hat u, \hat u \rangle +  \xi^2 \langle [(\delta_1 K_1 + (1+\xi^2) K_4)A_m]^{\rm sy} \hat u, \hat u \rangle \\
 & +  \xi (1+\xi^2) \langle i [\mathcal{S}A_m]^{\rm asy} \hat u, \hat u \rangle  
  -  \xi \langle i [(\delta_1 K_1 + (1+\xi^2) K_4)L_m]^{\rm asy} \hat u, \hat u \rangle= 0.
  \end{split} 
\end{equation*}
Here we define $\mathcal{S} = \tilde{S} + S_4$.
We next multiply \eqref{FsysK3-1} with $\ell = 6$ and the above equation by $(1+\xi^2)^2$ and $\delta_2 \xi^2$, respectively, 
and combining the resultant equations, we obtain
\begin{equation*}
\begin{split}
&\frac{1}{2} \partial_t
\big\langle \big\{\delta_2 \xi^2 ((1+\xi^2)\mathcal{S}  - i \xi( \delta_1 K_1 + (1+\xi^2) K_4)) -\xi (1+\xi^2)^2 K_5 \big\} \hat u, \hat u \big\rangle  \\
& + \delta_2 \xi^2  (1+\xi^2) \langle [\mathcal{S}L_m]^{\rm sy} \hat u, \hat u \rangle  +  \delta_2 \xi^3 (1+\xi^2) \langle i [\mathcal{S}A_m]^{\rm asy} \hat u, \hat u \rangle   \\
&+  \xi^2 \langle [(\delta_2 \xi^2 (\delta_1 K_1 + (1+\xi^2) K_4) + (1+\xi^2)^2K_5)A_m]^{\rm sy} \hat u, \hat u \rangle \\
  &-  \xi \langle i[(\delta_2 \xi^2 (\delta_1 K_1 + (1+\xi^2) K_4) + (1+\xi^2)^2K_5)L_m]^{\rm asy} \hat u, \hat u \rangle= 0.
  \end{split} 
\end{equation*}
Moreover, multiplying \eqref{FsysKL-1} and the above equation by $(1+\xi^2)^3$ and $\delta_3 \xi^2$, respectively, 
and combining the resultant equations, we get
\begin{equation*}
\begin{split}
&\frac{1}{2} \partial_t
\big\langle \big\{ \delta_3 \xi^2 (\delta_2 \xi^2 ((1+\xi^2)\mathcal{S}  - i \xi( \delta_1 K_1 + (1+\xi^2) K_4)) \\
&\hskip30mm -i \xi (1+\xi^2)^2 K_5 ) - i \xi (1+\xi^2)^3K_6\big\} \hat u, \hat u \big\rangle  \\
& + \delta_2 \delta_3 \xi^4  (1+\xi^2) \langle [\mathcal{S}L_m]^{\rm sy} \hat u, \hat u \rangle  
+  \delta_2 \delta_3 \xi^5 (1+\xi^2) \langle i [\mathcal{S}A_m]^{\rm asy} \hat u, \hat u \rangle   \\
&+  \xi^2 \langle [(\delta_3 \xi^2 (\delta_2 \xi^2 (\delta_1 K_1 + (1+\xi^2) K_4) + (1+\xi^2)^2K_5) + (1+\xi^2)^3 K_6)A_m]^{\rm sy} \hat u, \hat u \rangle \\
  &-  \delta_3 \xi^3  \langle i[(\delta_2 \xi^2 (\delta_1 K_1 + (1+\xi^2) K_4) + (1+\xi^2)^2K_5)L_m]^{\rm asy} \hat u, \hat u \rangle= 0.
  \end{split} 
\end{equation*}
Consequently, by the induction argument with respect to $\ell$ in \eqref{FsysKL-1}, we have
\begin{multline}\label{FsysSK1-1}
\frac{1}{2} \partial_t
\Big\langle \Big\{ \prod_{j=2}^{\ell-3} \delta_j \xi^{2(\ell-4)}(1+\xi^2)\mathcal{S}  - i \xi \mathcal{K}_{\ell} \Big\} \hat u, \hat u \Big\rangle  \\
 + \prod_{j=2}^{\ell-3} \delta_j \xi^{2(\ell-4)} (1+\xi^2) \langle [\mathcal{S}L_m]^{\rm sy} \hat u, \hat u \rangle 
 +  \prod_{j=2}^{\ell-3} \delta_j \xi^{2(\ell-4) + 1}  (1+\xi^2) \langle i [\mathcal{S}A_m]^{\rm asy} \hat u, \hat u \rangle   \\
+  \xi^2 \langle [\mathcal{K}_{\ell}A_m]^{\rm sy} \hat u, \hat u \rangle 
 - \prod_{j=3}^{\ell-3} \delta_j \xi^{2(\ell-5) + 1}  \langle i[\mathcal{K}_5 L_m]^{\rm asy} \hat u, \hat u \rangle= 0.
\end{multline}
for $5 \le \ell \le m-1$, where the last term of left hand side is replaced by 
$\xi  \langle i[\mathcal{K}_5 L_m]^{\rm asy} \hat u, \hat u \rangle$ for $\ell = 5$.
Here we define $\mathcal{K}_{\ell}$ as $\mathcal{K}_4 =   \delta_1 K_1 + (1+\xi^2) K_4$ and 
\begin{equation*}
\mathcal{K}_{\ell} 
= \delta_{\ell-3} \xi^2 \mathcal{K}_{\ell -1} +  (1+\xi^2)^{\ell-3} K_{\ell} 
\end{equation*}
for $\ell \ge 5$.
Therefore, we make the combination of \eqref{FsysKM-1} and \eqref{FsysSK1-1} with $\ell = m-1$.
Then we can obtain
\begin{multline}\label{FsysSK2-1}
\frac{1}{2} \partial_t
\Big\langle \Big\{ \prod_{j=2}^{m-4} \delta_j \xi^{2(m-4)}(1+\xi^2)\mathcal{S}  - i \xi \mathcal{K}_{m} \Big\} \hat u, \hat u \Big\rangle  
+ \xi^2 \langle [\mathcal{K}_{m}A_m]^{\rm sy} \hat u, \hat u \rangle \\
 + \prod_{j=2}^{m-4} \delta_j \xi^{2(m-4)} (1+\xi^2) \langle [\mathcal{S}L_m]^{\rm sy} \hat u, \hat u \rangle 
 +  \prod_{j=2}^{m-4} \delta_j \xi^{2(m-4) + 1}  (1+\xi^2) \langle i [\mathcal{S}A_m]^{\rm asy} \hat u, \hat u \rangle   \\
  - \prod_{j=3}^{m-4} \delta_j \xi^{2(m-5) + 1}  \langle i[\mathcal{K}_5 L_m]^{\rm asy} \hat u, \hat u \rangle 
   -  \xi(1+\xi^2)^{m-3}  \langle i [K_{m} L_m]^{\rm asy} \hat u, \hat u \rangle= 0.
\end{multline}
Finally,  multiplying \eqref{FsysSK2-1} by $\delta_{m-3}/(1+\xi^2)^{m-2}$, 
and combining \eqref{eq} and the resultant equations, we can obtain
\begin{multline}\label{FsysSKfinal-1}
\frac{1}{2} \partial_t
\Big\langle \Big[ I + \frac{\delta_{m-3}}{(1+\xi^2)^{m-2}} \Big\{ \prod_{j=2}^{m-4} \delta_j \xi^{2(m-4)}(1+\xi^2)\mathcal{S}  - i \xi \mathcal{K}_{m} \Big\}\Big] \hat u, \hat u \Big\rangle  \\
 + \langle L_m \hat u, \hat u \rangle 
+ \prod_{j=2}^{m-3} \delta_j \frac{\xi^{2(m-4)}}{(1+\xi^2)^{m-3}} \langle [\mathcal{S}L_m]^{\rm sy} \hat u, \hat u \rangle \\
+ \delta_{m-3} \frac{\xi^2}{(1+\xi^2)^{m-2}} \langle [\mathcal{K}_{m}A_m]^{\rm sy} \hat u, \hat u \rangle
 +  \prod_{j=2}^{m-3} \delta_j \frac{\xi^{2(m-4) + 1}}{(1+\xi^2)^{m-3}} \langle i [\mathcal{S}A_m]^{\rm asy} \hat u, \hat u \rangle   \\
  - \prod_{j=3}^{m-3} \delta_j \frac{\xi^{2(m-5) + 1}}{(1+\xi^2)^{m-2}}  \langle i[\mathcal{K}_5 L_m]^{\rm asy} \hat u, \hat u \rangle  
- \delta_{m-3} \frac{\xi}{1+\xi^2}  \langle i [K_{m} L_m]^{\rm asy} \hat u, \hat u \rangle= 0.
\end{multline}
where $I$ denotes an identity matrix.
%
Letting $\delta_1, \cdots, \delta_{m-3}$ suitably small,
then \eqref{FsysSKfinal-1} derives energy estimate \eqref{energy-eq1}.
More precisely, noting that
\begin{equation*}
\begin{split}
\mathcal{K}_{m} 
& =\prod_{j=2}^{m-3} \delta_{j} \xi^{2(m-4)}(   \delta_1 K_1 + (1+\xi^2) K_4)+  (1+\xi^2)^{m-3} K_{m} \\  
&\qquad  +  \sum_{k=3}^{m-3} \prod_{j=k}^{m-3} \delta_{j} \xi^{2(m-k-2)}(1+\xi^2)^{k-1} K_{k+2} \\
\end{split}
\end{equation*}
for $m \ge 6$, 
we can estimate the dissipation terms as
\begin{multline}\label{1estSK}
  \langle L_m \hat u, \hat u \rangle 
+ \prod_{j=2}^{m-3} \delta_j \frac{\xi^{2(m-4)}}{(1+\xi^2)^{m-3}} \langle [\mathcal{S}L_m]^{\rm sy} \hat u, \hat u \rangle \\
+ \delta_{m-3} \frac{\xi^2}{(1+\xi^2)^{m-2}} \langle [\mathcal{K}_{m}A_m]^{\rm sy} \hat u, \hat u \rangle \\
 \ge c \Big\{
\frac{\xi^{2(m-4)}}{(1+\xi^2)^{m-3}} |\hat u_{1}|^2  +  \frac{\xi^{2(m-3)}}{(1+\xi^2)^{m-2}} |\hat u_{2}|^2 
 + \sum_{j=3}^{m}  \frac{\xi^{2(m-j)}}{(1+\xi^2)^{m-j}} |\hat u_{j}|^2  \Big\} 
\end{multline}
for suitably small $\delta_1, \cdots, \delta_{m-3}$.
Consequently we conclude that
our desired symmetric matrix $S$  
and skew-symmetric matrix $K$
are described as
\begin{equation*}
S = \frac{\xi^{2(m-4)}}{(1+\xi^2)^{m-3}} \mathcal{S},
\qquad
K = \frac{\xi^2}{(1+\xi^2)^{m-2}}\mathcal{K}_m.
\end{equation*}




\section{Model II}

\subsection{Main result II}
In this section, we treat the Cauchy problem \eqref{sys1}, \eqref{ID} with
\begin{equation}\label{2Tmat}
\begin{split}
A_m &=
\left(
\begin{array}{cc:cc:ccccc:cc}
 {0} & {1} & {0} & {0} &       &       &      &       &       &       &       \\
 {1} & {0} & {0} & {0} &       &       &      &       &       &       &       \\
 \hdashline 
 {0} & {0} & {0} & {a_4} & {0} &     
     &       &     &         &      &       \\
 {0} & {0} & {a_4} & {0} & {0} & {0} &     &   &  \mbox{\smash{\huge\textit{O}}}   &       &      \\
 \hdashline
       &      & {0} & {0} & {0} & {a_6} &    &     &       &       &      \\
       &       &       & {0} & {a_6} & {0} &   &   &     &     &      \\
       &       &       &       &      &    &   \ddots  &    & &    &       \\
       &       &       &       &      &       &       &  0 & {a_{m-2}}  &    &      \\
       &       &      \mbox{\smash{\huge\textit{O}}}  &     &    
       &       &       & {a_{m-2}} & {0} & {0} & {0} \\
\hdashline
       &       &       &       &       &       &       & {0} & {0} & {0} & {a_{m}} \\
       &       &       &       &       &       &       &       & {0} & {a_{m}} & {0} \\
\end{array}
\right),  \\
L_m &=
\left(
\begin{array}{c:cc:ccccc:cc:c}
 {0} & {0} & {0} & {0} &       &       &      &       &       &       &       \\
  \hdashline  
 {0} & {\gamma} & {1} & {0} &       &       &      &       &       &       &       \\
 {0} & {-1} & {0} & {0} & {0} &     
     &       &   &  \mbox{\smash{\huge\textit{O}}}     &      &       \\
\hdashline
 {0} & {0} & {0} & {0} & {a_5} &      &     &       &       &       &      \\
        &      & {0} & {-a_5} & {0} &      &      &     &       &       &      \\
       &       &       &       &      & \ddots &     &     &     &     &      \\
       &       &       &       &      &      &  0  &  {a_{m-3}} &  0  &    &       \\
       &       &       &       &      &       &   {-a_{m-3}}   &   0 & 0  & 0   &      \\
\hdashline
       &       &   &     &     &       &  {0}     & {0} & {0} & {a_{m-1}} & {0} \\
       &       &         \mbox{\smash{\huge\textit{O}}}      &       &       &       &       & {0} & {-a_{m-1}} & {0} & {0} \\
\hdashline
       &       &       &       &       &       &       &       & {0} & {0} & {0} \\
\end{array}
\right),  
\end{split}
\end{equation}
where integer $m\geq 4$ is even, $\gamma > 0$, and all elements $a_j$ $(4\leq j\leq m)$ are nonzero.
We note that the system \eqref{sys1} with \eqref{2Tmat} for $m=4$ is the Timoshenko system (cf.~\cite{IHK08,IK08}). 
For this problem, we can derive the following decay structure.

\begin{thm}
\label{thm2}
%
The Fourier image $\hat u$ of the solution $u$ to the Cauchy problem
\eqref{sys1}-\eqref{ID} with \eqref{2Tmat} satisfies the pointwise estimate:
\begin{equation}\label{point2}
|\hat u(t,\xi)| \le Ce^{-c \lambda(\xi)t}|\hat u_0(\xi)|,
\end{equation}
where $\lambda(\xi) := \xi^{3m-10}/(1+\xi^2)^{2(m-3)}$.
Furthermore, let $s \ge 0$ be an integer and suppose that the initial data
$u_0$ belong to $H^s \cap L^1$.
Then the solution $u$ satisfies the decay estimate: 
\begin{equation}
\notag
\|\partial_x^{k} u(t)\|_{L^2} \le C(1+t)^{-\frac{1}{3m-10}(\frac{1}{2}+k)}\|u_0\|_{L^1} + C(1+t)^{-\frac{\ell}{m-2}}\|\partial_x^{k+\ell} u_0\|_{L^2}
\end{equation}
for $k+\ell \le s$. Here $C$ and $c$ are positive constants.
\end{thm}



\subsection{Energy method in the case $m=6$}
Ide-Hramoto-Kawashima \cite{IHK08} and Ide-Kawashima \cite{IK08} had already obtained the desired estimates in the case $m=4$.
Thus we consider the case $m=6$ in this subsection, which can shed light on the proof of the general case $m\geq 6$ to be given by Section 3.3.
Then we rewrite the system \eqref{sys1} with \eqref{2Tmat} as follows.
\begin{equation}\label{2equations6}
\begin{split}
&\partial_t \hat u_1 + i \xi \hat u_2 = 0, \\
&\partial_t \hat u_2 + i \xi \hat u_1 + \gamma \hat u_2 + \hat u_3  = 0, \\
&\partial_t \hat u_3 + i \xi a_4 \hat u_4 - \hat u_2 = 0, \\
&\partial_t \hat u_4 + i \xi a_{4} \hat u_{3} + a_{5} \hat u_{5} = 0, \\
&\partial_t \hat u_5 + i \xi a_{6} \hat u_{6}  - a_{5} \hat u_{4} = 0,  \\
&\partial_t \hat u_6 + i \xi a_{6} \hat u_{5} = 0. 
\end{split}
\end{equation}

\medskip

\noindent
{\bf Step 1.}\ \
We first derive the basic energy equality for the system \eqref{2equations6} in the Fourier space.
We multiply the all equations of \eqref{2equations6} by 
$\bar{\hat{u}} = (\bar{\hat{u}}_1, \bar{\hat{u}}_2,\bar{\hat{u}}_3,\bar{\hat{u}}_4, \bar{\hat{u}}_5,\bar{\hat{u}}_6)^T$, respectively,
and combine the resultant equations.
Furthermore, taking the real part for the resultant equality, we arrive at the basic energy equality
\begin{equation}\label{2eq6}
\frac{1}{2} \partial_t |\hat u|^2
+ \gamma |\hat u_2|^2 = 0.
\end{equation}
Next we create the dissipation terms by the following two steps.

\medskip

\noindent
{\bf Step 2.}\ \
We multiply the first and second equations in \eqref{2equations6} 
by $i \xi \bar{\hat{u}}_{2}$ and $ - i \xi \bar{\hat{u}}_{1}$, respectively.
Then, combining the resultant equations and taking the real part, we have
\begin{equation}\label{2dissipation6-4'}
 \xi \partial_t \Re (i \hat u_{1} \bar{\hat{u}}_2)
+ \xi^2 ( |\hat u_{1}|^2 - |\hat u_{2}|^2 )  
+ \gamma  \xi \,  \Re (i \hat{u}_{1} \bar{\hat{u}}_{2}) 
+ \xi \, \Re (i \hat{u}_{1} \bar{\hat{u}}_{3}) =  0.
\end{equation}
Next, we combine the fourth and sixth equations in \eqref{2equations6}, obtaining
\begin{equation*}
\partial_t (\xi a_6\hat u_4  +  i a_5 \hat{u}_6) + i \xi^2 a_{4} a_6 \hat u_{3}  = 0.
\end{equation*}
Then multiplying  the first equation in \eqref{2equations6} and the resultant equation by $\xi a_6\bar{\hat u}_4  - i a_5 \bar{\hat{u}}_6$ and $\bar{\hat{u}}_1$, 
and combining the resultant equations and taking the real part, we obtain
\begin{multline}\label{2dissipation6-2}
\partial_t \big\{ a_6 \xi \Re (\hat{u}_1 \bar{\hat{u}}_4) -  a_5 \Re (i \hat{u}_1 \bar{\hat{u}}_6) \big\} \\
- a_4 a_6 \xi^2  \, \Re (i \hat{u}_{1} \bar{\hat{u}}_3) +  a_6 \xi^2  \, \Re ( i \hat{u}_{2} \bar{\hat{u}}_4)  + a_5 \xi \, \Re (\hat{u}_{2} \bar{\hat{u}}_{6})  = 0.
\end{multline}
To eliminate $\Re (i \hat{u}_{1} \bar{\hat{u}}_3)$, we multiply \eqref{2dissipation6-4'} and \eqref{2dissipation6-2} by $a_4^2 a_6^2 \xi^2$ 
and $a_4 a_6 \xi$, add the resultant equations.
Then this yields
\begin{multline}\label{2dissipation6-3}
a_4 a_6 \xi \partial_t E_1^{(6)} 
+ a_4^2 a_6^2 \xi^4 (|\hat u_{1}|^2 - |\hat u_{2}|^2 ) \\
+ a_4 a_6^2 \xi^3  \, \Re ( i \hat{u}_{2} \bar{\hat{u}}_4)  + a_4 a_5 a_6 \xi^2 \, \Re (\hat{u}_{2} \bar{\hat{u}}_{6})  + \gamma a_4^2 a_6^2 \xi^3 \, \Re (i \hat{u}_{1} \bar{\hat{u}}_{2}) = 0,
\end{multline}
where $E_1^{(6)} = a_6 \xi \Re (\hat{u}_1 \bar{\hat{u}}_4) -  a_5 \Re (i \hat{u}_1 \bar{\hat{u}}_6)  + a_4 a_6 \xi^2 \Re (i\hat u_{1} \bar{\hat{u}}_2)$.

On the other hand, we multiply the second and third equations in \eqref{2equations6} 
by $\bar{\hat{u}}_{3}$ and $\bar{\hat{u}}_{2}$, respectively.
Then, combining the resultant equations and taking the real part, we have
\begin{equation}\label{2dissipation6-1'}
\partial_t  \Re (\hat u_{2} \bar{\hat{u}}_3) +   |\hat u_{3}|^2 - |\hat u_{2}|^2  
+ \xi \, \Re (i \hat{u}_{1} \bar{\hat{u}}_3) - a_{4} \xi \, \Re (i \hat{u}_{2} \bar{\hat{u}}_{4}) 
+ \gamma \, \Re ( \hat{u}_{2} \bar{\hat{u}}_{3}) =  0.
\end{equation}
By the Young inequality, the equation \eqref{2dissipation6-1'} is estimated as
\begin{equation}\label{2dissipation6-1}
\partial_t E_3 + \frac{1}{2} |\hat u_{3}|^2 
\le \xi^2 |\hat{u}_{1}|^2 +  (1+\gamma^2) |\hat u_{2}|^2  + a_{4} \xi \, \Re (i \hat{u}_{2} \bar{\hat{u}}_{4}),
\end{equation}
where $E_3 =  \Re (\hat u_{2} \bar{\hat{u}}_3)$.

Furthermore,
we multiply the third equation and fourth equation of \eqref{2equations6} 
by  $- i \xi a_4 \bar{\hat{u}}_{4}$ and $ i \xi a_4 \bar{\hat{u}}_{3}$, respectively.
Then, combining the resultant equations and taking the real part, we have
\begin{equation}\label{2dissipation6-5'}
- a_4 \xi \partial_t \Re (i \hat u_{3} \bar{\hat{u}}_4) 
+   a_4^2 \xi^2( |\hat u_{4}|^2 - |\hat u_{3}|^2 )  
+ a_4 \xi \, \Re (i \hat{u}_{2} \bar{\hat{u}}_4) - a_4 a_{5} \xi \, \Re (i \hat{u}_{3} \bar{\hat{u}}_{5}) =  0.
\end{equation}
By the Young inequality, the above equation is estimated as
\begin{equation}\label{2dissipation6-5}
\xi \partial_t E_4
+  \frac{1}{2}  a_4^2 \xi^2 |\hat u_{4}|^2 
\le \frac{1}{2} |\hat{u}_{2} |^2  + a_4^2 \xi^2 |\hat u_{3}|^2
+ a_4 a_{5} \xi \, \Re (i \hat{u}_{3} \bar{\hat{u}}_{5}),
\end{equation}
where $E_4 = - a_4 \Re (i \hat u_{3} \bar{\hat{u}}_4)$.

We multiply the fourth equation and fifth equation in \eqref{2equations6} 
by $a_{5} \bar{\hat{u}}_{5}$ and $ a_{5} \bar{\hat{u}}_{4}$, respectively.
Then, combining the resultant equations and taking the real part, we have
\begin{equation*}
a_{5} \partial_t \Re (\hat u_{4} \bar{\hat{u}}_{5}) 
+  a_{5}^2 ( |\hat u_{5}|^2 - |\hat u_{4}|^2)  \\
+ a_{4} a_{5} \xi \, \Re (i\hat{u}_{3} \bar{\hat{u}}_{5}) 
- a_{5} a_{6} \xi \, \Re (i \hat{u}_{4} \bar{\hat{u}}_{6}) = 0.
\end{equation*}
By using Young inequality, we obtain 
\begin{equation}\label{2dissipation6-6}
 \partial_t E_5 
+  \frac{1}{2} a_{5}^2 |\hat u_{5}|^2 \\
 \le  a_{5}^2 |\hat u_{4}|^2
+ \frac{1}{2} a_{5}^2 \xi^2 |\hat{u}_{3}|^2 
+ a_{5} a_{6} \xi \, \Re (i \hat{u}_{4} \bar{\hat{u}}_{6}).
\end{equation}
where $E_5 = a_{5} \partial_t \Re (\hat u_{4} \bar{\hat{u}}_{5})$.

Moreover, we multiply the last equation and the fifth equation in \eqref{2equations6} 
by  $ i \xi a_6 \bar{\hat{u}}_{5}$ and $- i \xi a_6 \bar{\hat{u}}_{6}$, respectively.
Then, combining the resultant equations and taking the real part, we have
\begin{equation*}
- a_6 \xi \partial_t \Re (i \hat u_{5} \bar{\hat{u}}_6) 
+   a_6^2 \xi^2( |\hat u_{6}|^2 - |\hat u_{5}|^2 )  
+ a_{5}a_6 \xi \, \Re (i \hat{u}_{4} \bar{\hat{u}}_6) =  0.
\end{equation*}
Using Young inequality, this yields
\begin{equation}\label{2dissipation6-8}
 \xi \partial_t E_6
+ \frac{1}{2}  a_6^2 \xi^2 |\hat u_{6}|^2 
\le a_6^2 \xi^2 |\hat u_{5}|^2 
+ \frac{1}{2} a_{5}^2  |\hat{u}_{4}|^2,
\end{equation}
where $E_6 = - a_6 \Re(i \hat u_{5} \bar{\hat{u}}_6)$.

\medskip

\noindent
{\bf Step 3.}\ \
In this step, we sum up the energy inequalities 
and derive the desired energy inequality.
For this purpose, 
we first multiply \eqref{2dissipation6-6} and \eqref{2dissipation6-8} 
by $\xi^2$ and $\beta_1$, respectively.
Then we combine the resultant equation, obtaining
\begin{multline*}
 \partial_t 
\big\{ \xi^2 E_5 +\beta_1  \xi E_6 \big\} + \frac{1}{2} \beta_1 a_6^2 \xi^2 |\hat u_{6}|^2 
+ \Big( \frac{1}{2} a_{5}^2  - \beta_1 a_6^2 \Big) \xi^2 |\hat u_{5}|^2 \\
 \le \Big( \frac{1}{2} \beta_1 + \xi^2 \Big) a_{5}^2  |\hat u_{4}|^2
+ \frac{1}{2} a_{5}^2 \xi^4 |\hat{u}_{3}|^2 
+ |a_{5}| |a_{6}| |\xi|^3 |\hat{u}_{4}| |\hat{u}_{6}|.
\end{multline*}
Letting $\beta_1$ suitably small and using Young inequality, we get
\begin{equation*}
\begin{split}
& \partial_t \big\{ \xi^2 E_5 +\beta_1  \xi E_6 \big\}
 + c \xi^2 (|\hat u_{5}|^2 + |\hat u_{6}|^2) 
 \le C (1 + \xi^2)^2  |\hat u_{4}|^2
+ \frac{1}{2} a_{5}^2 \xi^4 |\hat{u}_{3}|^2.
\end{split}
\end{equation*}
Moreover, combining the above estimate and 
\eqref{2dissipation6-6}, we get 
\begin{multline}\label{2eneq6-1}
 \partial_t \big\{ (1+\xi^2) E_5 +\beta_1  \xi E_6 \big\}
 + c (1+ \xi^2) |\hat u_{5}|^2 + c \xi^2 |\hat u_{6}|^2 \\
 \le C (1 + \xi^2)^2  |\hat u_{4}|^2
+ \frac{1}{2} a_{5}^2 \xi^2(1+\xi^2) |\hat{u}_{3}|^2.
\end{multline}
Second, we multiply \eqref{2dissipation6-5} and \eqref{2eneq6-1} by $(1+\xi^2)^2$ and $\beta_2 \xi^2$, respectively,
and combine the resultant equations.
Then we obtain
\begin{multline*}
 \partial_t \big\{ \beta_2 \xi^2((1+\xi^2) E_5 +\beta_1  \xi E_6) + \xi (1+\xi^2)^2 E_4 \big\} \\
 + \beta_2 c \xi^2 (1+ \xi^2) |\hat u_{5}|^2 + \beta_2 c \xi^4 |\hat u_{6}|^2 
 + \Big( \frac{1}{2} a_{4}^2 - \beta_2 C \Big) \xi^2 (1 + \xi^2)^2  |\hat u_{4}|^2  \\
\le C(1+\xi^2)^2 |\hat{u}_{2} |^2 + C \xi^2(1+\xi^2)^2 |\hat{u}_{3}|^2 +  C \xi (1+\xi^2)^2 \, \Re (i \hat{u}_{3} \bar{\hat{u}}_{5}).
\end{multline*}
Letting $\beta_2$ suitably small and using Young inequality, we get
\begin{multline}\label{2eneq6-2}
 \partial_t \big\{ \beta_2 \xi^2((1+\xi^2) E_5 +\beta_1  \xi E_6) + \xi (1+\xi^2)^2 E_4 \big\} 
 +  c \xi^2 (1+ \xi^2) |\hat u_{5}|^2 \\
 +  c \xi^4 |\hat u_{6}|^2 + c \xi^2 (1 + \xi^2)^2  |\hat u_{4}|^2  
\le C(1+\xi^2)^2 |\hat{u}_{2} |^2 + C (1+\xi^2)^3 |\hat{u}_{3}|^2.
\end{multline}
Third, we multiply \eqref{2dissipation6-1} and \eqref{2eneq6-2} by $(1+\xi^2)^3$ and $\beta_3$
and combine the resultant equations.
Then we obtain
\begin{multline*}
 \partial_t \big\{ \beta_3 (\beta_2 \xi^2((1+\xi^2) E_5 +\beta_1  \xi E_6) + \xi (1+\xi^2)^2 E_4)  + (1+\xi^2)^3 E_3\big\} 
 + \beta_3 c \xi^4 |\hat u_{6}|^2  \\
+  \beta_3 c \xi^2 (1+ \xi^2) |\hat u_{5}|^2 + \beta_3 c \xi^2 (1 + \xi^2)^2  |\hat u_{4}|^2 
 + \Big(\frac{1}{2} - \beta_3 C \Big)(1+\xi^2)^3 |\hat{u}_{3}|^2 \\
 \le C(1+\xi^2)^3 |\hat{u}_{2} |^2 +  \xi^2 (1+\xi^2)^3 |\hat{u}_{1}|^2 +  a_4 \xi (1+\xi^2)^3 \, \Re ( i \hat{u}_{2} \bar{\hat{u}}_4).
\end{multline*}
Therefore, letting $\beta_3$ suitably small and using Young inequality, we get
\begin{multline}\label{2eneq6-3}
 \partial_t \big\{ \beta_3 (\beta_2 \xi^2((1+\xi^2) E_5 +\beta_1  \xi E_6) + \xi (1+\xi^2)^2 E_4)  +  (1+\xi^2)^3 E_3 \big\} \\
 + c \xi^4 |\hat u_{6}|^2 +  c \xi^2 (1+ \xi^2) |\hat u_{5}|^2 + c \xi^2 (1 + \xi^2)^2  |\hat u_{4}|^2 + c (1+\xi^2)^3 |\hat{u}_{3}|^2 \\
 \le C (1+\xi^2)^4 |\hat{u}_{2} |^2 +  \xi^2 (1+\xi^2)^3 |\hat{u}_{1}|^2.
\end{multline}
Fourth, we multiply \eqref{2dissipation6-3} and \eqref{2eneq6-3} by $(1+\xi^2)^3$ and $ \beta_{4} \xi^{2}$, respectively, 
and combine the resultant equalities. 
Moreover, letting $\beta_4$ suitably small and using Young inequality, then we obtain
\begin{multline}\label{2eneq6-4}
 \partial_t \tilde{E}
 + c \xi^6 |\hat u_{6}|^2 +  c \xi^4 (1+ \xi^2) |\hat u_{5}|^2 + c \xi^4 (1 + \xi^2)^2  |\hat u_{4}|^2+ c \xi^2 (1+\xi^2)^3 |\hat{u}_{3}|^2 \\
+ c \xi^4 (1+\xi^2)^3 |\hat{u}_{1}|^2 
 \le C (1+\xi^2)^5 |\hat{u}_{2} |^2 + a_4 a_5 a_6 \xi^2 (1+\xi^2)^3\, \Re (\hat{u}_{2} \bar{\hat{u}}_{6}),
\end{multline}
where we have defined
\begin{multline}
\notag
\tilde{E} = \beta_4 \xi^2 (\beta_3 (\beta_2 \xi^2((1+\xi^2) E_5 +\beta_1  \xi E_6) + \xi (1+\xi^2)^2 E_4)  +  (1+\xi^2)^3 E_3)\\
+ a_4 a_6 \xi (1+\xi^2)^3 E_1^{(6)}.
\end{multline}
Moreover, to estimate $\Re (\hat{u}_{2} \bar{\hat{u}}_{6})$, we multiply \eqref{2eneq6-4} by $\xi^2$ and use Young inequality again.
Then this yields
\begin{multline}\label{2eneq6-4'}
 \xi^2 \partial_t \tilde{E}
 + c \xi^8 |\hat u_{6}|^2 +  c \xi^6 (1+ \xi^2) |\hat u_{5}|^2 + c \xi^6 (1 + \xi^2)^2  |\hat u_{4}|^2\\
 + c \xi^4 (1+\xi^2)^3 |\hat{u}_{3}|^2 
+ c \xi^6 (1+\xi^2)^3 |\hat{u}_{1}|^2 
 \le C (1+\xi^2)^6 |\hat{u}_{2} |^2.
 \end{multline}
Finally, multiplying the basic energy \eqref{2eq6} and  \eqref{2eneq6-4'} by $(1+\xi^2)^6$ and  $\beta_{5}$, respectively,
combining the resultant equations and letting $\beta_{5}$ suitably small,
then this yields
\begin{multline}\label{2eneq6-5}
\partial _t \Big\{ \frac{1}{2}(1+\xi^2)^{6}|\hat{u}|^2 + \beta_{5}\xi^2  \tilde{E} \Big\} 
  +    c \xi^{6} (1 +  \xi^2)^3 |\hat u_{1}|^2  + c (1+ \xi^2)^{6}  |\hat u_{2}|^2   \\
+    c \xi^{4} (1 +  \xi^2)^3 |\hat u_{3}|^2 + c \xi^{6} (1 +  \xi^2)^2 |\hat u_{4}|^2  
  +    c \xi^{6} (1 +  \xi^2) |\hat u_{5}|^2 + c \xi^{8} |\hat u_{6}|^2  \le 0.
\end{multline}
Thus, integrating the above estimate with respect to $t$,
we obtain the following energy estimate
\begin{multline}
\notag
 |\hat{u}(t,\xi)|^2 
+ \int^t_0   \Big\{
\frac{\xi^{6}}{(1+\xi^2)^{3}} |\hat u_{1}|^2 + |\hat u_{2}|^2 
+\frac{\xi^{4}}{(1+\xi^2)^{3}} |\hat u_{3}|^2  
+  \frac{\xi^{6}}{(1+\xi^2)^{4}} |\hat u_{4}|^2 \\
+ \frac{\xi^{6}}{(1+\xi^2)^{5}} |\hat u_{5}|^2  
+  \frac{\xi^{8}}{(1+\xi^2)^{6}} |\hat u_{6}|^2 
  \Big\} d\tau \le  C|\hat{u}(0,\xi)|^2.
\end{multline}
Here we have used the following inequality
\begin{equation}\label{2eneq6-7}
c |\hat{u}|^2 \le 
\frac{1}{2} |\hat{u}|^2  +   \frac{\beta_{5}\xi^2 }{(1+\xi^2)^{6}} \tilde{E}
\le C |\hat{u}|^2
\end{equation}
for suitably small $\beta_{5}$.
Furthermore the estimate \eqref{2eneq6-5} with \eqref{2eneq6-7} give us the following pointwise estimate
\begin{equation}
\notag
|\hat{u}(t,\xi)|
\le C e^{- c \lambda(\xi)t} |\hat{u}(0,\xi)|,
\qquad
\lambda(\xi) = 
\frac{\xi^{8}}{(1+\xi^2)^{6}}.
\end{equation}
This therefore proves \eqref{point2} in the case $m=6$ for Theorem \ref{thm2}.


\subsection{Energy method for model II}

Inspired by the concrete calculation in Subsection 3.2, we consider the more general situation $m\geq 6$.
Then we rewrite our system \eqref{Fsys1} with \eqref{2Tmat} as follows:
\begin{equation}\label{2equations}
\begin{split}
&\partial_t \hat u_1 + i \xi \hat u_2 = 0, \\
&\partial_t \hat u_2 + i \xi \hat u_1 + \gamma \hat u_2 + \hat u_3  = 0, \\
&\partial_t \hat u_3 + i \xi a_4 \hat u_4 - \hat u_2 = 0, \\
&\partial_t \hat u_j + i \xi a_{j} \hat u_{j-1} + a_{j+1} \hat u_{j+1} = 0, \qquad j = 4, 6, \cdots, m-2, \ {\rm (for \ even)}  \\
&\partial_t \hat u_j + i \xi a_{j+1} \hat u_{j+1}  - a_{j} \hat u_{j-1} = 0, \qquad j = 5, 7, \cdots, m-1, \ {\rm (for \ odd)} \\
&\partial_t \hat u_m + i \xi a_{m} \hat u_{m-1} = 0. 
\end{split}
\end{equation}

\medskip

\noindent
{\bf Step 1.}\ \
We first derive the basic energy equality for the system \eqref{Fsys1} in the Fourier space.
Taking the inner product of \eqref{Fsys1} with $\hat{u}$, we have
\begin{equation*}
	\langle \hat u_t, \hat{u} \rangle
	+ i \xi \langle A_m \hat u, \hat u \rangle
	+ \langle L_m \hat u, \hat u \rangle = 0.
\end{equation*}
Taking the real part, we get the basic energy equality
\begin{equation*}
\frac{1}{2} \partial_t |\hat u|^2
+ \langle L_m \hat u, \hat u \rangle = 0, 
\end{equation*}
and hence
\begin{equation}\label{2eq}
\frac{1}{2} \partial_t |\hat u|^2
+ \gamma |\hat u_2|^2 = 0.
\end{equation}
Next we create the dissipation terms by the following three steps.

\medskip

\noindent
{\bf Step 2.}\ \
We note that we had already derived some useful equations in Subsection 3.2.
Indeed the equations \eqref{2dissipation6-4'}, \eqref{2dissipation6-1}, \eqref{2dissipation6-5} 
and \eqref{2dissipation6-6} are valid for our general problem.
Therefore we adopt these equations in this subsection.

To eliminate $\Re(i \hat{u}_{1} \bar{\hat{u}}_3) $ in \eqref{2dissipation6-4'}, 
we first prepare the useful equation.
We combine the fourth equations with $j= 4, \cdots, 2 \ell$ in \eqref{2equations} inductively.
Then we obtain
\begin{equation}\label{2mathcalU}
\partial_t \mathcal{U}_{2\ell} + i \xi (-i \xi)^{\ell - 2} \prod_{j=2}^{\ell} a_{2j} \hat u_{3} +  \prod_{j=2}^{\ell} a_{2j+1} \hat u_{2\ell + 1} = 0,
\end{equation}
for $4 \le 2 \ell \le m-2$, 
where we have defined $\mathcal{U}_4 = \hat{u}_4$ and 
\begin{equation}
\notag
\begin{split}
\mathcal{U}_{2 \ell} &= -  i \xi a_{2\ell} \mathcal{U}_{2\ell-2} + \prod_{j=2}^{\ell-1} a_{2j+1} \hat u_{2\ell}.  \\
  \end{split}
 \end{equation}
Moreover, combining the last equation in \eqref{2equations} and \eqref{2mathcalU}, this yields
\begin{equation}\label{2mathcalU-2}
i^{m/2} \partial_t \mathcal{U}_{m} - i \xi^{m/2 - 1} \prod_{j=2}^{m/2} a_{2j} \hat u_{3}  = 0.
\end{equation}
Multiplyig \eqref{2mathcalU-2} by $- \bar{\hat{u}}_1$ and the first equation in \eqref{2equations} by $- \overline{i^{m/2}\mathcal{U}_m}$, 
combining the resultant equations and taking the real part, we obtain
\begin{equation}\label{2estU}
- \partial_t \Re ( i^{m/2}\mathcal{U}_m \bar{\hat{u}}_1 ) 
-  \prod_{j=2}^{m/2} a_{2j}\xi^{m/2-1}  \Re (i \hat{u}_{1} \bar{\hat{u}}_3)  + \xi \, \Re (i^{m/2+1}\mathcal{U}_m \bar{\hat{u}}_{2})  = 0.
\end{equation}
In order to eliminate $\Re (i \hat{u}_{1} \bar{\hat{u}}_3)$,
we multiply \eqref{2dissipation6-4'} by $\prod_{j=2}^{m/2} a_{2j} \xi^{m/2-2} $ and 
combine the resultant equation and \eqref{2estU}.
Then we obtain
\begin{multline}\label{2estU2}
 \partial_t E_1^{(m)}
+  \prod_{j=2}^{m/2} a_{2j} \xi^{m/2}   ( |\hat u_{1}|^2 - |\hat u_{2}|^2 ) \\
+ \gamma \prod_{j=2}^{m/2} a_{2j}\xi^{m/2-1} \Re (i \hat{u}_{1} \bar{\hat{u}}_{2}) 
+ \xi \, \Re (i^{m/2+1}\mathcal{U}_m \bar{\hat{u}}_{2})  = 0,
\end{multline}
where we have defined 
$$
E_1^{(m)} = \prod_{j=2}^{m/2} a_{2j} \xi^{m/2-1}  \Re (i \hat u_{1} \bar{\hat{u}}_2) -\Re ( i^{m/2}\mathcal{U}_m \bar{\hat{u}}_1 ).
$$

For $\ell = 4,6, \cdots , m-2$,
we multiply the fourth equation and fifth equation with $j=\ell$ and $j=\ell+1$ in \eqref{2equations} 
by $a_{\ell + 1} \bar{\hat{u}}_{\ell+1}$ and $ a_{\ell+1} \bar{\hat{u}}_{\ell}$, respectively.
Then, combining the resultant equations and taking the real part, we have
\begin{multline}\label{2dissipation-6'}
a_{\ell+1} \partial_t \Re (\hat u_{\ell} \bar{\hat{u}}_{\ell+1}) 
+  a_{\ell+1}^2 ( |\hat u_{\ell+1}|^2 - |\hat u_{\ell}|^2)  \\
+ a_{\ell} a_{\ell+1} \xi \, \Re (i\hat{u}_{\ell-1} \bar{\hat{u}}_{\ell+1}) 
- a_{\ell+1} a_{\ell+2} \xi \, \Re (i \hat{u}_{\ell} \bar{\hat{u}}_{\ell+2}) = 0.
\end{multline}
By using Young inequality, we obtain 
\begin{equation}\label{2dissipation-6}
\partial_t E_{\ell+1} 
+  \frac{1}{2} a_{\ell+1}^2 |\hat u_{\ell+1}|^2 
 \le  a_{\ell+1}^2 |\hat u_{\ell}|^2
+ \frac{1}{2} a_{\ell+1}^2 \xi^2 |\hat{u}_{\ell-1}|^2 
+ a_{\ell+1} a_{\ell+2} \xi \, \Re(i \hat{u}_{\ell} \bar{\hat{u}}_{\ell+2}).
\end{equation}
where $E_{\ell + 1} = a_{\ell+1} \Re(\hat u_{\ell} \bar{\hat{u}}_{\ell+1})$.

On the other hand, for $\ell = 4, \cdots , m-4$,
we multiply the fourth and fifth equations with $j=\ell+2$ and $j = \ell +1$ in \eqref{2equations} 
by  $ i \xi a_{\ell+2} \bar{\hat{u}}_{\ell+1}$ and $- i \xi a_{\ell+2} \bar{\hat{u}}_{\ell+2}$, respectively.
Then, combining the resultant equations and taking the real part, we have
\begin{equation}\label{2dissipation-7'}
	\begin{split}
&- a_{\ell+2} \xi \partial_t \Re(i \hat u_{\ell+1} \bar{\hat{u}}_{\ell+2}) 
+   a_{\ell+2}^2 \xi^2( |\hat u_{\ell+2}|^2 - |\hat u_{\ell+1}|^2 )  \\
&\hskip20mm + a_{\ell+1}a_{\ell+2} \xi \, \Re (i \hat{u}_{\ell} \bar{\hat{u}}_{\ell+2}) 
- a_{\ell+2} a_{\ell+3} \xi \, \Re (i \hat{u}_{\ell+1} \bar{\hat{u}}_{\ell+3}) =  0.
	\end{split}
\end{equation}
Here, by using Young inequality, we obtain 
\begin{multline}\label{2dissipation-7}
 \xi \partial_t E_{\ell+2} 
+  \frac{1}{2} a_{\ell+2}^2 \xi^2 |\hat u_{\ell+2}|^2 \\
\hskip10mm \le  a_{\ell+2}^2 \xi^2 |\hat u_{\ell+1}|^2 
+ \frac{1}{2} a_{\ell+1}^2  | \hat{u}_{\ell}|^2  
+ a_{\ell+2} a_{\ell+3} \xi \, \Re (i \hat{u}_{\ell+1} \bar{\hat{u}}_{\ell+3}),
\end{multline}
where
$E_{\ell+2} = - a_{\ell+2} \Re (i \hat u_{\ell+1} \bar{\hat{u}}_{\ell+2})$.

Moreover, we multiply the last equation and the fifth equation with $j=m-1$ in \eqref{2equations} 
by  $ i \xi a_m \bar{\hat{u}}_{m-1}$ and $- i \xi a_m \bar{\hat{u}}_{m}$, respectively.
Then, combining the resultant equations and taking the real part, we have
\begin{equation}\label{2dissipation-8'}
- a_m \xi \partial_t \Re (i \hat u_{m-1} \bar{\hat{u}}_m) 
+   a_m^2 \xi^2( |\hat u_{m}|^2 - |\hat u_{m-1}|^2 )  
+ a_{m-1}a_m \xi \, \Re (i \hat{u}_{m-2} \bar{\hat{u}}_m) =  0.
\end{equation}
Using Young inequality, this yields
\begin{equation}\label{2dissipation-8}
\xi \partial_t E_m 
+ \frac{1}{2}  a_m^2 \xi^2 |\hat u_{m}|^2 
\le a_m^2 \xi^2 |\hat u_{m-1}|^2 
+ \frac{1}{2} a_{m-1}^2  |\hat{u}_{m-2}|^2.
\end{equation}
where $E_m = - a_m \Re (i \hat u_{m-1} \bar{\hat{u}}_m)$.

\medskip

\noindent
{\bf Step 3.}\ \
In this step, we sum up the energy inequalities constructed in the previous step
and then make the desired energy inequality.
The strategy is essentially the same as in Subsection 3.2. 

For this purpose, 
we first multiply \eqref{2dissipation-6} with $\ell = m-2$ and \eqref{2dissipation-8} 
by $\xi^2$ and $\beta_1$, respectively.
Then we combine the resultant equation, obtaining
\begin{multline*}
 \partial_t 
\big\{\xi^2 E_{m-1}  + \beta_1 \xi  E_m \big\} 
+ \frac{1}{2} \beta_1 a_m^2 \xi^2 |\hat u_{m}|^2 
+ \Big( \frac{1}{2} a_{m-1}^2  - \beta_1 a_m^2 \Big) \xi^2 |\hat u_{m-1}|^2 \\
 \le \Big( \frac{1}{2} \beta_1 + \xi^2 \Big) a_{m-1}^2  |\hat u_{m-2}|^2
+ \frac{1}{2} a_{m-1}^2 \xi^4 |\hat{u}_{m-3}|^2 
+ |a_{m-1}| |a_{m}| |\xi|^3 |\hat{u}_{m-2}| |\hat{u}_{m}|.
\end{multline*}
Letting $\beta_1$ suitably small and using Young inequality, we get
\begin{multline*}
 \partial_t \big\{\xi^2 E_{m-1}  + \beta_1 \xi  E_m \big\} 
 + c \xi^2 (|\hat u_{m}|^2 + |\hat u_{m-1}|^2) \\
 \le C (1 + \xi^2)^2  |\hat u_{m-2}|^2
+ \frac{1}{2} a_{m-1}^2 \xi^4 |\hat{u}_{m-3}|^2.
\end{multline*}
Moreover, combining the above estimate and 
\eqref{2dissipation-6} with $\ell = m-2$, we get 
\begin{multline}\label{2eneq-1}
 \partial_t \big\{(1+\xi^2) E_{m-1}  + \beta_1 \xi  E_m \big\} 
 + c \xi^2 |\hat u_{m}|^2  + c (1+\xi^2)|\hat u_{m-1}|^2 \\
 \le C (1 + \xi^2)^2  |\hat u_{m-2}|^2
+ \frac{1}{2} a_{m-1}^2 \xi^2(1+\xi^2) |\hat{u}_{m-3}|^2.
\end{multline}

Second, we multiply \eqref{2eneq-1} and \eqref{2dissipation-7} with $\ell = m-4$
by $\beta_2 \xi^2$ and $(1+\xi^2)^2$
and combine the resultant equations. 
Then  we obtain
\begin{multline*}
 \partial_t \big\{\beta_2 \xi^2((1+\xi^2) E_{m-1}  + \beta_1 \xi  E_m) + \xi (1+\xi^2)^2 E_{m-2} \big\} \\
 + \beta_2 c \xi^4 |\hat u_{m}|^2  + \beta_2 c \xi^2(1+\xi^2)|\hat u_{m-1}|^2  
+ \Big( \frac{1}{2} a_{m-2}^2 - \beta_2 C \Big) \xi^2 (1 + \xi^2)^2  |\hat u_{m-2}|^2 \\
\le C \xi^2(1+\xi^2)^2 |\hat{u}_{m-3}|^2
+ \frac{1}{2} a_{m-3}^2 (1+\xi^2)^2  | \hat{u}_{m-4}|^2  
+ C |\xi| (1+\xi^2)^2 |\hat{u}_{m-3}| |\hat{u}_{m-1}|,
\end{multline*}
Letting $\beta_2$ suitably small and using Young inequality, we get
\begin{multline}\label{2eneq-2}
 \partial_t \big\{\beta_2 \xi^2((1+\xi^2) E_{m-1}  + \beta_1 \xi  E_m) + \xi (1+\xi^2)^2 E_{m-2} \big\} \\
 +  c \xi^4 |\hat u_{m}|^2  +  c \xi^2(1+\xi^2)|\hat u_{m-1}|^2 + c \xi^2 (1 + \xi^2)^2  |\hat u_{m-2}|^2 \\
\le C (1+\xi^2)^3 |\hat{u}_{m-3}|^2
+ \frac{1}{2} a_{m-3}^2 (1+\xi^2)^2  | \hat{u}_{m-4}|^2.  
\end{multline}

Third, we multiply \eqref{2eneq-2}  and \eqref{2dissipation-6} with $\ell = m-4$
by $\beta_3$ and $(1+\xi^2)^3$, respectively,
and combine the resultant equations.
Then we obtain
\begin{equation*}
\begin{split}
& \partial_t \big\{\beta_3(\beta_2 \xi^2((1+\xi^2) E_{m-1}  + \beta_1 \xi  E_m) + \xi (1+\xi^2)^2 E_{m-2}) + (1+\xi^2)^3 E_{m-3} \big\} \\
& + \beta_3 c \xi^4 |\hat u_{m}|^2  + \beta_3 c \xi^2(1+\xi^2)|\hat u_{m-1}|^2 + \beta_3 c \xi^2 (1 + \xi^2)^2  |\hat u_{m-2}|^2 \\
&+  \Big( \frac{1}{2} a_{m-3}^2 -\beta_3 C\Big) (1+\xi^2)^3 |\hat{u}_{m-3}|^2 \\
&\le C (1+\xi^2)^3  | \hat{u}_{m-4}|^2 + \frac{1}{2} a_{m-3}^2 \xi^2  (1+\xi^2)^3 |\hat{u}_{m-5}|^2
+ C |\xi| (1+\xi^2)^3 |\hat{u}_{m-4}| |\hat{u}_{m-2}|.
\end{split}
\end{equation*}
Therefore, letting $\beta_3$ suitably small and using Young inequality, we get
\begin{multline}\label{2eneq-3}
 \partial_t \big\{\beta_3(\beta_2 \xi^2((1+\xi^2) E_{m-1}  + \beta_1 \xi  E_m) + \xi (1+\xi^2)^2 E_{m-2}) \\
 + (1+\xi^2)^3 E_{m-3} \big\} 
 + c \xi^4 |\hat u_{m}|^2  + c \xi^2(1+\xi^2)|\hat u_{m-1}|^2 +  c \xi^2 (1 + \xi^2)^2  |\hat u_{m-2}|^2 \\
 +  c (1+\xi^2)^3 |\hat{u}_{m-3}|^2 
\le C (1+\xi^2)^4  | \hat{u}_{m-4}|^2 + \frac{1}{2} a_{m-3}^2 \xi^2  (1+\xi^2)^3 |\hat{u}_{m-5}|^2.
\end{multline}

Inspired by the derivation of \eqref{2eneq-1}, \eqref{2eneq-2} and \eqref{2eneq-3},
we can conclude that the following inequality
\begin{multline}\label{2eneq-4}
 \partial_t \mathcal{E}_{m-5}
+ c \sum_{\ell = 5}^{m} \xi^{2([\ell/2]-2)}(1+\xi^2)^{m - \ell}| \hat{u}_{\ell}|^2  \\
 \le C(1+\xi^2)^{m-4} | \hat{u}_{4}|^2
+ \frac{1}{2} a_{5}^2 \xi^2(1+\xi^2)^{m-5} |\hat{u}_{3}|^2,
\end{multline}
is derived by the induction argument.
Here
$[ \  ]$ denotes the greatest integer function, 
and $\mathcal{E}_1 = \beta_1 \xi  E_m + (1+\xi^2) E_{m-1}$ and 
\begin{equation}\label{2defE}
\begin{split}
\mathcal{E}_\ell &= \beta_\ell \xi^2 \mathcal{E}_{\ell - 1} + \xi (1+\xi^2)^{\ell} E_{m - \ell}, \\
\mathcal{E}_{\ell+1} &= \beta_{\ell+1} \mathcal{E}_{\ell} + (1+\xi^2)^{\ell + 1} E_{m- (\ell +1)},
\end{split}
\end{equation}
for $\ell$ are even integers with $\ell \ge 2$.

Furthermore, we multiply \eqref{2eneq-4} and \eqref{2dissipation6-5} 
by $\beta_{m-4} \xi^2$ and $(1+\xi^2)^{m-4}$, respectively,
and combine the resultant equation.
Then we obtain
\begin{equation*}
\begin{split}
& \partial_t \mathcal{E}_{m-4}
+ \beta_{m-4} c \sum_{\ell = 5}^{m} \xi^{2([\ell/2]-1)}(1+\xi^2)^{m - \ell}| \hat{u}_{\ell}|^2 \\
&+ \Big( \frac{1}{2}  a_4^2 - \beta_{m-4} C\Big) \xi^2 (1+\xi^2)^{m-4} | \hat{u}_{4}|^2 \\
& \le C \xi^2(1+\xi^2)^{m-4} |\hat{u}_{3}|^2 + \frac{1}{2} (1+\xi^2)^{m-4} |\hat{u}_{2} |^2 
+ C |\xi| (1+\xi^2)^{m-4} |\hat{u}_{3}| |\hat{u}_{5}|,
\end{split}
\end{equation*}
where $\mathcal{E}_{m-4}$ is defined by \eqref{2defE} with $\ell = m-4$.
Thus, letting $\beta_{m-4}$ suitably small and using Young inequality, we obtain
\begin{multline}\label{2eneq-5}
\partial_t \mathcal{E}_{m-4}
+  c \sum_{\ell = 4}^{m} \xi^{2([\ell/2]-1)}(1+\xi^2)^{m - \ell}| \hat{u}_{\ell}|^2 \\
\le C  (1+\xi^2)^{m-3} |\hat{u}_{3}|^2 + \frac{1}{2} (1+\xi^2)^{m-4} |\hat{u}_{2} |^2. 
\end{multline}

Similarly,  we multiply \eqref{2eneq-5} and \eqref{2dissipation6-1} by $\beta_{m-3}$ and $(1+\xi^2)^{m-3}$, 
combine the resultant equalities, and take $\beta_{m-3}$ suitably small.
Then we have
\begin{multline}\label{2eneq-5'}
\partial_t \mathcal{E}_{m-3}
+  c \sum_{\ell = 3}^{m} \xi^{2([\ell/2]-1)}(1+\xi^2)^{m - \ell}| \hat{u}_{\ell}|^2 \\
\le  C (1+\xi^2)^{m-2} |\hat{u}_{2} |^2 + \xi^2 (1+\xi^2)^{m-3} |\hat{u}_{1}|^2,
\end{multline}
where $\mathcal{E}_{m-3}$ is defined by \eqref{2defE} with $\ell = m-3$.

To estimate $|\hat{u}_1|^2$ in \eqref{2eneq-5'}, we next employ \eqref{2estU2}.
Namely, we multiply \eqref{2estU2} and \eqref{2eneq-5'} 
by $(1+\xi^2)^{m-3}$ and $\beta_{m-2} \alpha_m \xi^{m/2-2}$, respectively.
Then we combine the resultant equation, obtaining
\begin{multline*}
 \partial_t \big\{ \beta_{m-2} \alpha_m \xi^{m/2-2} \mathcal{E}_{m-3} + (1+\xi^2)^{m-3} E_1^{(m)} \big\} \\
+  \beta_{m-2} \alpha_m c \xi^{m/2-2}  \sum_{\ell = 3}^{m} \xi^{2([\ell/2]-1)}(1+\xi^2)^{m - \ell}| \hat{u}_{\ell}|^2 
+ \alpha_m (1  - \beta_{m-2})\xi^{m/2} (1+\xi^2)^{m-3} |\hat{u}_{1}|^2 \\
  \le C  \xi^{m/2-2} (1+\xi^2)^{m-2} |\hat{u}_{2} |^2 +  \gamma \alpha_m \xi^{m/2-1} (1+\xi^2)^{m-3} \Re (i \hat{u}_{1} \bar{\hat{u}}_{2}) \\
+ \xi (1+\xi^2)^{m-3}  \, \Re (i^{m/2+1}\mathcal{U}_m \bar{\hat{u}}_{2}) ,
\end{multline*}
where we have defined $\alpha_m = \prod_{j=2}^{m/2} a_{2j}$.
Here, taking $\beta_{m-2}$ suitably small and using Young inequality, we get
\begin{multline}\label{2eneq-6}
 \partial_t \big\{ \beta_{m-2} \alpha_m \xi^{m/2-2} \mathcal{E}_{m-3} + (1+\xi^2)^{m-3} E_1^{(m)} \big\} \\
+ c \xi^{m/2-2}  \sum_{\ell = 3}^{m} \xi^{2([\ell/2]-1)}(1+\xi^2)^{m - \ell}| \hat{u}_{\ell}|^2
+ c \xi^{m/2} (1+\xi^2)^{m-3} |\hat{u}_{1}|^2 \\
  \le C  \xi^{m/2-2} (1+\xi^2)^{m-2} |\hat{u}_{2} |^2  
+ \xi (1+\xi^2)^{m-3}  \, \Re (i^{m/2+1}\mathcal{U}_m \bar{\hat{u}}_{2}) .
\end{multline}
For the last term of the right hand side in \eqref{2eneq-6},
we note that
\begin{multline*}
\mathcal{U}_{m} =  \Big( \prod_{j=0}^{m/2-3} a_{m- 2j}\Big) (-i\xi)^{m/2-2} \hat u_{4} + \Big( \prod_{j=2}^{m/2-1} a_{2j+1}\Big) \hat u_{m} \\
+ \sum_{k=3}^{m/2-1} \Big( \prod_{j=2}^{k-1} a_{2j+1}\Big) \Big( \prod_{j=0}^{m/2-1-k} a_{m- 2j}\Big) (-i\xi)^{m/2-k} \hat u_{2 k},
\end{multline*}
for $m \ge 6$,
where the last term of the right hand side is neglected in the case $m=6$.
Then, substituting the above equality into \eqref{2eneq-6}, we obtain
\begin{multline}\label{2eneq-7}
 \partial_t \big\{ \beta_{m-2} \alpha_m \xi^{m/2-2} \mathcal{E}_{m-3} + (1+\xi^2)^{m-3} E_1^{(m)} \big\} \\
+ c \xi^{m/2-2}  \sum_{\ell = 3}^{m} \xi^{2([\ell/2]-1)}(1+\xi^2)^{m - \ell}| \hat{u}_{\ell}|^2
+ c \xi^{m/2} (1+\xi^2)^{m-3} |\hat{u}_{1}|^2 \\
  \le C  \xi^{m/2-2} (1+\xi^2)^{m-2} |\hat{u}_{2} |^2  
+ C \sum_{k=2}^{m/2}|\xi|^{m/2 + 1-k} (1+\xi^2)^{m-3} |\hat u_{2}|   |\hat u_{2 k}|. 
\end{multline}
In order to control the term of $|\hat{u}_m|$ on the right hand side of \eqref{2eneq-7}
we introduce the following  inequality
\begin{equation*}
|\xi|^{3m/2-5 } (1+\xi^2)^{m-3}|\hat{u}_{2}||\hat{u}_{m}| 
\le \varepsilon \xi^{3m-10}| \hat{u}_{m}|^2 
+ C_\varepsilon (1+\xi^2)^{2(m-3)} |\hat{u}_{2}|^2.
\end{equation*}
Inspired by the above inequality,
we multiply \eqref{2eneq-7} by $\xi^{3m/2-6}$ and employ this inequality.
Then we obtain
\begin{multline*}
\xi^{3m/2-6}  \partial_t \big\{ \beta_{m-2} \alpha_m \xi^{m/2-2} \mathcal{E}_{m-3} +  (1+\xi^2)^{m-3} E_1^{(m)} \big\} 
+ (c - \varepsilon) \xi^{3m-10}| \hat{u}_{m}|^2 \\
+ c \xi^{2 m- 10}  \sum_{\ell = 3}^{m-1} \xi^{2[\ell/2]}(1+\xi^2)^{m - \ell}| \hat{u}_{\ell}|^2  + c \xi^{2m-6} (1+\xi^2)^{m-3} |\hat{u}_{1}|^2\\
  \le \{ C  \xi^{2 m-8}  + C_\varepsilon (1+\xi^2)^{m-3}\} (1+\xi^2)^{m-3} |\hat{u}_{2} |^2 \\
+ C \sum_{k=2}^{m/2-1} |\xi|^{2 m - 5 -k} (1+\xi^2)^{m-3} |\hat u_{2}| |\hat u_{2 k}|. 
\end{multline*}
Therefore, letting $\varepsilon$ suitably small, we have
\begin{multline}\label{2eneq-8}
\xi^{3m/2-6}  \partial_t \big\{ \beta_{m-2} \alpha_m \xi^{m/2-2} \mathcal{E}_{m-3} +  (1+\xi^2)^{m-3} E_1^{(m)} \big\}  \\
+ c \xi^{2 m- 10}  \sum_{\ell = 3}^{m} \xi^{2[\ell/2]}(1+\xi^2)^{m - \ell}| \hat{u}_{\ell}|^2  + c \xi^{2m-6} (1+\xi^2)^{m-3} |\hat{u}_{1}|^2\\
  \le C (1+\xi^2)^{2(m-3)} |\hat{u}_{2} |^2 
+ C \sum_{k=2}^{m/2-1} |\xi|^{2 m - 5 -k} (1+\xi^2)^{m-3} |\hat u_{2}| |\hat u_{2 k}|. 
\end{multline}
Moreover, applying the inequality
\begin{multline*}
 |\xi|^{2m-5-k}  (1+\xi^2)^{m-3} |\hat{u}_{2}||\hat{u}_{2k}| \\
\le   \varepsilon \xi^{2m -10 + 2k}(1+\xi^2)^{m- 2k}| \hat{u}_{2k}|^2 
+ C_\varepsilon \xi^{2m-4k} (1+\xi^2)^{m -6 + 2k} |\hat{u}_{2}|^2
\end{multline*}
to \eqref{2eneq-8}, 
we can get 
\begin{multline}\label{2eneq-9}
 \partial_t \mathcal{E}_{m-2}
+ c \xi^{2 m- 10}  \sum_{\ell = 3}^{m} \xi^{2[\ell/2]}(1+\xi^2)^{m - \ell}| \hat{u}_{\ell}|^2 \\
 + c \xi^{2m-6} (1+\xi^2)^{m-3} |\hat{u}_{1}|^2  \le C (1+\xi^2)^{2(m-3)} |\hat{u}_{2} |^2,
\end{multline}
where we have defined 
$\mathcal{E}_{m-2} = 
\xi^{3m/2-6} ( \beta_{m-2} \alpha_m \xi^{m/2-2} \mathcal{E}_{m-3} +  (1+\xi^2)^{m-3} E_1^{(m)} )$.

Finally, multiplying the basic energy \eqref{2eq6} and  \eqref{2eneq-9} by $(1+\xi^2)^{2(m-3)}$ and  $\beta_{m-1}$, respectively,
combining the resultant equations and letting $\beta_{m-1}$ suitably small,
then this yields
\begin{multline}\label{2eneq-10}
 \partial_t \Big\{ \frac{1}{2}(1+\xi^2)^{2(m-3)}|\hat{u}|^2 + \beta_{m-1} \mathcal{E}_{m-2} \Big\} 
+ c \xi^{2 m- 6}  (1+\xi^2)^{m - 3}| \hat{u}_{1}|^2 \\
+ c (1+\xi^2)^{2(m-3)} |\hat{u}_{2} |^2
+ c \xi^{2 m- 10} \sum_{\ell = 3}^{m} \xi^{2[\ell/2]}(1+\xi^2)^{m - \ell}| \hat{u}_{\ell}|^2  \le 0.
\end{multline}
Thus, integrating the above estimate with respect to $t$,
we obtain the following energy estimate
\begin{multline}\label{2eneq-11}
 |\hat{u}(t,\xi)|^2 
+ \int^t_0   \Big\{
\frac{\xi^{2m-6}}{(1+\xi^2)^{m-3}} |\hat u_{1}|^2 + |\hat u_{2}|^2 \\
+  \frac{\xi^{2 m- 10}}{(1+\xi^2)^{m-3}} \sum_{\ell = 3}^{m} \frac{\xi^{2[\ell/2]}}{(1+\xi^2)^{\ell -3}}| \hat{u}_{\ell}|^2  
  \Big\} d\tau \le  C|\hat{u}(0,\xi)|^2.
\end{multline}
Here we have used the following inequality
\begin{equation}\label{2eneq-12}
c |\hat{u}|^2 \le 
\frac{1}{2} |\hat{u}|^2  +   \frac{\beta_{m-1}}{(1+\xi^2)^{2(m-3)}} \mathcal{E}_{m-2}
\le C |\hat{u}|^2
\end{equation}
for suitably small $\beta_{m-1}$.
Furthermore the estimate \eqref{2eneq-10} with \eqref{2eneq-12} give us the following pointwise estimate
\begin{equation}
\notag
|\hat{u}(t,\xi)|
\le C e^{- c \lambda(\xi)t} |\hat{u}(0,\xi)|,
\qquad
\lambda(\xi) = 
\frac{\xi^{3m-10}}{(1+\xi^2)^{2(m-3)}}.
\end{equation}
This therefore proves \eqref{point2} and completes the proof of Theorem \ref{thm2}.


\subsection{Construction of the matrices $K$ and  $S$}

In this section, inspired by the energy method stated in Sections 3.2 and 3.3, 
we derive the desired matrices $K$ and $S$.

Based on the energy method of Step 2 in Subsection 3.2, 
we first introduce the following $m \times m$ matrices:
\begin{equation*}
K_1 =
\left(
\begin{array}{cccc:c}
 {0} & {1} & {0} & {0}  &    \\
 {-1} & {0} & {0} & {0}  &    \\
 {0} & {0} & {0} & {0}  &    \mbox{\smash{\huge\textit{O}}}  \\
 {0} & {0} & {0} & {0}  &    \\
 \hdashline
     &      &      &      &       \\
       &     \mbox{\smash{\huge\textit{O}}}    &     &     &  \mbox{\smash{\huge\textit{O}}}     \\
 \end{array}
\right),
\qquad
K_4 = a_4
\left(
\begin{array}{cccc:c}
 {0} & {0} & {0} & {0}  &    \\
 {0} & {0} & {0} & {0}  &    \\
 {0} & {0} & {0} & {-1}  &    \mbox{\smash{\huge\textit{O}}}  \\
 {0} & {0} & {1} & {0}  &    \\
 \hdashline
     &      &      &      &       \\
       &     \mbox{\smash{\huge\textit{O}}}    &     &     &  \mbox{\smash{\huge\textit{O}}}     \\
 \end{array}
 \right).
\end{equation*}
Then, we multiply \eqref{Fsys1} by $-i\xi K_1$ and take the inner product with $\hat u$.
Moreover, taking the real part of the resultant equation, we have 
\begin{equation}\label{2FsysK1-1}
-\frac{1}{2} \xi \partial_t
\langle i K_1 \hat u, \hat u \rangle 
 + \xi^2 \langle [K_1A_m]^{\rm sy} \hat u, \hat u \rangle
 -  \xi \langle i [K_1L_m]^{\rm asy} \hat u, \hat u \rangle = 0, 
\end{equation}
where 
\begin{equation*}
K_1A_m =  
\left(
\begin{array}{cccc:cc}
 {1} & {0} & {0} & {0}   & \\
 {0} & {-1} & {0} & {0}   &  \\
 {0} & {0} & {0} & {0}   & \mbox{\smash{\huge\textit{O}}}  \\
 {0} & {0} & {0} & {0}  &  \\
 \hdashline
     &      &      &      &     &  \\
       &   \mbox{\smash{\huge\textit{O}}}  &   &   &   \mbox{\smash{\huge\textit{O}}}      \\
 \end{array}
\right), \qquad
K_1 L_m =  
\left(
\begin{array}{cccc:c}
 {0} & {\gamma} & {1} & {0}  &    \\
 {0} & {0} & {0} & {0}  &     \\
 {0} & {0} & {0} & {0}  &     \mbox{\smash{\huge\textit{O}}}  \\
 {0} & {0} & {0} & {0}  &     \\
 \hdashline
     &      &      &      &   \\
       &     \mbox{\smash{\huge\textit{O}}}    &     &     &  \mbox{\smash{\huge\textit{O}}}      \\
 \end{array}
\right). \\
\end{equation*}
The equality \eqref{2FsysK1-1} is equivalent to \eqref{2dissipation6-4'}. 
Similarly, by using the matrix $K_4$, we can obtain
\begin{equation}\label{2FsysK4-1}
-\frac{1}{2} \xi \partial_t
\langle i K_4 \hat u, \hat u \rangle 
 + \xi^2 \langle [K_4A_m]^{\rm sy} \hat u, \hat u \rangle
 -  \xi \langle i [K_4L_m]^{\rm asy} \hat u, \hat u \rangle = 0, 
\end{equation}
where
\begin{equation*}
\begin{split}
K_4A_m & = a_4^2  
\left(
\begin{array}{cccc:c}
 {0} & {0} & {0} & {0}   & \\
 {0} & {0} & {0} & {0}    &  \\
 {0} & {0} & {-1} & {0}    & \mbox{\smash{\huge\textit{O}}}  \\
 {0} & {0} & {0} & {1}    &  \\
 \hdashline
     &      &            &     &  \\
       &   \mbox{\smash{\huge\textit{O}}}  &   &      & \mbox{\smash{\huge\textit{O}}}      \\
 \end{array}
\right), \qquad
K_4 L_m = - a_4
\left(
\begin{array}{cccc:cc}
 {0} & {0} & {0} & {0}  &   {0}  & \\
 {0} & {0} & {0} & {0 }  & {0} &    \\
 {0} & {0} & {0} & {0}  &    {a_5}  & \mbox{\smash{\huge\textit{O}}}  \\
 {0} & {1} & {0} & {0}  &   {0}  & \\
 \hdashline
     &      &      &      &  &   \\
       &     \mbox{\smash{\huge\textit{O}}}    &     &  &   &  \mbox{\smash{\huge\textit{O}}}      \\
 \end{array}
\right). \\
\end{split}
\end{equation*}
The equality \eqref{2FsysK4-1} is equivalent to \eqref{2dissipation6-5'}. 
%
\medskip

We next introduce 
\begin{equation*}
S_3 =
\left(
\begin{array}{cccc:c}
 {0} & {0} & {0} & {0}  &    \\
 {0} & {0} & {1} & {0}  &    \\
 {0} & {1} & {0} & {0}  &    \mbox{\smash{\huge\textit{O}}}  \\
 {0} & {0} & {0} & {0}  &    \\
 \hdashline
     &      &      &      &       \\
       &     \mbox{\smash{\huge\textit{O}}}    &     &     &  \mbox{\smash{\huge\textit{O}}}     \\
 \end{array}
\right), 
\qquad
\tilde{S}_{\ell} = \hskip -3mm
\bordermatrix*{
    &						  	 &		& 1 		&							& 	&	&	&	\cr
    &    \mbox{\smash{\huge\textit{O}}} &		& 0 		&\mbox{\smash{\huge\textit{O}}} 	&	&	& 	&	\cr
    &  							 &		&\vdots 	&							&	&	&	&	\cr
1  &    0  						&\cdots	& 0 		& \cdots 						&0	&\ell	&	&	\cr
    & \mbox{\smash{\huge\textit{O}}}   &		& \vdots 	&\mbox{\smash{\huge\textit{O}}}	&	&	& 	&	\cr
    &         						&		& 0         	&							&	&	&	&	\cr
   &     						&		&  \ell     	&							&	&	&	&	
}\\[2mm]
\end{equation*}
for $2 \le \ell \le m-1$.
Then, by using the same argument, 
we can show that the equality
\begin{equation}\label{Fsys2S}
\frac{1}{2} \partial_t
  \langle  S_{3} \hat u, \hat u \rangle 
 + \xi \langle i [S_{3} A_m]^{\rm asy} \hat u, \hat u \rangle
 +  \langle [S_{3}L_m]^{\rm sy} \hat u, \hat u \rangle = 0,
\end{equation}
which satisfies
\begin{equation*}
S_3A_m =  
\left(
\begin{array}{cccc:cc}
 {0} & {0} & {0} & {0}   & \\
 {0} & {0} & {0} & {a_4}   &  \\
 {1} & {0} & {0} & {0}   & \mbox{\smash{\huge\textit{O}}}  \\
 {0} & {0} & {0} & {0}  &  \\
 \hdashline
     &      &      &      &     &  \\
       &   \mbox{\smash{\huge\textit{O}}}  &   &   &   \mbox{\smash{\huge\textit{O}}}      \\
 \end{array}
\right), \qquad
S_3 L_m =  
\left(
\begin{array}{cccc:c}
 {0} & {0} & {0} & {0}  &    \\
 {0} & {-1} & {0} & {0}  &     \\
 {0} & {\gamma} & {1} & {0}  &     \mbox{\smash{\huge\textit{O}}}  \\
 {0} & {0} & {0} & {0}  &     \\
 \hdashline
     &      &      &      &   \\
       &     \mbox{\smash{\huge\textit{O}}}    &     &     &  \mbox{\smash{\huge\textit{O}}}      \\
 \end{array}
\right) 
\end{equation*}
is equivalent to \eqref{2dissipation6-1'}.
Similarly, we derive that
\begin{equation}\label{FsysS-ell}
\frac{1}{2} \partial_t
  \langle  \tilde{S}_{2j} \hat u, \hat u \rangle 
 + \xi \langle i [\tilde{S}_{2j} A_m]^{\rm asy} \hat u, \hat u \rangle
 +  \langle [\tilde{S}_{2j}L_m]^{\rm sy} \hat u, \hat u \rangle = 0,
\end{equation}
which satisfies
\begin{equation*}
\tilde{S}_{2j} A_m = \hskip-12mm
\bordermatrix*{
  						  	 &		& a_{2j}	&							& 	&	&	&	\cr
     &	  \mbox{\smash{\huge\textit{O}}}	& 0 		&\mbox{\smash{\huge\textit{O}}} 	&	&	& 	&	\cr
  							 &		&\vdots 	&							&	&	&	&	\cr
0 \ \ 1\ \ 0  						&\cdots  &0 		& \cdots 						&0	&\  2j	&	&	\cr
     &		 \mbox{\smash{\huge\textit{O}}}   & \vdots 	&\mbox{\smash{\huge\textit{O}}}	&	&	& 	&	\cr
            						&		& 0         	&							&	&	&	&	\cr
        						&		&  2j-1     	&							&	&	&	&	
}\\[2mm]
\end{equation*}
and
\begin{equation*}
\tilde{S}_{2j} L_m =  \hskip-3mm
\bordermatrix*{
0 &   					  \cdots	 	&0  \ \ a_{2j+1} \ \ 0			&  \cdots	&	0&	&	&\cr
   &                                                             		& 0 		&	&		& 	&	\cr
&    \mbox{\smash{\huge\textit{O}}}   		& \vdots 	&\mbox{\smash{\huge\textit{O}}}	&	&	& 	&&	\cr
&             								& 0         	&							&	&	&&	&	\cr
&        								&  2j+1     	&							&	&	&&	&	
},\\[2mm]
\end{equation*}
is equivalent to 
\begin{equation}
\notag
	\begin{split}
& \partial_t \Re(\hat u_{1} \bar{\hat{u}}_{2j}) 
 -  a_{2j} \xi \, \Re (i \hat{u}_{1} \bar{\hat{u}}_{2j-1}) 
 +  a_{2j+1} \, \Re ( \hat{u}_{1} \bar{\hat{u}}_{2j+1}) 
 + \xi \, \Re (i \hat{u}_{2} \bar{\hat{u}}_{2j}) =  0,
	\end{split}
\end{equation}
for $2 \le j \le (m-2)/2$.
Therefore, to construct \eqref{2estU}, we sum up \eqref{FsysS-ell} with respect to $j$ with $2 \le j \le (m-2)/2$,
and find that 
\begin{equation}\label{FsysS-ell2}
\frac{1}{2} \partial_t
  \langle  \tilde{\mathcal{S}}_{m-2} \hat u, \hat u \rangle 
 + \xi \langle i [\tilde{\mathcal{S}}_{m-2} A_m]^{\rm asy} \hat u, \hat u \rangle
 +  \langle [\tilde{\mathcal{S}}_{m-2}L_m]^{\rm sy} \hat u, \hat u \rangle = 0
\end{equation}
is equivalent to \eqref{2estU}.
Here we define $\tilde{\mathcal{S}}_{2\ell}$ as $\tilde{\mathcal{S}}_4 =  \tilde{S}_4$ and 
$$
\tilde{\mathcal{S}}_{2 \ell} 
= a_{2 \ell} \xi \tilde{\mathcal{S}}_{2 \ell -2} + \prod^{\ell-1}_{j=2} a_{2j +1} \tilde{S}_{2 \ell} 
$$
for $\ell \ge 3$.
Consequently, multiplying \eqref{2FsysK1-1} by $\prod_{j=2}^{m/2} a_{2j} \xi^{m/2-2} $
and combining the resultant equality and \eqref{FsysS-ell2}, we obtain
\begin{equation}\label{FsysS-ell3}
\begin{split}
&\frac{1}{2} \partial_t
  \Big\langle \Big( \tilde{\mathcal{S}}_{m-2}  - i \prod_{j=2}^{m/2} a_{2j} \xi^{m/2-1} K_1 \Big) \hat u, \hat u \Big\rangle  \\
& +  \Big\langle \Big[\tilde{\mathcal{S}}_{m-2}L_m + \prod_{j=2}^{m/2} a_{2j} \xi^{m/2} K_1A_m\Big]^{\rm sy} \hat u, \hat u \Big\rangle \\
 &+ \xi \langle i [\tilde{\mathcal{S}}_{m-2} A_m]^{\rm asy} \hat u, \hat u \rangle - \prod_{j=2}^{m/2} a_{2j} \xi^{m/2-1} \langle i [K_1L_m]^{\rm asy} \hat u, \hat u \rangle  = 0.
 \end{split}
\end{equation}
%
This equality is the same as \eqref{2estU2}.

\medskip

Based on the energy method of Step 3 in Subsection 3.3, 
we next introduce the following $m \times m$ matrices:
$$
S_{\ell+1} = a_{\ell+1} \hskip-3mm
\bordermatrix*[{(  )}]{
       &           & &0         & 0          &  &  & &  \cr
       & \mbox{\smash{\huge\textit{O}}} 
       & &\vdots &\vdots  &  & \mbox{\smash{\huge\textit{O}}}  & & \cr
       &           & &0         & 0          &  &  & &  \cr
0       & \cdots  &0 &0         & 1          & 0  &\cdots   &0\ \  &  \ell \cr
0       & \cdots  &0 &1        &0           &0  & \cdots   &0\ \ & \ell + 1\cr
       &           & &0          & 0         &     &  &   & \cr
       & \mbox{\smash{\huge\textit{O}}}& &\vdots  & \vdots &    
       & \mbox{\smash{\huge\textit{O}}}    &   &  \cr
       &           & &0         & 0 &   & & &  \cr
 &     &   &\ell  &  \ell+1      & & &   &
}
$$
for $\ell = 4, 6, \cdots , m-2$.
Then, we multiply \eqref{Fsys1} by $S_{\ell+1}$ and take the inner product with $\hat u$.
Furthermore, taking the real part of the resultant equation, we obtain 
\begin{equation}\label{2FsysSL}
\frac{1}{2} \partial_t
  \langle  S_{\ell+1} \hat u, \hat u \rangle 
 + \xi \langle i [S_{\ell+1} A_m]^{\rm asy} \hat u, \hat u \rangle
 +  \langle [S_{\ell+1} L_m]^{\rm sy} \hat u, \hat u \rangle = 0
\end{equation}
for $\ell = 4, 6, \cdots , m-2$, 
where 
$$
\hspace{-6mm}S_{\ell+1} A_m =  a_\ell \hskip-3mm
\bordermatrix*[{(  )}]{
       &           & &0 &0         & 0      &0    &  &  & &  \cr
       & \mbox{\smash{\huge\textit{O}}} 
       & &\vdots &\vdots &\vdots  &\vdots &  & \mbox{\smash{\huge\textit{O}}}  & & \cr
       &           & &0 &0         & 0   &0       &  &  & &  \cr
0       & \cdots  &0 &0 &0   & 0          &a_{\ell+2} &0&\cdots   &0  \ \ &  \ell \cr
0       & \cdots  &0 & a_{\ell} &0         &0        &0  &0& \cdots   &0 \ \ & \ell+1 \cr
       &           & &0 &0          & 0    &0     &     &  &   & \cr
       & \mbox{\smash{\huge\textit{O}}}& &\vdots &\vdots  & \vdots & \vdots &    
       & \mbox{\smash{\huge\textit{O}}}    &   &  \cr
       &           & &0 &0         & 0 &0  &   & & &  \cr
 &     &  &\ell-1  &\ell &  \ell+1   & \ell+2 & &&   &
}
$$
and
$$
S_{\ell+1} L_m= a_{\ell+1}^2 \hskip-3mm
\bordermatrix*[{(  )}]{
       &           & &0         & 0          &  &  & &  \cr
       & \mbox{\smash{\huge\textit{O}}} 
       & &\vdots &\vdots  &  & \mbox{\smash{\huge\textit{O}}}  & & \cr
       &           & &0         & 0          &  &  & &  \cr
0       & \cdots  &0 &-1         & 0          & 0  &\cdots   &0\ \   &  \ell \cr
0       & \cdots  &0 &0        &1           &0  & \cdots   &0 \ \  & \ell + 1\cr
       &           & &0          & 0         &     &  &   & \cr
       & \mbox{\smash{\huge\textit{O}}}& &\vdots  & \vdots &    
       & \mbox{\smash{\huge\textit{O}}}    &   &  \cr
       &           & &0         & 0 &   & & &  \cr
 &     &   &\ell  &  \ell+1      & & &   &
} \qquad\quad .
$$
We note that the equalities \eqref{2FsysSL} is equivalent to \eqref{2dissipation-6'}.

On the other hand,
we introduce the following $m \times m$ matrices:
$$
K_{\ell+2} = a_{\ell + 2} \hskip-3mm
\bordermatrix*[{(  )}]{
       &           & &0         & 0          &  &  & &  \cr
       & \mbox{\smash{\huge\textit{O}}} 
       & &\vdots &\vdots  &  & \mbox{\smash{\huge\textit{O}}}  & & \cr
       &           & &0         & 0          &  &  & &  \cr
0       & \cdots  &0 &0         & -1          & 0  &\cdots   &0\ \  &  \ell+1 \cr
0       & \cdots  &0 &1        &0           &0  & \cdots   &0\ \  & \ell + 2 \cr
       &           & &0          & 0         &     &  &   & \cr
       & \mbox{\smash{\huge\textit{O}}}& &\vdots  & \vdots &    
       & \mbox{\smash{\huge\textit{O}}}    &   &  \cr
       &           & &0         & 0 &   & & &  \cr
 &     &   &\ell + 1 &  \ell +2     & & &   &
}
$$

\vskip2mm

\noindent
for $\ell = 4, 6, \cdots , m-2$.
Then, we multiply \eqref{Fsys1} by $-i K_{\ell+2}$ and take the inner product with $\hat u$.
Furthermore, taking the real part of the resultant equation, we obtain 
\begin{equation}\label{2FsysKL}
-\frac{1}{2} \xi \partial_t
\langle i K_{\ell+2} \hat u, \hat u \rangle 
 + \xi^2 \langle [K_{\ell+2}A_m]^{\rm sy} \hat u, \hat u \rangle
	 -  \xi \langle i[K_{\ell+2}L_m]^{\rm asy} \hat u, \hat u \rangle = 0, 
\end{equation}
for $\ell = 4,6, \cdots , m-4$, 
where 
$$
K_{\ell+2} A_m= a_{\ell+2}^2 \hskip-3mm
\bordermatrix*[{(  )}]{
       &           & &0         & 0          &  &  & &  \cr
       & \mbox{\smash{\huge\textit{O}}} 
       & &\vdots &\vdots  &  & \mbox{\smash{\huge\textit{O}}}  & & \cr
       &           & &0         & 0          &  &  & &  \cr
0       & \cdots  &0 &-1         & 0          & 0  &\cdots   &0\ \   &  \ell +1 \cr
0       & \cdots  &0 &0        &1           &0  & \cdots   &0\ \   & \ell + 2\cr
       &           & &0          & 0         &     &  &   & \cr
       & \mbox{\smash{\huge\textit{O}}}& &\vdots  & \vdots &    
       & \mbox{\smash{\huge\textit{O}}}    &   &  \cr
       &           & &0         & 0 &   & & &  \cr
 &     &   &\ell + 1  &  \ell+2      & & &   &
}
$$
and
$$
K_{\ell+2} L_m = a_{\ell+2} \hskip-3mm
\bordermatrix*[{(  )}]{
       &           & &0 &0         & 0      &0    &  &  & &  \cr
       & \mbox{\smash{\huge\textit{O}}} 
       & &\vdots &\vdots &\vdots  &\vdots &  & \mbox{\smash{\huge\textit{O}}}  & & \cr
       &           & &0 &0         & 0   &0       &  &  & &  \cr
0       & \cdots  &0 &0 &0  & 0          & - a_{\ell+3} &0&\cdots   &0\ \   &  \ell+1 \cr
0       & \cdots  &0 & -a_{\ell+1} &0         &0           &0  &0& \cdots   &0\ \  & \ell+2 \cr
       &           & &0 &0          & 0    &0     &     &  &   & \cr
       & \mbox{\smash{\huge\textit{O}}}& &\vdots &\vdots  & \vdots & \vdots &    
       & \mbox{\smash{\huge\textit{O}}}    &   &  \cr
       &           & &0 &0         & 0 &0  &   & & &  \cr
 &     &  &\ell-1  &\ell &  \ell + 1  & \ell+3 & &&   &
}\qquad\quad.
$$

\vskip2mm

\noindent
Moreover we have
\begin{equation}\label{2FsysKM}
-\frac{1}{2} \xi \partial_t
\langle i K_{m} \hat u, \hat u \rangle 
 + \xi^2 \langle [K_{m} A_m]^{\rm sy} \hat u, \hat u \rangle
 -  \xi \langle i [K_{m} L_m]^{\rm asy} \hat u, \hat u \rangle = 0,
\end{equation}
where
\begin{equation*}
K_{m}A_m = a_m^2
\left(
\begin{array}{cccccc}
       &               &       &  {0}  &  {0}  \\
       &   \mbox{\smash{\huge\textit{O}}}  &       & {\vdots}  &  {\vdots}   \\
       &               &       & {0}                              &  {0}   \\
 {0} & {\cdots} & {0}              & {-1}  &  {0}    \\
 {0} & {\cdots} & {0}   & {0}  &   {1}   \\
 \end{array}
\right), \quad
K_{m}L_m = a_{m-1}a_m
\left(
\begin{array}{cccccc}
       &                     &  &  {0} & {0} & {0}\\
       &   \mbox{\smash{\huge\textit{O}}}         &  &  {\vdots} &{\vdots} & {\vdots} \\
       &                                           &        &  {0}  &{0} & {0}\\
 {0} & {\cdots}                        & {0}  &  {0}  &{0} & {0} \\
 {0} & {\cdots}   & {0}  &   {-1}  &{0} & {0}\\
 \end{array}
\right).
\end{equation*}
The equalities \eqref{2FsysKL} and  \eqref{2FsysKM} 
are equivalent to \eqref{2dissipation-7'} and \eqref{2dissipation-8'}, respectively.



For the rest of this subsection, we construct the desired matrices.
According to the strategy of Step 3 in Subsection 3.2, we first combine \eqref{2FsysSL} and \eqref{2FsysKM}.
More precisely, multiplying \eqref{2FsysSL} with $\ell = m-2$ and  \eqref{2FsysKM} by $(1+\xi^2)$ and $\delta_1$, respectively, 
and combining the resultant equations, we obtain
\begin{multline*}
\frac{1}{2} \partial_t
\big\langle \big\{(1+\xi^2)S_{m-1}  -  \delta_1 i \xi K_m \big\} \hat u, \hat u \big\rangle  \\
 +(1+\xi^2) \langle [S_{m-1}L_m]^{\rm sy} \hat u, \hat u \rangle +  \delta_1 \xi^2 \langle [K_mA_m]^{\rm sy} \hat u, \hat u \rangle \\
  +  \xi (1+\xi^2) \langle i [S_{m-1}A_m]^{\rm asy} \hat u, \hat u \rangle  
  - \delta_1 \xi \langle i [K_1L_m]^{\rm asy} \hat u, \hat u \rangle= 0.
\end{multline*}
We next multiply \eqref{2FsysKL} with $\ell = m-4$ and the above equation by $(1+\xi^2)^2$ and $\delta_2 \xi^2$, respectively, 
and combining the resultant equations, we obtain
\begin{multline*}
\frac{1}{2} \partial_t
\big\langle \big\{\delta_2 \xi^2((1+\xi^2)S_{m-1}  -  \delta_1 i \xi K_m) - i \xi  (1+\xi^2)^2 K_{m-2}\big\} \hat u, \hat u \big\rangle  \\
 + \delta_2 \xi^2(1+\xi^2) \langle [S_{m-1}L_m]^{\rm sy} \hat u, \hat u \rangle + \xi^2 \langle [(\delta_1\delta_2 \xi^2 K_m +  (1+\xi^2)^2 K_{m-2})A_m]^{\rm sy} \hat u, \hat u \rangle \\
  +\delta_2 \xi^3 (1+\xi^2) \langle i [S_{m-1}A_m]^{\rm asy} \hat u, \hat u \rangle  \\
  -  \xi \langle i [(\delta_1\delta_2 \xi^2 K_m +  (1+\xi^2)^2 K_{m-2})L_m]^{\rm asy} \hat u, \hat u \rangle= 0.
\end{multline*}

Furthermore, multiplying \eqref{2FsysSL} with $\ell = m-4$ and the above equation by $(1+\xi^2)^3$ and $\delta_3$, respectively, 
and combining the resultant equations, we get
\begin{multline}\label{2FsysSK2-1}
\frac{1}{2} \partial_t
\big\langle \big\{\delta_3(\delta_2 \xi^2((1+\xi^2)S_{m-1}  -  \delta_1 i \xi K_m) \\
- i \xi  (1+\xi^2)^2 K_{m-2}) + (1+\xi^2)^3 S_{m-3}\big\} \hat u, \hat u \big\rangle  \\
 + (1+\xi^2) \langle [(\delta_2 \delta_3 \xi^2S_{m-1} +  (1+\xi^2)^2S_{m-3} )L_m]^{\rm sy} \hat u, \hat u \rangle \\
+ \delta_3 \xi^2 \langle [(\delta_1\delta_2 \xi^2 K_m +  (1+\xi^2)^2 K_{m-2})A_m]^{\rm sy} \hat u, \hat u \rangle \\
  +\xi (1+\xi^2) \langle i  [(\delta_2 \delta_3 \xi^2S_{m-1} +  (1+\xi^2)^2S_{m-3} )A_m]^{\rm asy}  \hat u, \hat u \rangle  \\
  - \delta_3 \xi \langle i [(\delta_1\delta_2 \xi^2 K_m +  (1+\xi^2)^2 K_{m-2})L_m]^{\rm asy} \hat u, \hat u \rangle= 0.
\end{multline}
Now, we introduce the new matrices $\mathcal{K}_{\ell}$ and $\mathcal{S}_\ell$ as $\mathcal{K}_0 =  K_m$ and 
$$
\mathcal{K}_{\ell} 
= \delta_{\ell-1}\delta_{\ell} \xi^2 \mathcal{K}_{\ell -2} +  (1+\xi^2)^{\ell} K_{m-\ell} 
$$
for $\ell \ge 2$, and 
$\mathcal{S}_1 =  S_{m-1}$ and 
$$
\mathcal{S}_{\ell} 
= \delta_{\ell-1}\delta_{\ell} \xi^2 \mathcal{S}_{\ell -2} +  (1+\xi^2)^{\ell-1} S_{m-\ell} 
$$
for $\ell \ge 3$.
Then the equation \eqref{2FsysSK2-1} is rewritten as
\begin{multline*}
\frac{1}{2} \partial_t
\big\langle \big\{(1+\xi^2)\mathcal{S}_{3}  -  \delta_3 i \xi \mathcal{K}_2 \big\} \hat u, \hat u \big\rangle 
+ (1+\xi^2) \langle [\mathcal{S}_{3}L_m]^{\rm sy} \hat u, \hat u \rangle 
+ \delta_3 \xi^2 \langle [\mathcal{K}_{2}A_m]^{\rm sy} \hat u, \hat u \rangle \\
  +\xi (1+\xi^2) \langle i  [\mathcal{S}_{3} A_m]^{\rm asy}  \hat u, \hat u \rangle  
  - \delta_3 \xi \langle i [\mathcal{K}_{2}L_m]^{\rm asy} \hat u, \hat u \rangle= 0.
\end{multline*}
Consequently, by the induction argument with respect to $\ell$ in \eqref{2FsysSL} and \eqref{2FsysKL}, 
we arrive at
\begin{multline}\label{2FsysSK2-2}
\frac{1}{2} \partial_t
\big\langle \big\{(1+\xi^2)\mathcal{S}_{m-5}  -  \delta_{m-5} i \xi \mathcal{K}_{m-6} \big\} \hat u, \hat u \big\rangle 
+ (1+\xi^2) \langle [\mathcal{S}_{m-5}L_m]^{\rm sy} \hat u, \hat u \rangle \\
+ \delta_{m-5} \xi^2 \langle [\mathcal{K}_{m-6}A_m]^{\rm sy} \hat u, \hat u \rangle 
  +\xi (1+\xi^2) \langle i  [\mathcal{S}_{m-5} A_m]^{\rm asy}  \hat u, \hat u \rangle  \\
  - \delta_{m-5} \xi \langle i [\mathcal{K}_{m-6}L_m]^{\rm asy} \hat u, \hat u \rangle= 0.
\end{multline}
Applying Young inequality to \eqref{2FsysSK2-2}, we can obtain \eqref{2eneq-4}.

Moreover, we multiply \eqref{2FsysK4-1} and \eqref{2FsysSK2-2} by $(1+\xi^2)^{m-4}$ and $\delta_{m-4}\xi^2$, respectively,
and combine the resultant equations. Then this yields
\begin{multline*}
\frac{1}{2} \partial_t
\big\langle \big\{\delta_{m-4}\xi^2(1+\xi^2)\mathcal{S}_{m-5}  -  i \xi \mathcal{K}_{m-4} \big\} \hat u, \hat u \big\rangle \\
+ \delta_{m-4}\xi^2(1+\xi^2) \langle [\mathcal{S}_{m-5}L_m]^{\rm sy} \hat u, \hat u \rangle 
+ \xi^2 \langle [\mathcal{K}_{m-4}A_m]^{\rm sy} \hat u, \hat u \rangle \\
  +\delta_{m-4}\xi^3 (1+\xi^2) \langle i  [\mathcal{S}_{m-5} A_m]^{\rm asy}  \hat u, \hat u \rangle  
  - \xi \langle i [\mathcal{K}_{m-4}L_m]^{\rm asy} \hat u, \hat u \rangle= 0.
\end{multline*}
Similarly, Moreover, we multiply \eqref{Fsys2S} and the above equation by $(1+\xi^2)^{m-3}$ and $\delta_{m-3}$, respectively,
and combine the resultant equations. Then we get
\begin{multline}\label{2FsysSK2-2'}
\frac{1}{2} \partial_t
\big\langle \big\{(1+\xi^2)\mathcal{S}_{m-3}  -  \delta_{m-3} i \xi \mathcal{K}_{m-4} \big\} \hat u, \hat u \big\rangle 
+ (1+\xi^2) \langle [\mathcal{S}_{m-3}L_m]^{\rm sy} \hat u, \hat u \rangle \\ 
+ \delta_{m-3} \xi^2 \langle [\mathcal{K}_{m-4}A_m]^{\rm sy} \hat u, \hat u \rangle 
  + \xi (1+\xi^2) \langle i  [\mathcal{S}_{m-3} A_m]^{\rm asy}  \hat u, \hat u \rangle \\  
  - \delta_{m-3} \xi \langle i [\mathcal{K}_{m-4}L_m]^{\rm asy} \hat u, \hat u \rangle= 0.
\end{multline}
By Young inequality to \eqref{2FsysSK2-2'}, we can derive \eqref{2eneq-5'}.

We next employ \eqref{FsysS-ell3} constructed before.
Multiplying \eqref{FsysS-ell3} and \eqref{2FsysSK2-2'} by $(1+\xi^2)^{m-3}$ and $\delta_{m-2} \alpha_m \xi^{m/2-2}$, respectively,
and combining the resultant equations, we get
\begin{multline}\label{2FsysSK2-3}
\frac{1}{2} \partial_t
\big\langle \big\{(1+\xi^2)\mathcal{S}'  -  \alpha_{m} i \xi^{m/2-1} \mathcal{K}' \big\} \hat u, \hat u \big\rangle 
+ (1+\xi^2) \langle [\mathcal{S}'L_m]^{\rm sy} \hat u, \hat u \rangle \\ 
+ \alpha_{m} \xi^{m/2} \langle [\mathcal{K}' A_m]^{\rm sy} \hat u, \hat u \rangle 
  + \xi (1+\xi^2) \langle i  [\mathcal{S}' A_m]^{\rm asy}  \hat u, \hat u \rangle \\  
  - \alpha_{m} \xi^{m/2-1} \langle i [\mathcal{K}' L_m]^{\rm asy} \hat u, \hat u \rangle= 0,
\end{multline}
where we have defined 
\begin{equation*}
\begin{split}
\mathcal{S}' &= \delta_{m-2} \alpha_m \xi^{m/2-2}\mathcal{S}_{m-3} + (1+\xi^2)^{m-4} \tilde{\mathcal{S}}_{m-2}, \\
\mathcal{K}' &= \delta_{m-2} \delta_{m-3} \mathcal{K}_{m-4} + (1+\xi^2)^{m-3}K_1,
\end{split}
\end{equation*}
and had already defined $\alpha_m = \prod_{j=2}^{m/2} a_{2j}$.
By \eqref{2FsysSK2-3}, we can get \eqref{2eneq-9}.


Finally,  multiplying \eqref{2FsysSK2-3} by $\delta_{m-1}\xi^{3m/2-6}/(1+\xi^2)^{2(m-3)}$, 
and combining \eqref{2eq} and the resultant equations, we can obtain
\begin{multline}\label{2FsysSKfinal-1}
\frac{1}{2} \partial_t
\Big\langle \Big[ I + \frac{\delta_{m-1}}{(1+\xi^2)^{2(m-3)}}  \big\{\xi^{3m/2-6}(1+\xi^2)\mathcal{S}'  -  \alpha_m i \xi^{2m-7}  \mathcal{K}' \big\}\Big] \hat u, \hat u \Big\rangle  \\
 + \langle L_m \hat u, \hat u \rangle 
+ \delta_{m-1} \frac{\xi^{3m/2-6}}{(1+\xi^2)^{2m-7}} \langle [\mathcal{S}'L_m]^{\rm sy} \hat u, \hat u \rangle \\
+ \alpha_m \delta_{m-1} \frac{\xi^{2(m-3)}}{(1+\xi^2)^{2(m-3)}} \langle [\mathcal{K}'A_m]^{\rm sy} \hat u, \hat u \rangle 
  -  \alpha_m \delta_{m-1} \frac{\xi^{2m-7}}{(1+\xi^2)^{2(m-3)}}  \langle i [\mathcal{K}' L_m]^{\rm asy} \hat u, \hat u \rangle  \\
 + \delta_{m-1} \frac{\xi^{3m/2-5}}{(1+\xi^2)^{2m-7}} \langle i [\mathcal{S}'A_m]^{\rm asy} \hat u, \hat u \rangle   = 0,
\end{multline}
where $I$ denotes an identity matrix.
%
Letting $\delta_1, \cdots, \delta_{m-1}$ suitably small,
then \eqref{2FsysSKfinal-1} derives the energy estimate \eqref{2eneq-11}.
To be more precise, we introduce
\begin{equation*}
\mathcal{K}_{m-4} 
 = (1+\xi^2)^{m-4} K_{4} 
 +  \sum_{k=3}^{m/2} \prod_{j=2}^{k-1} \delta_{m-2j}\delta_{m-2j-1} \xi^{2(k-2)}(1+\xi^2)^{m-2k} K_{2k} 
\end{equation*}
for $m \ge 6$, 
and hence
\begin{equation}\label{2estK}
\begin{split}
\mathcal{K}' & =  (1+\xi^2)^{m-3}K_1
+ \delta_{m-2} \delta_{m-3}  (1+\xi^2)^{m-4} K_{4} \\
 &\qquad +  \delta_{m-2} \delta_{m-3} \sum_{k=3}^{m/2} \prod_{j=2}^{k-1} \delta_{m-2j}\delta_{m-2j-1} \xi^{2(k-2)}(1+\xi^2)^{m-2k} K_{2k}.
\end{split}
\end{equation}
Moreover, we find that
\begin{equation*}
\mathcal{S}_{m-3} 
 = (1+\xi^2)^{m-4} S_{3} 
 +  \sum_{k=3}^{m/2} \prod_{j=2}^{k-1} \delta_{m-2j}\delta_{m-2j+1} \xi^{2(k-2)}(1+\xi^2)^{m-2k} S_{2k-1}
\end{equation*}
for $m \ge 6$, and $\tilde{\mathcal{S}}_{4} = \tilde{S}_{4}$, $\tilde{\mathcal{S}}_{6} = a_5 \tilde{S}_{6} + a_6 \xi \tilde{S}_4$ and
\begin{equation*}
\begin{split}
\tilde{\mathcal{S}}_{m-2} 
&=  \prod^{m/2 -2}_{j=2} a_{2j +1} \tilde{S}_{m-2} + \prod^{m/2 -3}_{j=1} a_{m-2j} \xi^{m/2-3} \tilde{S}_{4} \\
&\qquad +  \sum_{k=2}^{m/2-3}\Big( \prod_{j=2}^{m/2-k-1} a_{2j+1} \Big)  \Big( \prod_{j=1}^{k-1} a_{m-2j} \Big) \xi^{k-1} \tilde{S}_{m-2k}
\end{split}
\end{equation*}
for $m \ge 10$, and also
\begin{equation}\label{2estS}
\begin{split}
\mathcal{S}'  
& =  \delta_{m-2} \alpha_m \xi^{m/2-2}(1+\xi^2)^{m-4} S_{3} \\
& +  \alpha_m  \sum_{k=3}^{m/2} \prod_{j=1}^{k-1} \delta_{m-2j}\delta_{m-2j+1} \xi^{m/2 + 2(k-3)}(1+\xi^2)^{m-2k} S_{2k-1} \\ 
&+  \prod^{m/2 -2}_{j=2} a_{2j +1}(1+\xi^2)^{m-4} \tilde{S}_{m-2} + \prod^{m/2 -3}_{j=1} a_{m-2j} \xi^{m/2-3} (1+\xi^2)^{m-4} \tilde{S}_{4} \\
&\qquad +  \sum_{k=2}^{m/2-3}\Big( \prod_{j=2}^{m/2-k-1} a_{2j+1} \Big)  \Big( \prod_{j=1}^{k-1} a_{m-2j} \Big) \xi^{k-1}(1+\xi^2)^{m-4} \tilde{S}_{m-2k}
\end{split}
\end{equation}
Therefore, by using \eqref{2estK} and \eqref{2estS},
we can estimate the dissipation terms as
\begin{multline}\label{2estSK}
 \langle L_m \hat u, \hat u \rangle 
+ \delta_{m-1} \frac{\xi^{3(m-4)/2}}{(1+\xi^2)^{2m-7}} \langle [\mathcal{S}'L_m]^{\rm sy} \hat u, \hat u \rangle \\
+ \delta_{m-1} \frac{\xi^{2(m-3)}}{(1+\xi^2)^{2(m-3)}} \langle [\mathcal{K}'A_m]^{\rm asy} \hat u, \hat u \rangle \\
 \ge c \Big\{
\frac{\xi^{2(m-3)}}{(1+\xi^2)^{m-3}} |\hat u_{1}|^2  +  |\hat u_{2}|^2 
 + \sum_{j=2}^{m/2}  \frac{\xi^{2(m+j-6)}}{(1+\xi^2)^{m+2j-7}} |\hat u_{2j-1}|^2  \\
 + \sum_{j=2}^{m/2}  \frac{\xi^{2(m+j-5)}}{(1+\xi^2)^{m+2j-6}} |\hat u_{2j}|^2  \Big\}, 
\end{multline}
for suitably small $\delta_1, \cdots, \delta_{m-1}$.
We note that this estimate is the same as the dissipation part of \eqref{2eneq-11}.
Consequently we conclude that
our desired symmetric matrix $S$  
and skew-symmetric matrix $K$
are described as
\begin{equation*}
S = \frac{\xi^{3(m-4)/2}}{(1+\xi^2)^{2m-7}} \mathcal{S}',
\qquad
K = \frac{\xi^{2(m-3)}}{(1+\xi^2)^{2(m-3)}}\mathcal{K}'.
\end{equation*}




\section{Alternative approach}


\subsection{General strategy}
In this section, by using the Fourier energy method, we provide an alternative way 
to justify the dissipative structure of the linear symmetric hyperbolic system with relaxation \eqref{sys1}.
%
The key point of the approach is to derive from the above system a new system of $m$ number of equations or inequalities 
\begin{equation*}
(I_1), (I_2), \cdots, (I_j),\cdots, (I_m),
\end{equation*}
in the Fourier space, such that their appropriate linear combination can capture the dissipation rate of all the degenerate components only over the frequency domain far from $|\xi|=0$ and $|\xi|=\infty$. Precisely,  for any $0<\epsilon<M<\infty$, by considering
\begin{equation}
\label{aag.p1}
\sum_{j=1}^mc_j I_j
\end{equation}
for an appropriate choice of constants $c_j>0$ $(1\leq j\leq m)$ which may depend on $\epsilon$ and $M$, we expect to obtain  that for $\epsilon\leq |\xi|\leq M$,
\begin{equation}
\label{aag.p2}
\pa_t \{|\hat u|^2 + \Re E^{int}_1 (\hat u)\} + c_{\epsilon,M} |\hat u|^2 \leq 0,
\end{equation}
where $c_{\epsilon,M}>0$  depending on $\epsilon$ and $M$ is a constant, and $E^{int}_1 (\hat u)$ is an interactive  functional such that $|\hat u|^2 + \Re E^{int}_1 (\hat u)\sim |\hat u|^2$ over $\epsilon\leq |\xi|\leq M$.  To deal with the dissipation rate around $|\xi|=0$ or $|\xi|=\infty$, instead of \eqref{aag.p1},  we re-consider the frequency weighted linear combination in the form of 
\begin{equation}
\label{aag.p3}
\sum_{j=1}^mc_j \frac{|\xi|^{\al_j}}{(1+|\xi|)^{\al_j +\be_j}} I_j.
\end{equation}
Here $\al_j\geq 0$ and $\be_j\geq 0$ $(1\leq j\leq m)$ are constants to be chosen such that the similar computations for deriving \eqref{aag.p2} can be applied so as to obtain a Lyapunov inequality  taking the form
\begin{equation}
\label{aag.p4}
\pa_t \{|\hat u|^2 + \Re E^{int} (\hat u)\} +c\sum_{j=1}^m \la_j (\xi) |\hat u_j|^2\leq 0,
\end{equation}
for all $t\geq 0$ and all $\xi\in \R$, where $c>0$ is a constant, $\la_j (\xi)$ $(j=1,2,\cdots,m)$ are nonnegative rational functions of $|\xi|$, and $E^{int} (\hat u)$ is an interactive functional such that $|\hat u|^2 + \Re E^{int} (\hat u)\sim |\hat u|^2$ for all $\xi \in \R$. If \eqref{aag.p4} was proved then by defining 
\begin{equation*}
\la_{min} (\xi) =\min_{1\leq j\leq m} \la_j (\xi),\quad \xi \in \R,
\end{equation*} 
it follows that 
\begin{equation*}
|\hat u (t,\xi)|^2 \leq C e^{-c\la_{min} (\xi) t}|\hat u (0,\xi)|^2,
\end{equation*}
for all $t\geq 0$ and all $\xi\in \R$, which thus implies the dissipative structure of the considered system \eqref{sys1}. Observe that $\la_j(\xi)$ $(1\leq j\leq m)$ and hence $\la_{min}(\xi)$ may depend on  $\al_j\geq 0$ and $\be_j\geq 0$ $(1\leq j\leq m)$. In general,  $\al_j$ and $\be_j$ are required to satisfy a series of inequalities such that \eqref{aag.p3} indeed can be applied to deduce \eqref{aag.p4} by using the Cauchy-Schwarz inequalities. Therefore we always expect to choose constants $\al_j$ and $\be_j$  such that $\la_{min} (\xi) $ is optimal in the sense that $\la_{min} (\xi) $ may tend to zero when $|\xi|\to 0$ or $|\xi|\to \infty$ in the slowest rate. Finally, we remark that due to \eqref{aag.p2} which holds over $\epsilon\leq |\xi|\leq M$, considering \eqref{aag.p3} is equivalent to considering   
$
\sum_{j=1}^mc_j|\xi|^{\al_j} I_j
$ over $|\xi|\leq \epsilon$ with $0<\epsilon\leq 1$, and 
$
\sum_{j=1}^mc_j |\xi|^{-\be_j} I_j
$ over $|\xi|\geq M$ with $M\geq 1$. In such way, it is more convenient to derive those inequalities satisfied by $\la_j(\xi)$ $(1\leq j\leq m)$.

\subsection{Revisit Model I}

By using the same strategy as in Subsection 2.2 and 2.3,
one can obtain $m$ number of identities $(I_j)$ with $j = 1,2,\cdots,m$ as follows:
\begin{eqnarray*}
&&(I_1):\ \pa_t \lag i\xi \hat{u}_2,\hat{u}_1\rag+|\xi|^2 |\hat{u}_2|^2=- \lag i\xi \hat{u}_2,\hat{u}_4\rag + |\xi|^2 |\hat{u}_1|^2.\\
&&(I_2):\ \dis \pa_t \lag -\hat u_1,\hat u_4\rag +|\hat u_1|^2= |\hat u_4|^2 +\lag i\xi \hat u_2,\hat u_4\rag +\lag \hat u_1, i\xi a_4 \hat u_3 +i\xi a_5 \hat u_5\rag.\\
&&(I_3):\ \dis \pa_t \{\lag i\xi a_4 \hat u_3,\hat u_4\rag -\lag a_4 \hat u_3,\hat u_2\rag\} +a_4^2 |\xi|^2 |\hat u_3|^2 =\\
\dis &&\qquad\qquad\qquad\qquad+a_4^2 |\xi|^2 |\hat u_4|^2 +\lag i\xi a_4 \hat u_3,-i\xi a_5 \hat u_5\rag+a_4^2 \lag i\xi \hat u_4,\hat u_2 \rag.\\
&&(I_4):\  \pa_t \lag i \xi a_5 \hat u_{4}, \hat u_5\rag +a_5^2 |\xi|^2 |\hat u_{4}|^2 =\lag i\xi a_5 \hat u_{4}, -i \xi a_{6} \hat u_{6} \rag \\
&&\dis \qquad\qquad\qquad\qquad+a_5^2 |\xi|^2 |\hat u_5|^2 + a_5 a_{4} |\xi|^2 \lag \hat u_{3}, \hat u_5\rag + \lag i\xi a_5 \hat u_1, \hat u_5\rag.\\
&&(I_{j-1}):\ \pa_t \lag i \xi a_j \hat u_{j-1}, \hat u_j \rag +a_j^2 |\xi|^2 |\hat u_{j-1}|^2 =\lag i\xi a_j \hat u_{j-1}, -i \xi a_{j+1} \hat u_{j+1} \rag \\
&&\dis  \qquad\qquad\qquad\qquad+a_j^2 |\xi|^2 |\hat u_j|^2 + a_j a_{j-1} |\xi|^2 \lag \hat u_{j -2}, \hat u_j \rag,\quad j = 6,7,\cdots,m-1.\\
&&(I_{m-1}):\ \pa_t \lag i \xi a_m \hat u_{m-1},\hat u_m\rag+ a_m^2 |\xi|^2 |\hat u_{m-1}|^2=\lag i\xi a_m \hat u_{m-1}, -\ga \hat u_m\rag\\
&&\dis \qquad\qquad\qquad\qquad+ a_m^2 |\xi|^2 |\hat u_m|^2 +a_{m-1}a_m |\xi|^2 \lag \hat u_{m-2}, \hat u_m\rag.\\
&&(I_m):\ \frac{1}{2} \pa_t |\hat u|^2 +\ga |\hat u_m|^2=0.
\end{eqnarray*}
We note that the equations 
$(I_1), (I_2), (I_3), (I_4), (I_{j-1}), (I_{m-1}), (I_m)$ are parallel to 
\eqref{dissipation6-7'}, \eqref{dissipation6-a}, \eqref{dissipation6-8'}, \eqref{dissipation6-9'}, 
\eqref{dissipation-1'}, \eqref{dissipation-1'}, \eqref{eq}, respectively.
Hence we omit the proof for the derivation of these equations.

\medskip

\noindent{\bf Step 1.} 
We claim that for any $0<\epsilon<M<\infty$, there is $c_{\epsilon,M}>0$ such that for all $\epsilon\leq |\xi|\leq M$,
\begin{equation}
\label{tm1.p2}
\pa_t \{|\hat u|^2 +\Re E^{int}_1(\hat u)\} +c_{\epsilon,M}|\hat u|^2 \leq 0,
\end{equation}
where $E^{int}_1(\hat u)$ is an interactive functional chosen such that 
\begin{equation}
\label{tm1.p3}
|\hat u|^2 +\Re E^{int}_1(\hat u)\sim |\hat u|^2.
\end{equation}

\begin{proof}[Proof of claim:] The key observation is that all the right-hand terms of identities $(I_j)$ $(1\leq j\leq m)$ can be absorbed by the left-hand dissipative terms after taking an appropriate linear combination of all identities.  In fact, let us define
\begin{eqnarray*}
E^{int}_1(\hat u)&=&c_1\lag i\xi \hat{u}_2,\hat{u}_1\rag+c_2\lag -\hat u_1,\hat u_4\rag \\
&&+c_3\{\lag i\xi a_4 \hat u_3,\hat u_4\rag -\lag a_4 \hat u_3,\hat u_2\rag\}+\sum_{j=4}^{m-1}c_j \lag i \xi a_j \hat u_{j-1}, \hat u_j\rag.
\end{eqnarray*}
By taking the real part of each identity $(I_j)$, taking the sum $\sum_{j=1}^mc_jI_j$ with an appropriate choice of constants $c_j$ $(1\leq j\leq m)$, and applying the Cauchy-Schwarz inequality to the right-hand product terms, one can obtain \eqref{tm1.p2}, where constants $c_j$ $(1\leq j\leq m)$ depending on $\epsilon$ and $M$ are chosen such that 
\begin{equation}
\notag
0<c_1\ll c_2\ll \cdots \ll c_{m-2}\ll c_{m-1}\ll 1=c_m.
\end{equation}
The detailed representation of the proof is omitted for brevity. \eqref{tm1.p3} holds true due to $|E^{int}(\hat u)|\leq C_M c_{m-1} |\hat u|^2$ for some constant $C_M$ depending on $M$ and also due to smallness of $c_{m-1}$.
\end{proof}

\noindent{\bf Step 2.}  
Let $|\xi|\geq M$ for $M\geq 1$. We consider the weighted linear combination of identities $(I_j)$ $(1\leq j\leq m)$ in the form of
\begin{equation}
\notag
I_m+\sum_{j=1}^{m-1}c_j |\xi|^{-\be_j} I_j, 
\end{equation}
where $c_j$ $(1\leq j\leq m-1)$ are chosen in terms of step 2, and $\be_j\geq 0$ are chosen such that all the right-hand product terms can be absorbed after using the Cauchy-Schwarz inequality.  In fact, multiplying $(I_j)$ by $|\xi|^{-\be_j}$, one has
\begin{eqnarray*}
&&(I_{\be_1}):\ \pa_t \lag i\xi |\xi|^{-\be_1}\hat{u}_2,\hat{u}_1\rag+|\xi|^{2-\be_1} |\hat{u}_2|^2=- \lag i\xi|\xi|^{-\be_1} \hat{u}_2,\hat{u}_4\rag + |\xi|^{2-\be_1} |\hat{u}_1|^2.\\
&&(I_{\be_2}):\ \dis \pa_t \lag -|\xi|^{-\be_2}\hat u_1,\hat u_4\rag +|\xi|^{-\be_2}|\hat u_1|^2= |\xi|^{-\be_2}|\hat u_4|^2 +\lag i\xi|\xi|^{-\be_2} \hat u_2,\hat u_4\rag \\
&&\qquad\qquad+\lag \hat u_1, i\xi |\xi|^{-\be_2}a_4 \hat u_3 +i\xi |\xi|^{-\be_2}a_5 \hat u_5\rag.\\
&&(I_{\be_3}):\ \dis \pa_t \{\lag i\xi |\xi|^{-\be_3}a_4 \hat u_3,\hat u_4\rag -\lag a_4 |\xi|^{-\be_3}\hat u_3,\hat u_2\rag\} +a_4^2 |\xi|^{2-\be_3} |\hat u_3|^2 =\\
\dis &&\qquad\qquad+a_4^2 |\xi|^{2-\be_3} |\hat u_4|^2 +\lag i\xi|\xi|^{-\be_3}  a_4 \hat u_3,-i\xi  a_5 \hat u_5\rag+a_4^2 \lag i\xi |\xi|^{-\be_3}\hat u_4,\hat u_3\rag.\\
&&(I_{\be_4}):\  \pa_t \lag i \xi|\xi|^{-\be_4} a_5 \hat u_{4}, \hat u_5\rag +a_5^2 |\xi|^{2-\be_4} |\hat u_{4}|^2 =\lag i\xi |\xi|^{-\be_4}a_5 \hat u_{4}, -i \xi a_{6} \hat u_{6} \rag \\
&&\dis \qquad\qquad\qquad\qquad+a_5^2 |\xi|^{2-\be_4} |\hat u_5|^2 + a_5 a_{4} |\xi|^{2-\be_4} \lag \hat u_{3}, \hat u_5\rag + \lag i\xi |\xi|^{-\be_4} a_5 \hat u_1, \hat u_5\rag.\\
&&(I_{\be_{j-1}}):\ \pa_t \lag i \xi |\xi|^{-\be_{j-1}}a_j \hat u_{j-1}, \hat u_j\rag +a_j^2 |\xi|^{2-\be_{j-1}} |\hat u_{j-1}|^2 \\
&&\qquad\qquad=\lag i\xi |\xi|^{-\be_{j-1}}a_j \hat u_{j-1}, -i \xi a_{j+1} \hat u_{j+1} \rag+a_j^2 |\xi|^{2-\be_{j-1}} |\hat u_j|^2 \\
&&\qquad\qquad\qquad\qquad+ a_j a_{j-1} |\xi|^{2-\be_{j-1}} \lag \hat u_{j-2}, \hat u_j\rag,\quad j=6,7,\cdots,m-1.\\
&&(I_{\be_{m-1}}):\ \pa_t \lag i \xi |\xi|^{-\be_{m-1}}a_m \hat u_{m-1},\hat u_m\rag+ a_m^2 |\xi|^{2-\be_{m-1}} |\hat u_{m-1}|^2\\
&&\qquad\qquad=\lag i\xi |\xi|^{-\be_{m-1}}a_m \hat u_{m-1}, -\ga \hat u_m\rag+ a_m^2 |\xi|^{2-\be_{m-1}} |\hat u_m|^2 \\
&&\qquad\qquad\qquad\qquad+a_{m-1}a_m |\xi|^{2-\be_{m-1}} \lag \hat u_{m-2}, \hat u_m\rag.
\end{eqnarray*}
We then require $\be_j$ $(1\leq j\leq m-1)$ to satisfy the following relations. From $(I_{\be_1})$,
\begin{eqnarray*}
&& \be_1-1\geq 0,\quad \be_1-2\geq 0,\\
&& 2(\be_1-1)\geq  (\be_1-2)+(\be_4-2),\quad \be_1-2\geq \be_2,
\end{eqnarray*}
where since $|\xi|\geq M$, $\be_1-1\geq 0$ is such that $\xi |\xi|^{-\be_1}$ in the left first product term of $(I_{\be_1})$ is bounded, $\be_1-2\geq 0$ is such that $|\xi|^{2-\be_1}$  in the left second product term of $(I_{\be_1})$ is bounded, $2(\be_1-1)\geq  (\be_1-2)+(\be_4-2)$ is such that  the product term $\lag i\xi|\xi|^{-\be_1} \hat{u}_2,\hat{u}_4\rag$ on the right first term of $(I_{\be_1})$ can be bounded by the linear combination of the dissipative term $|\xi|^{2-\be_1} |\hat{u}_2|^2$ in $(I_{\be_1})$ and $ |\xi|^{2-\be_4} |\hat u_{4}|^2 $ in $(I_{\be_4})$,  $ \be_1-2\geq \be_2$ is such that the term $|\xi|^{2-\be_1} |\hat{u}_1|^2$ on the right second  term of $(I_{\be_1})$ can be bounded by the dissipative term $|\xi|^{-\be_2}|\hat u_1|^2$ of $(I_{\be_2})$. In terms of the completely same way, from $(I_{\be_j})$ for $j=2,3,\cdots,m-1$, respectively, we require 
\begin{eqnarray*}
&& \be_2\geq 0,\\
&&\be_2\geq \be_4-2,\quad
2(\be_2-1)\geq (\be_1-2)+(\be_4-2),\quad \be_2\geq \be_3,\quad
\be_2\geq \be_5,
\end{eqnarray*}
and
\begin{eqnarray*}
&&\be_3-1\geq 0,\quad \be_3\geq 0,\quad \be_3-2\geq 0,\\
&& \be_3-2\geq \be_4-2,\quad
\be_3\geq \be_5,\quad
2(\be_3-1)\geq (\be_3-2)+(\be_4-2),
\end{eqnarray*}
and
\begin{eqnarray*}
&&\be_4\geq 1,\quad \be_4\geq 2,\quad
\be_4\geq \be_6,\quad 
\be_4\geq \be_5,\\
&&2( \be_4-2)\geq (\be_3-2)+(\be_5-2),\quad 
2( \be_4- 1)\geq \be_2+(\be_5-2),
\end{eqnarray*}
and for 
$j=6,\cdots, m-1$,
\begin{eqnarray*}
&& \be_{j-1}\geq 1,\quad \be_{j-1}\geq 2,\\
&& \be_{j-1}\geq \be_{j+1},\quad 
 \be_{j-1}\geq \be_{j},\quad 
2(\be_{j-1}-2)\geq (\be_{j-2}-2)+(\be_j-2),
\end{eqnarray*}
and
\begin{eqnarray*}
&& \be_{m-1}\geq 1,\quad \be_{m-1}\geq 2,\\
&& 2(\be_{m-1} -1)\geq \be_{m-1}-2,\quad  \be_{m-1} -2 \geq 0,\quad 
2( \be_{m-1} -2) \geq \be_{m-2}-2.
\end{eqnarray*}
Let us choose 
\begin{equation}
\notag
\be_1=4,\quad \be_2=\be_3=\cdots=\be_{m-1}=2,
\end{equation}
which satisfy all the above inequalities of $\be_j$ $(1\leq j\leq m-1)$. 

We now define
\begin{eqnarray*}
E^{int}_\infty(\hat u)&=&c_1\lag i\xi |\xi|^{-4}\hat{u}_2,\hat{u}_1\rag+c_2\lag - |\xi|^{-2}\hat u_1,\hat u_4\rag \\
&&+c_3\{\lag i\xi |\xi|^{-2}a_4 \hat u_3,\hat u_4\rag 
-\lag a_4 |\xi|^{-2}\hat u_3,\hat u_2\rag\}\\
&&+\sum_{j=4}^{m-1}c_j \lag i \xi|\xi|^{-2} a_j \hat u_{j-1}, \hat u_j\rag.
\end{eqnarray*}
Then, as in Step 1, one can show that for any $M\geq 1$, there is $c_{M}>0$ such that for all $|\xi|\geq M$,
\begin{equation}\notag
\pa_t \{|\hat u|^2 +\Re E^{int}_\infty(\hat u)\} +c_{M}\{|\xi|^{-2}(|\hat u_1|^2+|\hat u_2|^2)+\sum_{j=3}^m|\hat u_j|^2\} \leq 0.
\end{equation}

\noindent{\bf Step 3.}  
Let $|\xi|\leq \epsilon$ for $0<\epsilon\leq 1$. 
As in Step 2, we consider the weighted linear combination of identities $(I_j)$ $(1\leq j\leq m)$ in the form of
\begin{equation}
\notag
I_m+\sum_{j=1}^{m-1}c_j |\xi|^{\al_j} I_j, 
\end{equation}
where $c_j$ $(1\leq j\leq m-1)$ are chosen in terms of step 1, and $\al_j\geq 0$ are chosen such that all the right-hand product terms can be absorbed after using the Cauchy-Schwarz inequality.  
In fact, as in Step 2, multiplying $(I_j)$ by $|\xi|^{\al_j}$, one has
\begin{eqnarray*}
&&(I_{\al_1}):\ \pa_t \lag i\xi |\xi|^{\al_1}\hat{u}_2,\hat{u}_1\rag+|\xi|^{2+\al_1} |\hat{u}_2|^2=- \lag i\xi|\xi|^{\al_1} \hat{u}_2,\hat{u}_4\rag + |\xi|^{2+\al_1} |\hat{u}_1|^2.\\
&&(I_{\al_2}):\ \dis \pa_t \lag -|\xi|^{\al_2}\hat u_1,\hat u_4\rag +|\xi|^{\al_2}|\hat u_1|^2= |\xi|^{\al_2}|\hat u_4|^2 +\lag i\xi|\xi|^{\al_2} \hat u_2,\hat u_4\rag \\
&&\qquad\qquad+\lag \hat u_1, i\xi |\xi|^{\al_2}a_4 \hat u_3 +i\xi |\xi|^{\al_2}a_5 \hat u_5\rag.\\
&&(I_{\al_3}):\ \dis \pa_t \{\lag i\xi |\xi|^{\al_3}a_4 \hat u_3,\hat u_4\rag -\lag a_4 |\xi|^{\al_3}\hat u_3,\hat u_2\rag\} +a_4^2 |\xi|^{2+\al_3} |\hat u_3|^2 =\\
\dis &&\qquad\qquad+a_4^2 |\xi|^{2+\al_3} |\hat u_4|^2 +\lag i\xi|\xi|^{\al_3}  a_4 \hat u_3,-i\xi  a_5 \hat u_5\rag+a_4^2 \lag i\xi |\xi|^{\al_3}\hat u_4,\hat u_3\rag.\\
&&(I_{\al_4}):\  \pa_t \lag i \xi|\xi|^{\al_4} a_5 \hat u_{4}, \hat u_5\rag +a_5^2 |\xi|^{2+\al_4} |\hat u_{4}|^2 =\lag i\xi |\xi|^{\al_4}a_5 \hat u_{4}, -i \xi a_{6} \hat u_{6} \rag \\
&&\dis \qquad\qquad\qquad\qquad+a_5^2 |\xi|^{2+\al_4} |\hat u_5|^2 + a_5 a_{4} |\xi|^{2+\al_4} \lag \hat u_{3}, \hat u_5\rag + \lag i\xi |\xi|^{\al_4} a_5 \hat u_1, \hat u_5\rag.\\
&&(I_{\al_{j-1}}):\ \pa_t \lag i \xi |\xi|^{\al_{j-1}}a_j \hat u_{j-1}, \hat u_j\rag +a_j^2 |\xi|^{2+\al_{j-1}} |\hat u_{j-1}|^2 \\
&&\qquad\qquad=\lag i\xi |\xi|^{\al_{j-1}}a_j \hat u_{j-1}, -i \xi a_{j+1} \hat u_{j+1} \rag+a_j^2 |\xi|^{2+\al_{j-1}} |\hat u_j|^2 \\
&&\qquad\qquad\qquad\qquad+ a_j a_{j-1} |\xi|^{2+\al_{j-1}} \lag \hat u_{j-2}, \hat u_j\rag,\quad j=6,7,\cdots,m-1.\\
&&(I_{\al_{m-1}}):\ \pa_t \lag i \xi |\xi|^{\al_{m-1}}a_m \hat u_{m-1},\hat u_m\rag+ a_m^2 |\xi|^{2+\al_{m-1}} |\hat u_{m-1}|^2\\
&&\qquad\qquad=\lag i\xi |\xi|^{\al_{m-1}}a_m \hat u_{m-1}, -\ga \hat u_m\rag+ a_m^2 |\xi|^{2+\al_{m-1}} |\hat u_m|^2 \\
&&\qquad\qquad\qquad\qquad+a_{m-1}a_m |\xi|^{2+\al_{m-1}} \lag \hat u_{m-2}, \hat u_m\rag.
\end{eqnarray*}
As in the case of the large frequency domain, for $|\xi|\leq \epsilon$ with $\epsilon>0$, in order for all the right product terms to be bounded, from equations $(I_{\al_j})$ $(j=1,2,\cdots,m-1)$ above, respectively, we have to require 
\begin{eqnarray*}
&& \al_1+1\geq 0,\ 
2(\al_1+1)\geq (\al_1+2)+(\al_4+2), \ \al_1+2\geq \al_2,
\end{eqnarray*}
and
\begin{eqnarray*}
&& \al_2 \geq \al_4 +2,\
2(\al_2+1) \geq (\al_1+2)+(\al_4+2),\ 
 \al_2\geq \al_3,\
\al_2\geq \al_5,
\end{eqnarray*}
and
\begin{eqnarray*}
&&\al_3 \geq \al_4,\ 
\al_3\geq \al_5,\
2( \al_3+1) \geq (\al_4+2)+(\al_1+2),
\end{eqnarray*}
and 
\begin{eqnarray*}
&& \al_4\geq \al_6,\
\al_4\geq \al_5,\\
&&2( \al_4 +2) \geq (\al_3+2)+(\al_5+2),\ 
2( \al_4 +1) \geq \al_2+(\al_5+2),
\end{eqnarray*}
and for $j=6,\cdots, m-1$,
\begin{eqnarray*}
&& \al_{j-1}\geq \al_{j+1},\
 \al_{j-1}\geq \al_{j},\
2(\al_{j-1}+2) \geq (\al_{j-2}+2)+(\al_j+2),
\end{eqnarray*}
and
\begin{eqnarray*}
&& \al_{m-1}\geq 0,\
\al_{m-1}+2\geq 0,\
2(\al_{m-1} +2) \geq \al_{m-2}+2.
\end{eqnarray*}
To consider the best choice of $\{\al_j\}_{j=1}^{m-1}$, one can see
\begin{equation}
\notag
\al_1\geq \al_2\geq\cdots \geq \al_j\geq \al_{j+1}\geq \cdots\geq \al_{m-2}\geq \al_{m-1}\geq 0:=\al_m,
\end{equation}
with
\begin{eqnarray*}
&& \al_1-\al_4\geq 2,\\
&&\al_2-\al_4\geq 2,\\
&&\al_3-\al_4\geq 2,\\
&& \al_{j-1}-\al_{j}\leq \al_{j}-\al_{j+1},\quad  4\leq j\leq m-1.
\end{eqnarray*}
Therefore, the possible best choice satisfies
\begin{eqnarray*}
&& \al_1-\al_4= 2,\\
&&\al_2-\al_4= 2,\\
&&\al_3-\al_4= 2,\\
&&2=\al_3-\al_4\leq \al_4-\al_5\leq  \cdots\leq \al_{m-1}-\al_m=\al_{m-1}=2,
\end{eqnarray*}
which implies
\begin{eqnarray*}
&& \al_1=\al_2=\al_3=2(m-4),\\
&&\al_{j}=2(m-j-1), \quad 4\leq j\leq m-1.
\end{eqnarray*}

We now define
\begin{eqnarray*}
E^{int}_0(\hat u)&=&c_1\lag i\xi |\xi|^{2(m-4)}\hat{u}_2,\hat{u}_1\rag+c_2\lag - |\xi|^{2(m-4)}\hat u_1,\hat u_4\rag \\
&&+c_3\{\lag i\xi |\xi|^{2(m-4)}a_4 \hat u_3,\hat u_4\rag 
-\lag a_4 |\xi|^{2(m-4)}\hat u_3,\hat u_2\rag\}\\
&&+\sum_{j=4}^{m-1}c_j \lag i \xi|\xi|^{2(m-j-1)} a_j \hat u_{j-1}, \hat u_j\rag.
\end{eqnarray*}
Then, as in Step 1, one can show that for any $0<\epsilon\leq 1$, there is $c_{\epsilon}>0$ such that for all $|\xi|\leq \epsilon$,
\begin{equation*}
\pa_t \{|\hat u|^2 +\Re E^{int}_0(\hat u)\} 
+c_{\epsilon}\{|\xi|^{2m-8}|\hat u_1|^2+|\xi|^{2m-6}|\hat u_2|^2+\sum_{j=3}^m |\xi|^{2(m-j)}|\hat u_j|^2\} \leq 0,
\end{equation*}
which further implies that for $|\xi|\leq \epsilon$,
\begin{equation}
\notag
\pa_t \{|\hat u|^2 +\Re E^{int}_0(\hat u)\} +c_\epsilon |\xi|^{2m-6}|\hat u|^2\leq 0.
\end{equation}

\noindent{\bf Step 4.} For $\xi\in \R$ let us define
\begin{eqnarray*}
E^{int}(\hat u)&=&c_1\frac{|\xi|^{2(m-4)}}{(1+|\xi|)^{2m-4}}\lag i\xi \hat{u}_2,\hat{u}_1\rag+c_2\frac{|\xi|^{2(m-4)}}{(1+|\xi|)^{2m-6}}\lag - \hat u_1,\hat u_4\rag \\
&&+c_3\frac{|\xi|^{2(m-4)}}{(1+|\xi|)^{2m-6}}\{\lag i\xi  a_4 \hat u_3,\hat u_4\rag 
-\lag a_4 \hat u_3,\hat u_2\rag\}\\
&&+\sum_{j=4}^{m-1}c_j \frac{|\xi|^{2(m-j-1)}}{(1+|\xi|)^{2(m-j)}} \lag i \xi a_j \hat u_{j-1}, \hat u_j\rag.
\end{eqnarray*}
As in Step 2 and Step 3, we consider the weighted linear combination of identities $(I_j)$ $(1\leq j\leq m)$ in the form of
\begin{multline*}
I_m+c_1\frac{|\xi|^{2(m-4)}}{(1+|\xi|)^{2m-4}}I_1
+c_2\frac{|\xi|^{2(m-4)}}{(1+|\xi|)^{2m-6}}I_2\\
+c_3\frac{|\xi|^{2(m-4)}}{(1+|\xi|)^{2m-6}}I_3
+\sum_{j=4}^{m-1}c_j \frac{|\xi|^{2(m-j-1)}}{(1+|\xi|)^{2(m-j)}} I_j,
\end{multline*}
where $c_j$ $(1\leq j\leq m-1)$ are chosen in terms of Step 1. 
Thanks to computations in Step 1, Step 2, and Step 3, in the completely same way, one can deduce that for $\xi\in \R$,
 \begin{multline*}
\pa_t \{|\hat u|^2 +\Re E^{int}(\hat u)\} 
+c\{\frac{|\xi|^{2m-8}}{(1+|\xi|)^{2m-6}}|\hat u_1|^2\\
+\frac{|\xi|^{2m-6}}{(1+|\xi|)^{2m-4}}|\hat u_2|^2+\sum_{j=3}^m \frac{|\xi|^{2(m-j)}}{(1+|\xi|)^{2(m-j)}}|\hat u_j|^2\} \leq 0,
\end{multline*}
which further gives
\begin{equation}
\notag
\pa_t \{|\hat u|^2 +\Re E^{int}(\hat u)\} +c\frac{|\xi|^{2m-6}}{(1+|\xi|)^{2m-4}} |\hat u|^2\leq 0.
\end{equation}
Noticing 
$|\hat u|^2 +\Re E^{int}(\hat u)\sim |\hat u|^2$,
it follows that
\begin{equation}
\notag
|\hat u(t,\xi)|\leq C e^{-c\eta (\xi) t}|\hat u (0,\xi)|,\quad \eta (\xi)=\frac{|\xi|^{2m-6}}{(1+|\xi|)^{2m-4}},
\end{equation}
for all $t\geq 0$ and all $\xi\in \R$. Notice that the result here is consistent with \eqref{point1} proved in Section 2.3.

\subsection{Revisit Model II}

In this section we revisit the Model II \eqref{sys1} with coefficients matrices $A_m$ and $L_m$ defined in \eqref{2Tmat}.  For simplicity of representation, we rewrite $A_{m}$ with $m=2n$ as
\begin{equation}\notag
    A_{2n}=
    \left(
      \begin{array}{ccccccc}
        0 & a_{12} &   &   && & \\
        a_{21} &  0 &  & & &&  \\
          &  &\ 0 &\ a_{34} & && \\
          &   &\ a_{43} &\ 0 & & & \\
                 &   &  &  &\ddots & &                  \\
          &&&  && 0& a_{2n-1,2n}\\
          &&&&&a_{2n,2n-1} & 0
      \end{array}
    \right),
\end{equation}
with 
$a_{2j-1,2j}=a_{2j,2j-1}=a_{j}$ for $1\leq j\leq n$,
and also choose $L_m$ with $m=2n$ as 
\begin{equation}\notag
    L_{2n}=\left(
        \begin{array}{ccccccccc}
          0 &  &  &  &  &&&&\\
          & 1 & 1 &  &  &&&& \\
           & -1 & 0 & &  &&&&  \\
           &  &  & 0 & 1 &&&& \\
           &&&-1&0&&&& \\
            &&&&&\ddots&&& \\
            &&&&&&0&1&\\
             &&&&&&-1&0&\\
           &&&&&& &&0
        \end{array}
      \right).
\end{equation}
With notations above, system \eqref{sys1} can read
\begin{eqnarray*}
&&\pa_t \hat u_{2j-1}+i\xi a_j \hat u_{2j} - \hat u_{2j-2}=0,
\\
&&\pa_t \hat u_{2j} + i\xi a_j \hat u_{2j-1} + \hat u_{2j+1} +\delta_{2,2j}\hat u_2=0,\quad j=1,2,\cdots,n, 
\end{eqnarray*}
with the convention that $ \hat u_{2n+1}\equiv 0$ and $\hat u_0\equiv 0$.
As for the model I, we can obtain the following estimates
\begin{equation}\label{tm2.pe1}
    \frac{1}{2}\pa_t |\hat u|^2+ |\hat u_2|^2=0,
\end{equation}
and
\begin{multline}
 \pa_t \Re\lag i\xi a_1\hat u_1, \sum_{j=1}^n(-i\xi)^{1-j} (\prod_{k=2}^j a_k)^{-1} \hat u_{2j}\rag 
+c a_1^2 \xi^2|\hat u_1|^2\\
\lesssim (1+|\xi|)^2|\hat u_2|^2 
+\Re\lag \xi^2a_1^2 \hat u_2, \sum_{j=2}^n(-i\xi)^{1-j} (\prod_{k=2}^j a_k)^{-1} \hat u_{2j}\rag,
\label{tm2.pe2a}
\end{multline}
and
\begin{equation}
\pa_t \Re\lag \hat u_{2j-1},u_{2j-2}\rag +c |\hat u_{2j-1}|^2 
\lesssim |\hat u_{2j-2}|^2+\xi^2 |\hat u_{2j-3}|^2 +\Re\lag -i\xi a_j \hat u_{2j}, \hat u_{2j-2}\rag,
\label{tm2.pej1}
\end{equation}
and
\begin{multline}
\pa_t\Re \lag i\xi a_j\hat u_{2j},\hat u_{2j-1}\rag +c a_j^2 \xi^2|\hat u_{2j}|^2 \\
\lesssim |\hat u_{2j-2}|^2+a_j^2\xi^2 |\hat u_{2j-1}|^2 +\Re\lag -i\xi a_j \hat u_{2j+1}, \hat u_{2j-1}\rag,
\label{tm2.pej2}
\end{multline}
for $j=2,3,\cdots,n$. 
Indeed, by using the equations \eqref{2eq}, \eqref{2estU2}, \eqref{2dissipation-6}, \eqref{2dissipation-7} derived in Subsection 3.3, we can get
\eqref{tm2.pe1}, \eqref{tm2.pe2a}, \eqref{tm2.pej1}, \eqref{tm2.pej2}, immediately.

Let us denote \eqref{tm2.pe1}, \eqref{tm2.pe2a}, \eqref{tm2.pej1}, \eqref{tm2.pej2} by $(I_1)$, $(I_2)$, $(I_{2j-1})$ and $(I_{2j})$, respectively, where $j=2,3,\cdots,n$. Consider the linear combination  of all $2n$ number of equations
\begin{equation}
\notag
\sum_{j=1}^n (c_{2j-1} I_{2j-1} +c_{2j} I_{2j}),
\end{equation}
where $c_1=1$, and $c_k>0$ $(k=2,3,\cdots,2n)$ are constants to be properly chosen. It is straightforward to verify that for any $0<\epsilon<M<\infty$, one can choose 
 constants $c_k$ $(1\leq k\leq 2n)$ depending on $\epsilon$ and $M$, with
\begin{equation}
\notag
0<c_{2n}\ll c_{2n-1}\ll \cdots \ll c_{2j}\ll c_{2j-1}\ll\cdots \ll c_3\ll c_2\ll 1=c_1,
\end{equation}  
such that there is $c_{\epsilon,M}>0$ such that for all $\epsilon\leq |\xi|\leq M$,
\begin{equation}
\notag
\pa_t \{|\hat u|^2 +\Re E^{int}_1(\hat u)\} +c_{\epsilon,M}|\hat u|^2 \leq 0,
\end{equation}
where $E^{int}_1(\hat u)$ is an interactive functional given by
\begin{eqnarray}
E^{int}_1(\hat u) &=&c_2 \lag i\xi a_1\hat u_1, \sum_{j=1}^n(-i\xi)^{1-j} (\prod_{k=2}^j a_k)^{-1} \hat u_{2j}\rag\notag \\
&&+\sum_{j=2}^{n} \{c_{2j-1}\lag \hat u_{2j-1},u_{2j-2}\rag+c_{2j}\lag i\xi a_j\hat u_{2j},\hat u_{2j-1}\rag\},\notag
\end{eqnarray}
satisfying
\begin{equation}
\notag
|\hat u|^2 +\Re E^{int}_1(\hat u)\sim |\hat u|^2,\quad\text{for}\ \epsilon\leq |\xi|\leq M.
\end{equation}
Furthermore, we can consider the frequency weighted linear combination in the form of 
\begin{equation}
\label{tm2.pf1}
\sum_{j=1}^n \{c_{2j-1}\frac{|\xi|^{\al_{2j-1}}}{(1+|\xi|)^{\al_{2j-1}+\be_{2j-1}}} I_{2j-1} +c_{2j} \frac{|\xi|^{\al_{2j}}}{(1+|\xi|)^{\al_{2j}+\be_{2j}}}I_{2j}\},
\end{equation}
where $\al_1=\be_1=0$.
As for the model I, we use the same strategy to determine the choice of  constants 
\begin{equation}\notag
  \al_2,\al_3,\cdots,\al_{2n};\quad \be_2,\be_3,\cdots,\be_{2n}.
\end{equation}
In fact, by considering the low frequency domain $|\xi|\leq \epsilon$ with $\epsilon\leq 1$, $ \al_2,\al_3,\cdots,\al_{2n}$ are required to satisfy inequalities
\begin{eqnarray*}
&& 2-j +\al_2\geq 0, j=2,3,\cdots,n,\\
&& \al_2\geq 0,\\
&&2+\al_2\geq 0,\\
&& \al_3\geq 0, 2+\al_3\geq 2+\al_2,\\
&&\al_4\geq 0,2+\al_4\geq \al_3,\\
&&\al_{2j}\geq 2+\al_{2j-2},2+\al_{2j}\geq \al_{2j-1},\\
&&\al_{2j-1}\geq 2+\al_{2j-2},2+\al_{2j-1}\geq \al_{2j-3},\
j=3,4,\cdots,n,
\end{eqnarray*}
and
\begin{eqnarray*}
&& 2(3-j+\al_2)\geq \al_{2j}+2, j=2,\cdots,n,\\
&& 1+\al_3 \geq \frac{2+\al_4}{2},\\
&& \al_{2j}\geq \frac{1}{2}(\al_{2j+1} +\al_{2j-1})-1,\\
&&\al_{2j+1}\geq \frac{1}{2} (\al_{2j+2}+\al_{2j})-1, j=2,\cdots,n-1.
\end{eqnarray*}
One can take the best choice 
\begin{eqnarray*}
&\dis\al_2=4(n-2),\\
&\dis \al_{2j-1}=\al_{2j}=4(n-2)+2(j-2),\ j=2,3,\cdots,n.
\end{eqnarray*}
Similarly, by considering the high frequency domain $|\xi|\geq M$ with $M\geq 1$,  constants $\be_2,\be_3,\cdots,\be_{2n}$ are required to satisfy inequalities
\begin{eqnarray*}
&& \be_2-2\geq 0,\\
&& \be_3\geq 0,\be_3-2\geq \be_2-2,\\
&&\be_4\geq 0,\be_4-2\geq \be_3,\\
&& \be_{2j}\geq \be_{2j-2}, \be_{2j}-2\geq \be_{2j-1},\\
&& \be_{2j-1} \geq \be_{2j-2}-2,\be_{2j-1}-2\geq \be_{2j-3},\ j=3,\cdots,n,
\end{eqnarray*}
and
\begin{eqnarray*}
&& 2(\be_3-1)\geq \be_4-2,\\
&&\be_2+j-2\geq 0, 2(\be_2+j-3)\geq \be_{2j}-2,j=2,\cdots,n,\\
&& 2(\be_{2j}-1)\geq \be_{2j+1}+\be_{2j-1},\\
&&2 (\be_{2j+1}-1) \geq (\be_{2j+2}-2)+(\be_{2j}-2),j=2,\cdots,n-1.
\end{eqnarray*}
One can take the best choice 
\begin{equation}\notag
\be_{2j}=\be_{2j+1}=2j,\quad j=1,2,\cdots,n.
\end{equation}

Now, by \eqref{tm2.pf1}, let us define the interactive functional
\begin{eqnarray*}
E^{int}(\hat u) &=& c_2 \frac{|\xi|^{\al_2}}{(1+|\xi|)^{\al_2+\be_2}}
\lag i\xi a_1\hat u_1, \sum_{j=1}^n(-i\xi)^{1-j} (\prod_{k=2}^j a_k)^{-1} \hat u_{2j}\rag\\
&&+\sum_{j=2}^n\{c_{2j-1}\frac{|\xi|^{\al_{2j-1}}}{(1+|\xi|)^{\al_{2j-1}+\be_{2j-1}}}\lag \hat u_{2j-1},\hat u_{2j-2}\rag\\
&&\qquad\quad+c_{2j} \frac{|\xi|^{\al_{2j}}}{(1+|\xi|)^{\al_{2j}+\be_{2j}}}\lag i\xi a_j\hat u_{2j},\hat u_{2j-1}\rag\},
\end{eqnarray*}
that is,
\begin{eqnarray*}
E^{int}(\hat u) &=& c_2 \frac{|\xi|^{4n-8}}{(1+|\xi|)^{4n-6}}
\lag i\xi a_1\hat u_1, \sum_{j=1}^n(-i\xi)^{1-j} (\prod_{k=2}^j a_k)^{-1} \hat u_{2j}\rag\\
&&+\sum_{j=2}^n\{c_{2j-1}\frac{|\xi|^{4n+2j-12}}{(1+|\xi|)^{4n+4j-14}}\lag \hat u_{2j-1},\hat u_{2j-2}\rag\\
&&\qquad\quad+c_{2j} \frac{|\xi|^{4n+2j-12}}{(1+|\xi|)^{4n+4j-12}}\lag i\xi a_j\hat u_{2j},\hat u_{2j-1}\rag\},
\end{eqnarray*}
and also define the energy dissipation rate 
\begin{eqnarray*}
D(\hat u) &=& |\hat u_2|^2 +\frac{|\xi|^{2+\al_2}}{(1+|\xi|)^{\al_2+\be_2}} |\hat u_1|^2 \\
&&+\sum_{j=2}^n\{\frac{|\xi|^{\al_{2j-1}}}{(1+|\xi|)^{\al_{2j-1}+\be_{2j-1}}}|\hat u_{2j-1}|^2+\frac{|\xi|^{2+\al_{2j}}}{(1+|\xi|)^{\al_{2j}+\be_{2j}}} |\hat u_{2j}|^2\},
\end{eqnarray*}
that is,
\begin{eqnarray*}
D(\hat u) &=& |\hat u_2|^2 +\frac{|\xi|^{4n-6}}{(1+|\xi|)^{4n-6}} |\hat u_1|^2 \\
&&+\sum_{j=2}^n\{\frac{|\xi|^{4n+2j-12}}{(1+|\xi|)^{4n+4j-14}}|\hat u_{2j-1}|^2+\frac{|\xi|^{4n+2j-10}}{(1+|\xi|)^{4n+4j-12}} |\hat u_{2j}|^2\}.
\end{eqnarray*}
Then it follows that 
\begin{equation}
\notag
\pa_t \{|\hat u|^2 +\Re E^{int}(\hat u) \} +c D (\hat u)\leq 0,
\end{equation}
for all $t\geq 0$ and all $\xi\in \R$. Noticing
\begin{equation}
\notag
|\hat u|^2 +\Re E^{int}(\hat u)  \sim |\hat u|^2,
\end{equation}
and 
\begin{equation}
\notag
D (\hat u) \gtrsim \frac{|\xi|^{6n-10}}{(1+|\xi|)^{8n-12}}|\hat u|^2,
\end{equation}
one can see that the Model II \eqref{sys1}, where coefficients matrices $A_m$ and $L_m$ are defined in \eqref{2Tmat}  with $m=2n$, enjoys the dissipative structure 
\begin{equation}
\notag
|\hat u(t,\xi)|^2 \leq C e^{-c\eta (\xi) t} |\hat u(0,\xi)|,
\end{equation}
with 
\begin{equation}
\notag
\eta (\xi)= \frac{|\xi|^{6n-10}}{(1+|\xi|)^{8n-12}}=\frac{|\xi|^{3m-10}}{(1+|\xi|)^{4m-12}}.
\end{equation}
Hence the derived result here is consistent with \eqref{point2} proved in Theorem \ref{thm2}.

\bigskip

\noindent {\sc Acknowledgments:}\ \ 
The first author is partially supported by
Grant-in-Aid for Young Scientists (B) No.\,25800078
from Japan Society for the Promotion of Science.
The second author's research was supported by the General
Research Fund (Project No.~400912) from RGC of Hong Kong.
The third author is partially supported by
Grant-in-Aid for Scientific Research (A) No.\,22244009.



\end{document}